\documentclass{amsproc}
\usepackage{amsmath}
\usepackage{amssymb}

\title{Hecke algebras of classical type and 
their representation type}
\author{Susumu Ariki}
\address{RIMS, Kyoto university, Kyoto 606--8502, Japan}

\email{ariki@kurims.kyoto-u.ac.jp}

\date{}
\begin{document}
\begin{abstract}
The purpose of this article is to determine 
representation type for all of the Hecke algebras of classical 
type. To do this, we combine methods from our previous 
work, which is used to obtain information on their Gabriel quivers, and 
recent advances in the theory of finite dimensional algebras. 
Principal computation is for Hecke algebras of type $B$ with 
two parameters. Then, we show that the representation type of 
Hecke algebras is governed by their Poincar\'e polynomials. 
\end{abstract}

\maketitle

\newtheorem{thm}{Theorem}
\newtheorem{cor}[thm]{Corollary}
\newtheorem{prop}[thm]{Proposition}
\newtheorem{defn}[thm]{Definition}
\newtheorem{rem}[thm]{Remark}
\newtheorem{examp}[thm]{Example}
\newtheorem{lem}[thm]{Lemma}

\newcommand{\Prod}{\Pi}
\newcommand{\End}{\operatorname{End}}
\newcommand{\Hom}{\operatorname{Hom}}
\newcommand{\Ext}{\operatorname{Ext}}
\renewcommand{\Im}{\operatorname{Im}}
\newcommand{\Ker}{\operatorname{Ker}}
\newcommand{\Rad}{\operatorname{Rad}}
\newcommand{\Soc}{\operatorname{Soc}}
\newcommand{\Top}{\operatorname{Top}}
\newcommand{\Beta}{{\bf\beta}}
\renewcommand{\H}{{\mathcal H}}

\section{Introduction}
\subsection{}
Let $F$ be an algebraically closed field, $A$ an 
$F$--algebra which is finite dimensional as an $F$--vector space. 
$A-mod$ is the category of left $A$--modules which are 
finite dimensional over $F$. We say that 
$A$ is {\sf finite} if the number of the 
isomorphism classes of indecomposable $A$--modules is finite. 

Let $F[X]$ be the polynomial ring generated by the indeterminate $X$. 
If, $A$ is not finite and, for each number $d\in \mathbb N$, 
there are finitely 
many $(A, F[X])$--bimodules $M_1,\dots, M_{n_d}$ which are free of 
finite rank as right $F[X]$--modules such that 
all but a finite number of the isomorphism classes of indecomposable 
$A$--modules of dimension $d$ contain an $A$--module of the form 
$M_i\otimes_{F[X]}F[X]/(X-\lambda)$, for some $i$ and some 
$\lambda\in F$, then we say that $A$ is {\sf tame}. 

Let $F\langle X, Y\rangle$ be the free $F$--algebra generated by 
two indeterminates $X$ and $Y$. We say that $A$ is {\sf wild} 
if there is a $(A, F\langle X, Y\rangle)$--bimodule $M$ which 
is free of finite rank as a right $F\langle X, Y\rangle$--module 
such that the associated functor 
\[
\mathcal F_M=M\otimes_{F\langle X, Y\rangle}-\;:\;
F\langle X, Y\rangle-mod \longrightarrow A-mod
\]
respects indecomposability and isomorphism classes. 
A famous theorem of Drozd \cite[Theorem 2]{Dr1} (see also \cite{C1}) 
asserts that $A$ is finite, tame or wild, and that 
these are mutually exclusive. This is the 
{\sf representation type} of the algebra $A$. 

\subsection{}
Let $q\in F^\times$ and let 
$W$ be a finite Weyl group of classical type. First 
we consider the case where $W$ is an irreducible Weyl group. 
Then, for each of $W(A_{n-1})$, $W(B_n)$ and $W(D_n)$, 
we have the associated Hecke algebra. 
We denote them by $\H^A_n(q)$, $\H^B_n(q)$ and $\H^D_n(q)$ 
respectively. For type $B_n$, we can choose two parameters 
$q, Q\in F^\times$ and the associated Hecke algebra is 
denoted by $\H_n(q,Q)$. 

In 1992, Uno gave a criterion 
for $\H_n^A(q)$ to be finite \cite{U}, and conjectured that 
the criterion would be true for other Hecke algebras. 
This conjecture was settled affirmatively for Hecke algebras of 
classical type \cite{A4}. Crucial for proving the Uno conjecture 
was the result proven in \cite{AM2} and \cite{AM3} 
which tells when $\H_n(q,Q)$ is finite. 

The purpose of this paper is to determine representation type 
for all of the Hecke algebras of classical type. 
Main theorems are Theorem \ref{separated parameter case}, 
Theorem \ref{two parameter case} 
and Theorem \ref{one parameter case}. These combined gives 
Theorem \ref{final result}, which is for the general case. 

To give these final results, 
we begin by analysing the representation type of $\H_n(q,Q)$. 
To carry out this, 
we need plenty of results from the theory of finite dimensional algebras. 
These include standard techniques from the covering theory, 
the theory of special biserial algebras, 
and recent results such as a criterion for wildness using the complexity 
of modules 
by Rickard \cite{Ric} and classification of representation types of two--point 
algebras by Han \cite{Ha2}. As we get $\H_n^B(q)$ by setting $Q=q$ 
and an embedding of $\H_n^D(q)$ into $\H_n(q,1)$ by setting $Q=1$, 
which allows 
us to apply the Clifford theory, we can determine representation type 
for $\H_n^B(q)$ and $\H_n^D(q)$ also. In 
the Clifford theory argument used for $\H_n^D(q)$ we also need a recent 
result of Hu \cite{Hu2}. Finally, the general case follows from 
these results for indivisual Hecke algebras $\H_n^X(q)$ with 
$X=A$, $B$ or $D$. Note that 
the case of $\H_n^A(q)$ was already known by \cite{EN}. 

We remark that a comlete set of the isomorphism classes of 
simple modules, for each of the Hecke algebras of classical 
type, was already given. See \cite{DJ1} for type $A$, 
\cite{DJ2}, \cite{AM1} and \cite{A3} for type $B$, 
and \cite{P}, \cite{Hu1} and \cite{Hu2} for type $D$. 
In type $D$, we assume that the characteristic of $F$ is odd. 
Hence, it is natural to proceed further to the classification of 
the isomorphism classes of the 
indecomposable modules of these algebras. 
Our results show when this is achievable, as 
it is well--known that the wild case is the pathological case. 
See \cite{Ge2} for a different approach for classifying 
simple $\H_n^D(q)$--modules. 

The paper is structured as follows. In section 2 we review various 
results on finite dimensional algebras. Then we treat the case 
of $\H_n(q,Q)$ in section 3. The cases of $\H_n^A(q)$, $\H_n^B(q)$ and 
$\H_n^D(q)$ are handled in section 4. The final section is 
for the general case. In the appendix, I list corrections for \cite{A1}. 

The author is grateful to the people of Beijing Normal University, 
professors Changchang Xi and Bangming Deng in particular. 
His visit to a workshop there, which was held in September 2002, was fruitful. 
He is also grateful to East China Normal University for the hospitality. 
Professor H.~Asashiba explained him several materials 
on finite dimensional algebras, which was very helpful. 

Finally, I add a few words about the references. 
I assume that the reader is familar 
with the crystal basis theory and the canonical basis. So, I list few 
about them in the references. 
On the other hand, as it is not appropriate to assume 
that researchers in our field know well about the theory of finite dimensional 
algebras, I always cite references whenever I quote a result. 
In fact, most of the results we use are not available in book form and 
they are scattered around in a vast literature. 
If the reader is familiar with these results, I recommend starting 
with section 3 and return to lemmas in section 2 when they are used. 

\section{Preliminaries}

\subsection{}
Let $A$ be a finite dimensional $F$--algebra as before, and let 
$\{P_1,\dots, P_s\}$ be a complete set of the isomorphism classes 
of indecomposable projective $A$--modules. 
Then $A$ is Morita--equivalent to 
$\End_A(P_1\oplus\cdots\oplus P_s)^{opp}$, which can be written 
in the form $FQ/I$ where $Q=(Q_0,Q_1)$ is a directed graph 
with nodes $Q_0$ and directed edges $Q_1$, 
and $I$ is an admissible ideal of $FQ$. See \cite{ARS} 
for this basic result of Gabriel. 
The directed graph $Q$ is called 
the {\sf Gabriel quiver} of $A$. 
We identify $Q_0$ with the set of the isomorphism classes 
of simple $A$--modules. 

If $A$ is a symmetric algebra then 
$Q_1$ is described as follows. Let $S$ and $T$ be two simple $A$--modules, 
$P(S)$ and $P(T)$ their projective covers. Then 
we write $a_{ST}$ arrows from $S$ to $T$ where 
$a_{ST}=[\,\Rad P(S)/\Rad^2 P(S):T]$. Thus the Gabriel quiver is the 
directed graph whose adjacency matrix is $(a_{ST})_{S,T\in Q_0}$. 

In \cite{AM2}, we used the following lemma to show that the 
$F$--algebras $\H_n(q,Q)$ with $n=e$ or $n=2f+4$, where 
$0\le f\le\frac{e}{2}$ and $Q=-q^f$, are not finite. 

\begin{lem}
\label{basic local}
Let $A$ be a finite dimensional local basic $F$--algebra. 
Then $A$ is finite if and only if 
$A$ is isomorphic to a truncated polynomial ring $F[X]/(X^N)$ for some 
positive integer $N$. 
\end{lem}

\subsection{}
We say that $A$ is {\sf weakly tame} if 
$A$ is not finite and, for each number $d\in \mathbb N$, 
there are finitely 
many $(A, F[X])$--bimodules $M_1,\dots, M_{n_d}$ which are free of 
finite rank as right $F[X]$--modules such that 
all but a finite number of the isomorphism classes of indecomposable 
$A$--modules of dimension $d$ contain an $A$--module which is 
a direct summand of an $A$--module of the form 
$M_i\otimes_{F[X]}F[X]/(X-\lambda)$, for some $i$ and some 
$\lambda\in F$. 

The following is a result of de la Pe$\tilde{\rm n}$a 
\cite[Chap.I, Proposition 2.3]{dlP1}. 

\begin{prop}
Let $A$ be a finite dimensional $F$--algebra. Then 
$A$ is weakly tame if and only if it is tame. 
\end{prop}

Using this result, Erdmann and Nakano proved the following. 

\begin{prop}[{\cite[Proposition 2.3]{EN}}]
\label{reduction to critical rank}
Let $A$ and $B$ be finite dimensional $F$--algebras and suppose that 
there are functors 
\[
\mathcal F:\;A-mod \longrightarrow B-mod, \;\;
\mathcal G:\;B-mod \longrightarrow A-mod
\]
and a constant $C$ such that, for any $A$--module $M$, 
\begin{itemize}
\item[(1)]
$M$ is a direct summand of $\mathcal G\mathcal F(M)$ 
as an $A$--module, 
\item[(2)]
$\operatorname{dim}_F\mathcal F(M)\le C\operatorname{dim}_F M$. 
\end{itemize}
If $A$ is wild then so is $B$. 
\end{prop}

\begin{cor}
\label{critical rank}
\begin{itemize}
\item[(1)]
If $\H_n(q,Q)$ is wild then so is $\H_m(q,Q)$, for all $m\ge n$. 
\item[(2)]
Let $X$ be one of $A$, $B$, $D$. 
If $\H_n^X(q)$ is wild then so is $\H_m^X(q)$, for all $m\ge n$. 
\end{itemize}
\end{cor}
\begin{proof}
We prove (1). The proof of (2) is the same. 
Let $A=\H_n(q,Q)$ and $B=\H_m(q,Q)$ with $m\ge n$. Then, by taking 
$\mathcal F$ and $\mathcal G$ to be the induction and the restriction 
functors respectively, we can apply Proposition 
\ref{reduction to critical rank}. 
The assumption (1) is satisfied as we have the Mackey decomposition theorem for Hecke algebras, 
and the assumption (2) is obvious. 
\end{proof}

\subsection{}
We need finer notions for wild algebras. 
We say that $A$ is {\sf strictly wild} if 
there exists a fully faithful exact functor 
\[
\mathcal F:\;F\langle X, Y\rangle-mod \longrightarrow A-mod. 
\]
All functors in this paper are assumed to be 
$F$--functors. 

The strictly wildness is equivalent to 
the condition that we have a fully faithful exact functor 
\[
B-mod \longrightarrow A-mod, 
\]
for any finite dimensional $F$--algebra $B$. 

In fact, by \cite[Theorem 3]{Br} or \cite[Proposition 14.10]{S}, 
we know that $B-mod$, for any finite dimensional $F$--algebra $B$, 
can be realized as a full subcategory of 
$F\langle X, Y\rangle-mod$. Hence, strictly wildness implies this. 
To show the converse, we consider the functor for 
a particular $11$--dimensional algebra $\Lambda_6(F)$ which 
is the algebra of $6\times 6$ matrices with non--zero entries 
in the first column and the diagonal. Thus, 
we assume that there is a fully faithful exact functor 
\[
\Lambda_6(F)-mod \longrightarrow A-mod.
\]
For this algebra, we have a fully faithful exact functor 
\[
F\langle X, Y\rangle-mod \longrightarrow \Lambda_6(F)-mod
\]
by \cite[Theorem 2]{Br} or \cite[14.2 Example 11]{S}. 
Hence, the composition of the two 
gives us the desired fully faithful exact functor which shows that 
$A$ is strictly wild. 

Another notion we need is the notion of controlled wildness. 
An $F$--algebra $A$ is {\sf controlled wild} if 
there exist a faithful exact functor $\mathcal F$ and 
a full subcategory $\mathcal C$ of $A-mod$ which is closed under 
direct sums and direct summands such that, if we denote by 
$\Hom_A(-,-)_{\mathcal C}$ the 
morphisms of $A-mod$ which factor through $\mathcal C$, then 
for any $F\langle X, Y\rangle$--modules $M$ and $N$ we have 
\[
\Hom_A(\mathcal F(M), \mathcal F(N))=
\mathcal F(\Hom_{F\langle X, Y\rangle}(M, N))
\oplus \Hom_A(\mathcal F(M), \mathcal F(N))_{\mathcal C}
\]
and $\Hom_A(\mathcal F(M), \mathcal F(N))_{\mathcal C}
\subset\operatorname{Rad}\Hom_A(\mathcal F(M), \mathcal F(N))$. 

Note that strictly wildness implies controlled wildness. 
We also know that controlled wildness implies wildness 
\cite[Proposition 2.2]{Ha1}. This depends on the fact that 
flat $F\langle X, Y\rangle$--modules are free: 
\cite[Theorem 2.2.4, Corollary 2.3.2]{C} imply that 
$F\langle X, Y\rangle$ is 
a left and right free ideal ring \cite[Corollary 2.4.3]{C}. 
Hence, $F\langle X, Y\rangle$ is a semifir 
by \cite[Theorem 1.4.1, Corollary 1.4.2]{C}. Then 
\cite[Proposition 1.4.5]{C} says that 
an $F\langle X, Y\rangle$--module is flat if and only if 
every finitely generated $F\langle X, Y\rangle$--submodule 
is free. 

Using the results from \cite{C}, we can also prove the following. 

\begin{prop}
\label{strictly wild}
A finite dimensional $F$--algebra $A$ is strictly wild if and only if 
there exists an $(A,F\langle X, Y\rangle)$--bimodule $M$ which is free 
of finite rank as a right $F\langle X, Y\rangle$--module such that 
the associated functor $\mathcal F_M$ is fully faithful exact. 
\end{prop}
\begin{proof}
If part is obvious. So, we assume that $A$ is strictly wild. 
We consider the composition of the two fully faithful exact functors
\begin{gather*}
\Lambda_6(F)-mod \longrightarrow A-mod, \quad\text{and} \\[5pt]
F\langle X, Y\rangle-mod \longrightarrow \Lambda_6(F)-mod. 
\end{gather*}

Write $B$ for $\Lambda_6(F)$. 
Let $L$ be the image of the $(B,B)$--bimodule $B$ under the first functor. 
Then the functor is of the form $\mathcal F_L$ and 
$L$ is flat as a right $B$--module. As $B$ is right perfect, 
this implies that $L$ is projective as a right $B$--module by a 
theorem of Bass; see \cite[Theorem 24.25]{La}. Further, 
the proof of \cite[Theorem 2]{Br} shows that the second functor 
is of the form $\mathcal F_N$ where $N$ is free of finite rank 
as a right $F\langle X, Y\rangle$--module. Hence, the composition of these 
functors is of the form $\mathcal F_M$ such that, if we view $M$ as 
a right $F\langle X, Y\rangle$--module then $M$ is a direct summand of 
a free $F\langle X, Y\rangle$--module of finite rank. 
As $F\langle X, Y\rangle$ is a semifir, \cite[Theorem 1.4.1]{C} 
implies that $M$ is 
free of finite rank as a right $F\langle X, Y\rangle$--module. 
\end{proof}

\subsection{}
In this subsection, we review a criterion to tell the 
representation type of a path algebra. The statement (1) of the 
following theorem is due to Gabriel \cite{Ga1} 
and the statement (2) is by \cite{DF} and \cite{N}. 
See \cite[VIII.5.4, 5.5]{ARS} for the finiteness result 
and \cite[Theorem 14.15]{S} for the remaining cases. 
Another proof is explained in \cite[Proposition 4.7.1]{B2}, which is 
most recommendable. 

\begin{thm}
\label{path algebra case}
Let $A=FQ$ be a finite dimensional path algebra. Then 
\begin{itemize}
\item[(1)]
$A$ is finite if and only if the underlying graph of the quiver $Q$ is one of 
the Dynkin diagrams of finite type. 
\item[(2)]
$A$ is tame if and only if the underlying graph of the quiver $Q$ is one of 
the Dynkin diagrams of affine type. 
\item[(3)]
$A$ is wild otherwise. 
\end{itemize}
\end{thm}

In (1) and (2) we do not mean that all of the Dynkin diagrams of finite and 
affine types actually occur. By the assumption that $F$ is 
algebraically closed, only $A_n$, $D_n$, $E_n$ occur in (1), and only 
$A_n^{(1)}$, $D_n^{(1)}$, $E_n^{(1)}$ occur in (2). 

The following is well--known. 

\begin{prop}
\label{wild path algebra}
A path algebra is wild if and only if 
it is strictly wild. 
\end{prop}

See \cite[Theorem 1.6]{K2} for the proof which 
does not depend on \cite{DF} and \cite{N}. We say that 
a path algebra $A=FQ$ is {\sf minimal wild} if it is wild and 
$A/AeA$ is tame or finite, for any non--zero idempotent $e\in A$. 
The list of minimal wild path algebras is given in \cite{K1}. 
Note that we follow the older definition of minimal wildness here. 

\begin{rem}
Recall that each path algebra $A=FQ$ is associated with 
a quadratic form
\[
q_A(x)=\sum_{i\in Q_0}x_i^2-\sum_{(i\rightarrow j)\in Q_1}x_ix_j,
\]
which is the {\sf Tits form} of $A$. As is well--known, 
this coincides with the quadratic form associated with the Ringel--Euler form, 
\[
\langle X, Y\rangle=\operatorname{dim}\Hom_A(X,Y)
-\operatorname{dim}\Ext_A^1(X,Y)
\]
for $A$--modules $X$ and $Y$. 
We say that $q_A(x)$ is {\sf weakly nonnegative} if $q_A(x)\ge0$ for 
all $x$ with $x_i\ge0$ $(i\in Q_0)$. 
For those path algebras which we will meet below, we use Theorem 
\ref{path algebra case} to prove that they are wild. 
However, we also know that 
if $A$ is tame then $q_A(x)$ is weakly nonnegative. 
This is a special case of a more general theorem. See 
\cite[Proposition 1.3]{dlP2} for example. Thus, we can apply 
this result instead. 

In fact, the proof of Theorem \ref{path algebra case} given in 
\cite[Proposition 4.7.1]{B2} uses this and results from \cite{BGP}. 
\end{rem}

\subsection{}
Another important class of finite dimensional algebras 
is the class of tilted algebras. 
See \cite{HR}, \cite{Bo1} and \cite{K3}. 

Let $B$ be a tilted algebra. The Tits form of $B$ coincides with the 
quadratic form associated with the Ringel--Euler form again. 
As the global dimension of $B$ is 
smaller than or equal to $2$ by \cite[1.7, Corollary 1]{Bo1} 
or \cite[Theorem 5.2]{HR}, the Ringel--Euler form here is of the form 
\[
\langle X, Y\rangle=\operatorname{dim}\Hom_B(X,Y)
-\operatorname{dim}\Ext_B^1(X,Y)
+\operatorname{dim}\Ext_B^2(X,Y),
\]
for $B$--modules $X$ and $Y$. 
In \cite[Theorem 6.2]{K3}, it is proved that the Tits form of $B$ is 
weakly nonnegative if and only if $B$ is finite or tame, and that 
the Tits form is not weakly nonnegative if and only if $B$ is 
strictly wild. Hence, 

\begin{thm}
\label{wild tilted algebra}
A tilted algebra is wild if and only if it is strictly wild. 
\end{thm}

We will use a particular wild tilted algebra. 

\begin{defn}
Let $A$ be a finite dimensional path algebra, $M$ a tilting $A$--module. 
The tilted algebra $B=\End_A(M)$ is called a 
{\sf concealed algebra} if $M$ has preprojective direct summands 
(of the Auslander--Reiten quiver of $A$) only. 
\end{defn}

\begin{rem}
We can replace \lq\lq preprojective\rq\rq 
with \lq\lq preinjective\rq\rq in the definition above. 
See \cite[4.3(1)]{R2}.
\end{rem}

In \cite{Un}, there is a list of the concealed algebras 
associated with minimal wild path algebras. 
Among them, we need the following. See \cite[p.150]{Un} or 
XXVIII in the list of \cite[1.5, Theorem 2]{R1}. 

\begin{lem}
\label{concealed algebra case}
Let $Q$ be the directed graph

\setlength{\unitlength}{16pt}
\begin{picture}(15,4)(-2,0)
\put(1.55,2){\vector(1,0){1.4}}
\put(4.45,2){\vector(-1,0){1.4}}
\put(4.5,2){\vector(1,0){1.4}}
\put(6.1,2.5){$\alpha$}
\put(7,3){\vector(-1,-1){1}}
\put(7.5,2.5){$\gamma$}
\put(7,3){\vector(1,-1){1}}
\put(6.1,1.1){$\beta$}
\put(7.5,1.1){$\delta$}
\put(6,2){\vector(1,-1){1}}
\put(8,2){\vector(-1,-1){1}}
\put(8.05,2){\vector(1,0){1.4}}
\put(10.9,2){\vector(-1,0){1.4}}
\end{picture}

\noindent
and let $A=FQ/I$ be the $F$--algebra defined by the 
relation $\beta\alpha=\delta\gamma$.
Then $A$ is a wild concealed algebra. In particular, $A$ is strictly wild. 
\end{lem}

We have explained Lemma \ref{concealed algebra case} in the 
most natural framework, but it may be proved by a more direct method. 
Following the argument 
in \cite[2.3]{R1}, we consider the following quiver. 

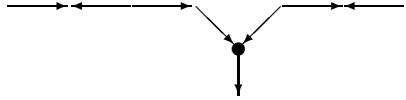
\setlength{\unitlength}{16pt}
\begin{picture}(15,3.5)(-2,-0.6)
\put(1.55,2){\vector(1,0){1.4}}
\put(4.45,2){\vector(-1,0){1.4}}
\put(4.5,2){\vector(1,0){1.4}}
\put(6,2){\vector(1,-1){0.9}}
\put(8,2){\vector(-1,-1){0.9}}
\put(7,0.9){\vector(0,-1){1}}
\put(7,1){\circle*{0.3}}
\put(8.05,2){\vector(1,0){1.4}}
\put(10.9,2){\vector(-1,0){1.4}}
\end{picture}

Denote this quiver by $T$. 
Given a representation of the quiver $Q$ of 
Lemma \ref{concealed algebra case}, 
we associate a representation of $T$ by placing 
the pushout of the linear maps of the arrows $\alpha$ and $\gamma$ 
on the black node, which is located in the middle of $T$. 
The linear maps of the three arrows around the black node are 
the obvious ones. Then, the full subcategory of $A-mod$ 
which is generated by those indecomposable $A$--modules which correspond 
to sincere $FT$--modules under the correspondence is equivalent to 
$(FT-mod)_s$. Thus, by constructing a fully faithful exact functor 
\[
F\langle X, Y\rangle-mod \longrightarrow (FT-mod)_s
\]
in a concrete manner, we know that $A$ is strictly wild. 

\begin{rem}
Concealment is defined in a different way in \cite{R1}, but 
it does not matter as long as the proof of Lemma \ref{concealed algebra case} 
is concerned. 
\end{rem}

\subsection{}
We turn to the theory of complexity. 
Let $A$ be a finite dimensional $F$--algebra, $M$ an $A$--module. 
Let 
\[
P^\bullet:\;\;
\cdots\longrightarrow P_i\longrightarrow\cdots\cdots
\longrightarrow P_1\longrightarrow P_0\longrightarrow M\longrightarrow 0
\]
be the minimal projective resolution of $M$. Then, 
the {\sf complexity} $c_A(M)$ of the $A$--module $M$ is defined as follows. 
\[
c_A(M)=\operatorname{min}
\{s\in\mathbb N\,|\,\exists C,\;\text{a constant, s.t.}\;
\operatorname{dim}_FP_i\le C(i+1)^{s-1},\;\text{for}\;\forall i.\}
\]

In the theorem below, 
the statement (1) is well--known and easy to see. 
(2) is almost obvious. (3) is a rather new result proven by Rickard 
\cite[Theorem 2]{Ric}. 

\begin{thm}
\label{complexity}
Let $A$ be a finite dimensional self--injective $F$--algebra. 
\begin{itemize}
\item[(1)]
An indecomposable $A$--module $M$ has complexity $c_A(M)=0$ if and 
only if $M$ is a projective $A$--module. In particular, 
$A$ is a semisimple algebra if and only if all indecomposable 
$A$--modules $M$ have complexity $c_A(M)=0$. 
\item[(2)]
If there is an $A$--module $M$ with $c_A(M)\ge 2$ then 
$A$ is not finite.
\item[(3)]
If there is an $A$--module $M$ with $c_A(M)\ge 3$ then 
$A$ is wild.
\end{itemize}
\end{thm}

Recall that the Hecke algebras $\H_n(q,Q)$ and $\H_n^X(q)$, for $X=A, B, D$, 
and their block algebras are symmetric algebras. So we can apply 
Theorem \ref{complexity} to them. 

\begin{rem}
If $A$ is a group algebra then 
an indecomposable $A$--module $M$ has $c_A(M)=1$ if and 
only if $M$ is a periodic $A$--module which is 
not projective. See \cite[Proposition 8.4.4]{Ev} for example. 
However, this is not true 
for general finite dimensional algebras. A counterexample is given by 
\cite{LS}. See \cite{R4}. 
\end{rem}

\begin{rem}
The converse does not hold in (2) and (3) above. 
To see this, we consider group algebras and 
use the theory of support varieties. 
References for this theory are \cite{Ev}, \cite{B1}, \cite{B3}. 
Let $G$ be a finite group, $A=FG$ the group algebra. 
Assume that the characteristic of $F$ is $2$. Then 
\[
X_G=\operatorname{Spec}(H^*(G,F))
\]
is, by definition, the variety of the group $G$. 
For each $FG$--module $M$, we have the {\sf support variety} 
$X_G(M)$, which is a closed subvariety of $X_G$. 
Then, the famous theorem of Alperin and Evens \cite{AE2}, 
Avrunin \cite{Av} 
asserts that $X_G(M)$ is covered by affine charts labelled by 
elementary abelian $2$--subgroups \cite[Theorem 8.3.1]{Ev}. 
In particular, we have the following theorem of Quillen. 
\[
\operatorname{dim} X_G=\max\{\operatorname{rank}\;E\,|\,
\text{$E$ is an elementary $2$--subgroup of $G$.}\}.
\]

An important fact relevant to us is that the dimension of $X_G(M)$ is 
equal to the complexity of $M$ \cite[Theorem 8.4.3]{Ev}. 
Hence, the complexity 
of any $FG$--module $M$ cannot exceed the maximal rank of 
elementary $2$--subgroups of $G$. 

Assume that $A=FG$ is tame. Then, by \cite{BD} or 
\cite[Theorem 14.17]{S}, 
the Sylow $2$--subgroup of $G$ is one of 
dihedral, semidihedral or generalized 
quartanion groups. Take $G$ to be the quartanion group $Q$, which has 
order $8$. 
Then the maximal rank of elementary $2$--subgroups of $G$ 
is $1$. 
Thus, $c_A(M)\le 1$ for all $A$--modules $M$ and 
the converse of (2) fails. 
If we consider $G=Q\times C_2$ then we see that 
the converse of (3) also fails. 
\end{rem}

\subsection{}
In this subsection, we collect results to show that 
an $F$--algebra is wild. We start with an easy lemma. 

\begin{lem}
\label{criterion for wildness-1}
Let $A$ be a finite dimensional $F$--algebra, $Q$ its 
Gabriel quiver. Assume that $Q$ 
contains the following quiver, or the quiver with the reversed arrows, 
as a subquiver. 

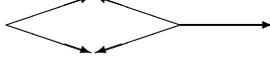
\setlength{\unitlength}{16pt}
\begin{picture}(15,2.4)(-3,0)
\put(3,1){\vector(3,1){2}}
\put(3,1){\vector(3,-1){2}}
\put(7.1,1){\vector(-3,1){2}}
\put(7.1,1){\vector(-3,-1){2}}
\put(7.1,1){\vector(1,0){2.2}}
\end{picture}

\noindent
Then $A$ is wild. 
\end{lem}
\begin{proof}
Let $T$ be the quiver given above. 
Assume that $Q$ contains $T$ as a subquiver. 
Then, we have a natural surjection 
\[
FQ \longrightarrow FT.
\]
Note that $FQ$ is an infinite dimensional $F$--algebra in general. 

Recall that 
$A$ is Morita--equivalent to the algebra $FQ/I$ such that 
$I$ is contained in the two--sided ideal consisting of pathes of 
length $\ge 2$ and contains all paths of length $\ge t$ for some $t\ge 2$. 
Because of the orientation of $T$, $FT$ is a square--zero algebra. Thus, 
the surjection induces 
\[
FQ/I \longrightarrow FT.
\]
As $FT$ is wild by Theorem \ref{path algebra case}(3), so is $FQ/I$. 
\end{proof}

A standard technique to show that an algebra is wild is the 
covering theory \cite{Gr}, \cite{Ga3}. See also \cite{W} and 
\cite{MP}. 

\begin{defn}
Let $Q=(Q_0,Q_1)$, $\tilde Q=(\tilde Q_0,\tilde Q_1)$ be 
directed graphs. We say that $\tilde Q$ is a {\sf covering} 
of $Q$ if there exist surjective maps of vertices and edges 
\[
\pi=(\pi_0, \pi_1):\;\tilde Q \longrightarrow Q
\]
such that, for any $\tilde x\in \tilde Q_0$, if we set $x=\pi(\tilde x)$ 
then $\pi$ induces bijection between 
$s(\tilde x)=\{\tilde y| (\tilde x, \tilde y)\in\tilde Q_1\}$ 
(resp. $e(\tilde x)=\{\tilde y| (\tilde y, \tilde x)\in\tilde Q_1\}$) and 
$s(x)=\{y| (x, y)\in Q_1\}$ (resp. $e(x)=\{y| (y, x)\in Q_1\}$). 

A covering $\tilde Q$ of $Q$ is a {\sf Galois covering} 
if there exists a subgroup $G$ of 
the automorphism group $\operatorname{Aut}(\tilde Q)$ such that, 
for any $x\in Q_0$, $\pi^{-1}(x)$ is a fixed point free $G$--orbit. 
\end{defn}

\begin{defn}
Let $\tilde Q$ be a Galois covering of $Q$ and 
$F\tilde Q \longrightarrow FQ$ the induced algebra 
homomorphism. Suppose that $\tilde I$ and $I$ are 
admissible ideals of $F\tilde Q$ and $FQ$ respectively. 
Then $\tilde A=F\tilde Q/\tilde I$ is a 
{\sf Galois covering} of $A=FQ/I$ if $\tilde I$ maps onto $I$. 
\end{defn}

In the definition above, we allow $\tilde A$ to be infinite dimensional. 
A useful result we use here is Han's covering criterion \cite{Ha1}. 
See \cite{D} for a different proof. 

\begin{thm}[{\cite[Theorem 3.3]{Ha1}}]
\label{Han}
Let $Q=(Q_0,Q_1)$ be a directed graph with a finite vertex set $Q_0$, 
$I$ an admissible ideal of $FQ$ such that 
$A=FQ/I$ is a finite dimensional $F$--algebra. Assume the following 
two conditions. 
\begin{itemize}
\item[(a)]
$\tilde A=F\tilde Q/\tilde I\longrightarrow A$ is a Galois covering 
whose Galois group $G$ is torsion free.
\item[(b)]
There is a subquiver $Q'$ of $\tilde Q$ such that, 
if we define an admissible ideal 
$I'$ of $FQ'$ by replacing each path which is not contained in $Q'$ 
with zero, in each element of $\tilde I$, then $A'=FQ'/I'$ is 
finite dimensional and strictly wild. 
\end{itemize}
Then, $A$ is controlled wild. 
\end{thm}

\begin{rem}
To know that $A$ is wild, it is not necessary to assume that 
the Galois group is torsion free, and it is enough to assume that 
$FQ'/I'$ is wild. This is because the pushdown functor of a Galois 
covering is a 
cleaving functor. See \cite[Chap.II, Lemma 2.1]{dlP1}. 
Then, we appeal to \cite[Chap.II, Proposition 2.1]{dlP1}. 
\end{rem}

Now we apply the covering criterion and prepare results which we need 
in later sections. We remark that 
there are proofs which do not use 
the covering criterion in the first two lemmas; 
we can construct a tensor functor 
from $F\langle X, Y\rangle-mod$ to the full subcategory 
$(FT-mod)_s$ of sincere $FT$--modules, where $T$ is a directed 
graph whose underlying graph is 
\[
\begin{cases}
\tilde{\tilde D}_4\; \text{in Lemma \ref{criterion for wildness-2},}\\
\tilde{\tilde E}_7\; \text{in Lemma \ref{criterion for wildness-3}.}
\end{cases}
\]
Then, the pushdown functor respects the indecomposability and the 
isomorphism classes of indecomposable $FT$--modules in 
$(FT-mod)_s$. In fact, this was the strategy in \cite[I.10.1--10.5]{Er}. 
However, the covering criterion gives us shorter proofs. 

\begin{lem}
\label{criterion for wildness-2}

Let $Q$ be the directed graph which is defined by the 
adjacency matrix $(a_{ij})_{1\le i,j\le 4}$ where 
\[
a_{ij}=\begin{cases} 1 \;\;\text{if $|i-j|=1$}\\
                     0 \;\;\text{otherwise.}\end{cases}
\]
We denote by $J^3$ the ideal of $FQ$ generated by paths of 
length $3$. Then $A=FQ/J^3$ is wild. 
\end{lem}
\begin{proof}
Let $\alpha, \beta, \gamma$ be the arrows $1\rightarrow 2$, 
$2\rightarrow 3$, $3\rightarrow 4$ and $\alpha', \beta', \gamma'$ 
the arrows with the opposite direction. 

We consider the following covering $\tilde Q$ of $Q$, which contains 
the quiver on the right hand side as a subquiver. 

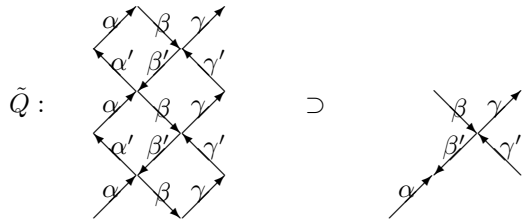
\setlength{\unitlength}{16pt}
\begin{picture}(15,7)(-3,-1)
\put(-2,2.5){$\tilde Q:$}
\put(0.2,4.5){$\alpha$}
\put(0,4){\vector(1,1){1}}
\put(1,3){\vector(-1,1){1}}
\put(0.4,3.5){$\alpha'$}
\put(0.2,2.5){$\alpha$}
\put(0,2){\vector(1,1){1}}
\put(1,1){\vector(-1,1){1}}
\put(0.4,1.5){$\alpha'$}
\put(0.2,0.5){$\alpha$}
\put(0,0){\vector(1,1){1}}

\put(1.05,5){\vector(1,-1){1}}
\put(1.45,4.4){$\beta$}
\put(1.25,3.5){$\beta'$}
\put(2.05,4){\vector(-1,-1){1}}
\put(1.05,3){\vector(1,-1){1}}
\put(1.45,2.4){$\beta$}
\put(1.25,1.5){$\beta'$}
\put(2.05,2){\vector(-1,-1){1}}
\put(1.05,1){\vector(1,-1){1}}
\put(1.45,0.4){$\beta$}

\put(2.1,4){\vector(1,1){1}}
\put(2.3,4.5){$\gamma$}
\put(2.6,3.5){$\gamma'$}
\put(3.1,3){\vector(-1,1){1}}
\put(2.1,2){\vector(1,1){1}}
\put(2.3,2.5){$\gamma$}
\put(2.6,1.5){$\gamma'$}
\put(3.1,1){\vector(-1,1){1}}
\put(2.3,0.5){$\gamma$}
\put(2.1,0){\vector(1,1){1}}
\put(5,2.5){$\supset$}

\put(7.2,0.5){$\alpha$}
\put(7,0){\vector(1,1){1}}

\put(8.05,3){\vector(1,-1){1}}
\put(8.45,2.4){$\beta$}
\put(8.25,1.5){$\beta'$}
\put(9.05,2){\vector(-1,-1){1}}

\put(9.1,2){\vector(1,1){1}}
\put(9.3,2.5){$\gamma$}
\put(9.6,1.5){$\gamma'$}
\put(10.1,1){\vector(-1,1){1}}
\end{picture}

$\tilde Q$ consists of $3$ zigzag lines, two going up and one going down, 
and the vertices with common $x$--coordinate is a fiber 
of $\tilde Q\rightarrow Q$. 
Define $\tilde I$ to be the ideal generated by paths 
of length $3$ in $\tilde Q$, and set $\tilde A=F\tilde Q/\tilde I$. 

Then, $\tilde A\rightarrow A$ is a Galois covering with Galois group 
$\mathbb Z$ and 
the condition (a) of Theorem \ref{Han} is satisfied. 
Denote the subquiver by $Q'$. Then, as there is no path 
of length greater than or equal to $3$ in $Q'$, we have $I'=0$ and 
$A'=FQ'$ is finite dimensional and strictly wild by 
Theorem \ref{path algebra case} and Proposition \ref{wild path algebra}. 
Thus, the condition (b) is also satisfied. Therefore, $A$ is wild by 
Theorem \ref{Han}. 
\end{proof}

\begin{lem}
\label{criterion for wildness-3}
Let $A$ be a finite dimensional $F$--algebra, 
$Q$ the Gabriel quiver of $A$. Suppose that $Q$ contains 
a subquiver with the adjacency matrix $(a_{ij})_{1\le i,j\le 3}$ 
where, 
\[
a_{ij}=
\begin{cases}
1\;\;\text{if $i\ge2$ and $j\ge2$ or $i=1$ and $j=2$,}\\
0\;\;\text{otherwise.}
\end{cases}
\]
Then $A$ is wild. 
\end{lem}
\begin{proof}
Let $\overline Q$ be the subquiver contained in $Q$. 
Define $\overline A=F\overline Q/J^2$ where $J^2$ is the ideal generated by 
paths of length $2$ in $\overline Q$. 
Then, there is a surjective algebra 
homomorphism $A\rightarrow \overline A$. Thus, it is enough to 
show that $\overline A$ is wild. We denote by 
$\rho, \mu, \nu$ the arrows 
$1\rightarrow 2$, $2\rightarrow 3$ and $2\leftarrow 3$, and 
the loops on the nodes $2$ and $3$ by $\alpha$ and $\beta$ 
respectively. 
Then, we can find a Galois covering $\tilde A=\tilde Q/\tilde I$ 
of $\overline A$ 
with Galois group $\mathbb Z$ such that, 
$\tilde Q$ contains the following quiver, which we define to be $Q'$, 
as a subquiver, 

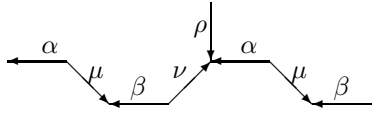
\setlength{\unitlength}{16pt}
\begin{picture}(15,4)(-5,0)

\put(1.4,2.2){$\alpha$}
\put(2,2){\vector(-1,0){1.4}}

\put(2.5,1.6){$\mu$}
\put(2,2){\vector(1,-1){1}}

\put(3.5,1.2){$\beta$}
\put(4.4,1){\vector(-1,0){1.4}}

\put(4.5,1.6){$\nu$}
\put(4.4,1){\vector(1,1){1}}

\put(6.1,2.2){$\alpha$}
\put(6.8,2){\vector(-1,0){1.4}}

\put(7.3,1.6){$\mu$}
\put(6.8,2){\vector(1,-1){1}}

\put(8.3,1.2){$\beta$}
\put(9.2,1){\vector(-1,0){1.4}}

\put(5,2.7){$\rho$}
\put(5.4,3.4){\vector(0,-1){1.4}}
\end{picture}

\noindent
and $\tilde I$ maps to $0$ under the surjective homomorphism 
$F\tilde Q\rightarrow FQ'$. 

Thus, $A'=FQ'$ and this is 
finite dimensional and strictly wild. 
Therefore, Theorem \ref{Han} implies that $\overline A$ is wild. 
\end{proof}

In the next lemma, we consider the following algebra. 
This is (32) of Han's Table W in \cite{Ha2}. 

\begin{defn}
\label{Lambda-32}
Let $Q$ be the directed graph defined by the adjacency matrix 
$\binom{1\; 1}{1\; 1}$, and let 
$FQ$ be the associated path algebra. 

We denote the loops on nodes $1$, $2$ by $\alpha$ and $\beta$, 
and the edges $1\rightarrow 2$, 
$1\leftarrow 2$ by $\mu$ and $\nu$ respectively. Then, 
$\Lambda_{32}$ is the $F$--algebra $FQ/I$ where the relations 
are given by 
\begin{gather*}
\mu\alpha=\beta\mu, \;\;\alpha\nu=0,\;\; \nu\beta=0,\\
\mu\nu\mu=0,\;\; \nu\mu\nu=0,\\
\alpha^2=0,\;\;\beta^2=0.
\end{gather*}
\end{defn}

\begin{lem}
\label{criterion for wildness-4}
$\Lambda_{32}$ is wild. 
\end{lem}
\begin{proof}
We consider the following covering $\tilde Q$ of $Q$, whose underlying 
graph is $\mathbb Z^2$. 

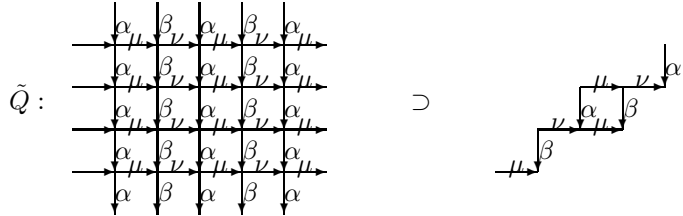
\setlength{\unitlength}{16pt}
\begin{picture}(15,7)(-3,-1)
\put(-1.5,2.5){$\tilde Q:$}
\put(1,4.3){$\alpha$}
\put(1,5){\vector(0,-1){1}}
\put(1,3.3){$\alpha$}
\put(1,4){\vector(0,-1){1}}
\put(1,2.3){$\alpha$}
\put(1,3){\vector(0,-1){1}}
\put(1,1.3){$\alpha$}
\put(1,2){\vector(0,-1){1}}
\put(1,0.3){$\alpha$}
\put(1,1){\vector(0,-1){1}}

\put(0,4){\vector(1,0){1}}
\put(0,3){\vector(1,0){1}}
\put(0,2){\vector(1,0){1}}
\put(0,1){\vector(1,0){1}}

\put(1,4){\vector(1,0){1}}
\put(1.3,4){$\mu$}
\put(1,3){\vector(1,0){1}}
\put(1.3,3){$\mu$}
\put(1,2){\vector(1,0){1}}
\put(1.3,2){$\mu$}
\put(1,1){\vector(1,0){1}}
\put(1.3,1){$\mu$}

\put(2,4){\vector(1,0){1}}
\put(2.3,4){$\nu$}
\put(2,3){\vector(1,0){1}}
\put(2.3,3){$\nu$}
\put(2,2){\vector(1,0){1}}
\put(2.3,2){$\nu$}
\put(2,1){\vector(1,0){1}}
\put(2.3,1){$\nu$}

\put(3,4){\vector(1,0){1}}
\put(3.3,4){$\mu$}
\put(3,3){\vector(1,0){1}}
\put(3.3,3){$\mu$}
\put(3,2){\vector(1,0){1}}
\put(3.3,2){$\mu$}
\put(3,1){\vector(1,0){1}}
\put(3.3,1){$\mu$}

\put(4,4){\vector(1,0){1}}
\put(4.3,4){$\nu$}
\put(4,3){\vector(1,0){1}}
\put(4.3,3){$\nu$}
\put(4,2){\vector(1,0){1}}
\put(4.3,2){$\nu$}
\put(4,1){\vector(1,0){1}}
\put(4.3,1){$\nu$}

\put(5,4){\vector(1,0){1}}
\put(5.3,4){$\mu$}
\put(5,3){\vector(1,0){1}}
\put(5.3,3){$\mu$}
\put(5,2){\vector(1,0){1}}
\put(5.3,2){$\mu$}
\put(5,1){\vector(1,0){1}}
\put(5.3,1){$\mu$}

\put(2,4.3){$\beta$}
\put(2,5){\vector(0,-1){1}}
\put(2,3.3){$\beta$}
\put(2,4){\vector(0,-1){1}}
\put(2,2.3){$\beta$}
\put(2,3){\vector(0,-1){1}}
\put(2,1.3){$\beta$}
\put(2,2){\vector(0,-1){1}}
\put(2,0.3){$\beta$}
\put(2,1){\vector(0,-1){1}}

\put(3,4.3){$\alpha$}
\put(3,5){\vector(0,-1){1}}
\put(3,3.3){$\alpha$}
\put(3,4){\vector(0,-1){1}}
\put(3,2.3){$\alpha$}
\put(3,3){\vector(0,-1){1}}
\put(3,1.3){$\alpha$}
\put(3,2){\vector(0,-1){1}}
\put(3,0.3){$\alpha$}
\put(3,1){\vector(0,-1){1}}

\put(4,4.3){$\beta$}
\put(4,5){\vector(0,-1){1}}
\put(4,3.3){$\beta$}
\put(4,4){\vector(0,-1){1}}
\put(4,2.3){$\beta$}
\put(4,3){\vector(0,-1){1}}
\put(4,1.3){$\beta$}
\put(4,2){\vector(0,-1){1}}
\put(4,0.3){$\beta$}
\put(4,1){\vector(0,-1){1}}

\put(5,4.3){$\alpha$}
\put(5,5){\vector(0,-1){1}}
\put(5,3.3){$\alpha$}
\put(5,4){\vector(0,-1){1}}
\put(5,2.3){$\alpha$}
\put(5,3){\vector(0,-1){1}}
\put(5,1.3){$\alpha$}
\put(5,2){\vector(0,-1){1}}
\put(5,0.3){$\alpha$}
\put(5,1){\vector(0,-1){1}}

\put(8,2.5){$\supset$}

\put(10,1){\vector(1,0){1}}
\put(10.3,1){$\mu$}

\put(11,1.3){$\beta$}
\put(11,2){\vector(0,-1){1}}

\put(11,2){\vector(1,0){1}}
\put(11.3,2){$\nu$}

\put(12,2.3){$\alpha$}
\put(12,3){\vector(0,-1){1}}

\put(12,3){\vector(1,0){1}}
\put(12.3,3){$\mu$}
\put(12,2){\vector(1,0){1}}
\put(12.3,2){$\mu$}

\put(13,2.3){$\beta$}
\put(13,3){\vector(0,-1){1}}

\put(13,3){\vector(1,0){1}}
\put(13.3,3){$\nu$}

\put(14,3.3){$\alpha$}
\put(14,4){\vector(0,-1){1}}
\end{picture}

We define $\tilde I$ by the same relations as for $\Lambda_{32}=FQ/I$, where 
monomials are understood as paths on $\tilde Q$, and 
$\mu\alpha=\beta\mu$ is understood as two paths with the same 
endpoints are equal. Denote the subquiver on the right hand side by $Q'$. 
Then, $I'=\langle \mu\alpha-\beta\mu\rangle$ and 
$FQ'/I'$ is the concealed algebra given in Lemma \ref{concealed algebra case}. 
Hence, $FQ'/I'$ is finite dimensional and strictly wild. 
Therefore, the covering criterion implies that $\Lambda_{32}$ is wild. 
\end{proof}

\subsection{}
There is a stable equivalence 
between an $F$--algebra with radical square zero and 
a path algebra. Further, 
given $F$--algebra with radical square zero, we can describe 
the directed graph of the corresponding path algebra. 
This is also a part of the Gabriel's theorem. 
See \cite[X.Theorem 2.4]{ARS} for example. 
We start with the following definition. 

\begin{defn}
Let $Q=(Q_0,Q_1)$ be a directed graph. Then, 
the {\sf associated directed bipartite graph} 
is the directed graph $\hat Q=(\hat Q_0,\hat Q_1)$ where 
$\hat Q_0$ is the disjoint union of two copies of $Q_0$, 
which we denote by $\{i',i''\}_{i\in Q_0}$, and $\hat Q_1$ 
is the set of arrows $i'\rightarrow j''$, one for each 
$i\rightarrow j$ in $Q_1$. 

We call the underlying 
graph of $\hat Q$ the {\sf associated bipartite graph}. 
\end{defn}

The following is the theorem of Gabriel. 

\begin{thm}
Let $A$ be a finite dimensional $F$--algebra, $Q$ 
the Gabriel quiver. Let $\hat Q$ be 
the directed bipartite graph associated with $Q$. Then, 
$A/\Rad^2A$ is stably equivalent to $F\hat Q$. 
\end{thm}

On the other hand, we have a theorem of Krause. This is 
a beautiful application of 
Crawley--Boevey's characterization of tameness in terms of 
the number of generic modules \cite{Ca2}. 
See \cite[Theorem 1]{Le} for the Krause's theorem. 
See also its final section for the proof and 
comments. The author thanks 
Professor Ringel for drawing his attension to \cite{Le}. 

The Krause's theorem asserts that stable equivalence preserves 
representation type, both tame and wild. 
Hence, Theorem \ref{path algebra case} implies the following. 

\begin{thm}
\label{stable equivalence}
Let $A$ be a finite dimensional $F$--algebra, $Q$ 
the Gabriel quiver. If the bipartite graph associated with $Q$ 
is not a Dynkin diagram of finite type, then $A$ is tame or wild. 
If it is not a Dynkin diagram of finite type nor affine type, 
then $A$ is wild. 
\end{thm}

Applying Theorem \ref{stable equivalence}, 
we have the following lemmas. The first lemma may be 
proved by using the covering criterion, by following the similar arguments 
as in Lemma \ref{criterion for wildness-2} and 
Lemma \ref{criterion for wildness-3}. 

\begin{lem}
\label{criterion for wildness-5}
Let $A$ be a finite dimensional $F$--algebra, and 
assume that the Gabriel quiver of $A$ 
contains the directed graph 
with adjacency matrix $\binom{2\;1}{0\;0}$ or 
$\binom{2\;0}{1\;0}$ as a subquiver. 
Then, $A$ is wild. 
\end{lem}
\begin{proof}
Let $Q$ be the subquiver. Then, 
there is a surjective algebra homomorphism 
$A\rightarrow FQ/J^2$ where $J^2$ is the ideal generated by 
paths of length $2$. 
As the associated bipartite graph is 
not a Dynkin diagram of finite type nor affine type, 
$FQ/J^2$ is wild. Thus, so is $A$. 
\end{proof}

\begin{lem}
\label{criterion for wildness-6}
Let $A$ be a finite dimensional $F$--algebra with the property that 
$\Ext_A^1(S,T)=0$ if and only if $\Ext_A^1(T,S)=0$, for any pair of simple 
$A$--modules $S$ and $T$. We denote by $\underline Q$ 
the underlying graph of the Gabriel quiver of $A$, and suppose that 
the following two conditions hold. 
\begin{itemize}
\item[(a)]
$\underline Q$ contains a cycle 
of length greater than or equal to $3$.
\item[(b)]
There is an edge of $\underline Q$ which is not contained in 
the cycle such that at least one of the 
endpoints of the edge belongs to the cycle. 
\end{itemize}
Then $A$ is wild. 
\end{lem}
\begin{proof}
Let $n$ be the length of the cycle. 
Assume that there is an edge which connects a 
node in the cycle and a node outside the cycle. 
If $n$ is even and $n\ge 4$ then 
the bipartite graph associated with $Q$ contains two copies of 
$A^{(1)}_{n-1}$ each of which is extended by an arrow. If $n$ is odd and 
$n\ge 3$ then the bipartite graph contains $A^{(1)}_{2n-1}$ which is 
extended by two arrows, whose endpoints on the cycle are distinct. 
Hence, the underlying graph of $\hat Q$ is not a Dynkin diagram of 
finite type nor affine type and $A$ is wild by 
Theorem \ref{stable equivalence}. 
The proof is entirely similar 
in the cases where there is an edge which 
connects two distinct nodes in the cycle, or 
there is a loop around a node of the cycle. 
\end{proof}

\subsection{}
Up to now, I have explained results to show that an algebra is wild. 
To show that an algebra is tame, we use the following result. 

\begin{defn}
A finite dimensional $F$--algebra is called {\sf special} 
if it is Morita--equivalent to 
a basic algebra $FQ/I$ with the following properties. 

\begin{itemize}
\item[(a1)]
For each $i\in Q_0$, the number of arrows $i\rightarrow j$ in $Q_1$ 
is at most $2$.
\item[(a2)]
For each $i\in Q_0$, the number of arrows $j\rightarrow i$ in $Q_1$ 
is at most $2$.
\item[(b1)]
For each directed edge $\alpha\in Q_1$, the number of arrows $\beta\in Q_1$ 
which satisfy $\alpha\beta\not\in I$ is at most $1$. 
\item[(b2)]
For each directed edge $\alpha\in Q_1$, the number of arrows $\beta\in Q_1$ 
which satisfy $\beta\alpha\not\in I$ is at most $1$. 
\end{itemize}

If, moreover, $I$ is given by a set of monomials, then 
$A$ is called a {\sf string algebra}. 
\end{defn}

\begin{thm}[{\cite[Corollary 2.4]{WW}}]
\label{special biserial}
Any special algebra is finite or tame. 
\end{thm}

The notion of special algebras first appeared in \cite{SW}. 
Recall that if, for any indecomposable projective $A$--module 
$P$, $P$ is uniserial or $\operatorname{Rad}P=N_1\oplus N_2/S$ 
with $N_1$ and $N_2$ uniserial and $S$ simple or zero, 
then $A$ is a {\sf biserial algebra}. It is known that 
special algebras are biserial \cite[Lemma 1]{SW}. Thus, we call 
special algebras {\sf special biserial algebras}. In fact, 
Theorem \ref{special biserial} is a special case of a more 
general theorem that biserial algebras are finite or tame. 
See \cite[Theorem A]{C3}. 

\medskip
Theorem \ref{special biserial} together with Lemma \ref{basic local} 
allow us to show that the special biserial algebras we will meet are 
tame. 

For special biserial algebras, 
we know their Auslander--Reiten quivers. As a corollary, 
we can prove the existence of a module with complexity $2$ 
in crucial cases. The author thanks 
Professor Erdmann for drawing his attention to this. 

Let $A$ be a finite dimensional basic $F$--algebra, 
$I=\{e_i|Ae_i\;\text{is injective.}\}$. 
Define $S=\oplus_{i\in I}\Soc Ae_i$. $S$ is a two--sided ideal 
of $A$. Define $\overline A=A/S$. The next lemma is standard. 
See \cite[I.8.11]{Er}. 

\begin{lem}
Let $A$, $I$, $S$ and $\overline A$ be as above. Then, 
the set of the isomorphism classes of indecomposable $A$--modules is the 
disjoint union of $\{Ae_i|\,i\in I\}$ and 
the set of the isomorphism classes of indecomposable 
$\overline A$--modules. The Auslander--Reiten quiver of 
$\overline A$ is obtained from that of $A$ by removing 
$\{Ae_i|\,i\in I\}$. 
\end{lem}
\begin{proof}
Fix $i\in I$ and write $P_i=Ae_i$. 
Let $M$ be an indecomposable $A$--module and assume that there is $m\in M$ 
such that $(\Soc P_i)m\ne 0$. Consider an $A$--module homomorphism 
\[
\phi\;:\;P_i \longrightarrow M
\]
defined by $ae_i\mapsto ae_im$, for $a\in A$. As $P_i$ is indecomposable 
and injective, $\Soc P_i$ is simple. Hence, if $\phi$ is not a monomorphism 
then $\Soc P_i\subset \Ker \phi$, contradicting to our assumption 
that $(\Soc P_i)m\ne 0$. Hence, we have $M\simeq P_i$. The first assertion 
is proved. 

Let $M$ and $N$ be $A$--modules on which $S$ act as zero. 
Suppose that there is an irreducible 
$\overline A$--module homomorphism $M\longrightarrow N$. 
Then, $M, N\not\simeq Ae_i$, for $i\in I$, implies that 
this homomorphism viewed as an $A$--module homomorphism 
is irreducible. 
Next assume that there is an irreducible 
$A$--module homomorphism $M\longrightarrow N$. Then, 
there exists an irreducible 
$\overline A$--module homomorphism $M\longrightarrow N$ 
by \cite[Corollary 3.3]{Ga2}. 
\end{proof}

Now assume that the basic algebra $A/S$ is a string algebra. 
For string algebras, we can classify 
their indecomposable modules. 
Among them, we only need the string modules for our purposes. 

Let $FQ/I$ be a string algebra. Consider 
a updown diagram of the following form 
with all the arrows directed downward. 

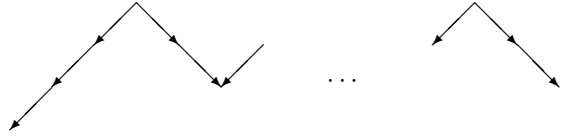
\setlength{\unitlength}{16pt}
\begin{picture}(15,5)(-3,-2)
\put(1,0){\vector(-1,-1){1}}
\put(2,1){\vector(-1,-1){1}}
\put(3,2){\vector(-1,-1){1}}
\put(3,2){\vector(1,-1){1}}
\put(4,1){\vector(1,-1){1}}
\put(6,1){\vector(-1,-1){1}}
\put(7.5,0){$\cdots$}
\put(11,2){\vector(-1,-1){1}}
\put(11,2){\vector(1,-1){1}}
\put(12,1){\vector(1,-1){1}}
\end{picture}

We see this as a collection of subpaths directed to southwest or southeast. 
Then, this updown diagram is called a {\sf string} if 
\begin{itemize}
\item[--]
edges of each subpath is labelled by elements of $Q_1$ so that 
the product is a monomial which is not in $I$, 
\item[--]
adjacent subpaths have a common source or a common target in $Q_0$ 
on a peak or in a deep, 
\item[--]
two arrows on a peak or in a deep have different labels.
\end{itemize}

Let $C$ be a string, $V(C)$ the set of its vertices. 
We define an $F$--vector space $M(C)$ by 
\[
M(C)=\bigoplus_{x\in V(C)}Fx. 
\]
Then, $M(C)$ becomes an $FQ/I$--module by the rule that 
$\alpha x=y$ if there is an arrow $x\rightarrow y$ 
with label $\alpha$ in $C$, and $\alpha x=0$ otherwise. 

\begin{defn}
Let $C$ be a string. We say that $C$ {\sf starts on a peak} 
(resp. {\sf ends on a peak}) if 
we cannot add an arrow directed to southwest 
(resp. southeast) to the right 
(resp. left) of $C$ so that the extended diagram may be a string again. 

Similarly, we say that $C$ {\sf starts in a deep} 
(resp. {\sf ends in a deep}) if 
we cannot add an arrow directed to southeast 
(resp. southwest) to the right 
(resp. left) of $C$ so that the extended diagram may be a string again. 
\end{defn}

Given $\alpha\in Q_1$, write $\alpha=i\rightarrow j$ and 
let $P_i$ be the indecomposable 
projective $FQ/I$--module corresponding to the node $i$, 
and let $P^j_i$ the $FQ/I$--submodule 
of $P_i$ generated by $\alpha$. 
Define $M_\alpha=P_i/P^j_i$. Note that $M_\alpha$ is uniserial. 
The following proposition gives an explicit rule for 
$\tau(M(C))$, where $\tau$ is the Auslander--Reiten translate. 
See \cite[II.3]{Er}, which is, in turn, based on \cite{BuR}. 

\begin{prop}
\label{AR translate}
Assume that $M(C)$ is not projective nor isomorphic to any of 
$M_\alpha$ $(\alpha\in Q_1)$. 
\begin{itemize}
\item[(1)]
If $C$ starts and ends in a deep then $\tau(M(C))=M(C')$ where 
$C'$ is obtained from $C$ by deleting 

\setlength{\unitlength}{16pt}
\begin{picture}(15,5)(1.2,-2)
\put(0,2.5){the leftmost}
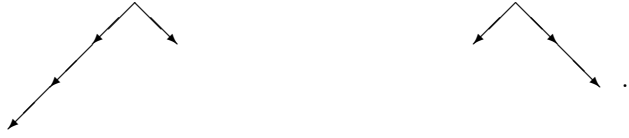
\put(5,0){\vector(-1,-1){1}}
\put(6,1){\vector(-1,-1){1}}
\put(7,2){\vector(-1,-1){1}}
\put(7,2){\vector(1,-1){1}}

\put(9,2.5){and the rightmost}

\put(16,2){\vector(-1,-1){1}}
\put(16,2){\vector(1,-1){1}}
\put(17,1){\vector(1,-1){1}}
\put(18.5,0){$.$}
\end{picture}

Thus, $C'$ does not start nor end on a peak. 
Note that the leftmost subpath directed to southwest and the rightmost 
subpath directed to southeast may be of length $0$. 

\item[(2)]
If $C$ starts in a deep but does not end in a deep 
then $\tau(M(C))=M(C')$ where 
$C'$ is obtained from $C$ by deleting the rightmost

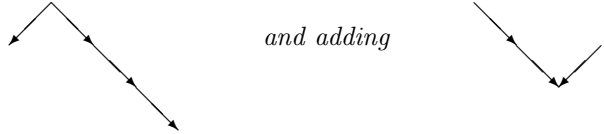
\setlength{\unitlength}{16pt}
\begin{picture}(15,5)(1.2,-2)
\put(2,2){\vector(-1,-1){1}}
\put(2,2){\vector(1,-1){1}}
\put(3,1){\vector(1,-1){1}}
\put(4,0){\vector(1,-1){1}}
\put(7,1){and adding}
\put(12,2){\vector(1,-1){1}}
\put(13,1){\vector(1,-1){1}}
\put(15,1){\vector(-1,-1){1}}
\end{picture}

\noindent
to the left of $C$ such that the length of the subpath directed to 
southeast is as large as possible. Thus, $C'$ does not start on a peak 
but ends on a peak. 
Note that the leftmost subpath of $C'$ directed to southeast and the rightmost 
subpath of $C$ directed to southeast may be of length $0$. 

\item[(3)]
If $C$ does not start in a deep but ends in a deep 
then $\tau(M(C))=M(C')$ where $C'$ is obtained from $C$ by deleting 
the leftmost

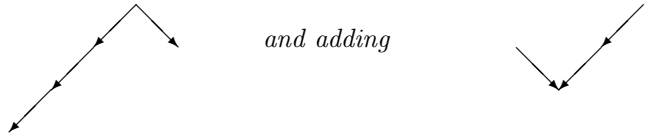
\setlength{\unitlength}{16pt}
\begin{picture}(15,5)(1.2,-2)
\put(2,0){\vector(-1,-1){1}}
\put(3,1){\vector(-1,-1){1}}
\put(4,2){\vector(-1,-1){1}}
\put(4,2){\vector(1,-1){1}}

\put(7,1){and adding}

\put(16,2){\vector(-1,-1){1}}
\put(15,1){\vector(-1,-1){1}}
\put(13,1){\vector(1,-1){1}}
\end{picture}

\noindent
to the right of $C$ such that the length of the subpath directed to 
southwest is as large as possible. Thus, $C'$ starts on a peak but 
does not end on a peak. 
Note that the rightmost subpath of $C'$ directed to 
southwest and the leftmost subpath of $C$ directed to 
southwest may be of length $0$. 

\item[(4)]
If $C$ does not start nor end in a deep then $\tau(M(C))=M(C')$ where 
$C'$ is obtained from $C$ by adding

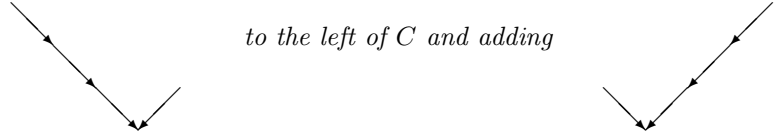
\setlength{\unitlength}{16pt}
\begin{picture}(15,5)(1.2,-2)
\put(1,2){\vector(1,-1){1}}
\put(2,1){\vector(1,-1){1}}
\put(3,0){\vector(1,-1){1}}
\put(5,0){\vector(-1,-1){1}}

\put(6.5,1){to the left of $C$ and adding}

\put(15,0){\vector(1,-1){1}}
\put(17,0){\vector(-1,-1){1}}
\put(18,1){\vector(-1,-1){1}}
\put(19,2){\vector(-1,-1){1}}
\end{picture}

\noindent
to the right of $C$ such that both the length of the subpath 
directed to southeast on the leftend and 
the length of the subpath directed to southwest on the 
rightend are as large as possible. Thus, $C'$ starts and ends on a peak. 
Note that these lengths may be of length $0$. 
\end{itemize}
\end{prop}

Let $Q$ be a directed graph with adjacency matrix 
$\binom{0\;1}{1\;1}$. We denote the arrows 
$1\rightarrow 2$, $1\leftarrow 2$ by $\mu$ and $\nu$ respectively, 
and the loop on the node $2$ by $\beta$. 
Define a string algebra $\Lambda=FQ/I$ by the relations
\[
\nu\beta=0,\;\;\beta\mu=0,\;\;\beta^2=0,\;\;(\mu\nu)^2=0,\;\;(\nu\mu)^2=0.
\]

\begin{lem}
\label{complexity=2(1)}
Let $S$ be the simple $\Lambda$--module corresponding to 
the node $2$. Then, 
$\operatorname{dim}_F\tau^{2n}(S)=6n+1$ and 
$\operatorname{dim}_F\tau^{2n+1}(S)=6n+6$, for $n\in\mathbb N$. 
\end{lem}
\begin{proof}
$\Lambda$ has the basis
\[
\{e_1,\;e_2,\;\beta,\;\mu,\;\nu,\;\mu\nu,\;\nu\mu,\;
\mu\nu\mu,\;\nu\mu\nu\}.
\]

$M_\mu$ is the simple $\Lambda$--module corresponding to the node $1$. Other 
$M_\nu$, $M_\beta$ and indecomposable projective $\Lambda$--modules are 
given by $M(C)$ with $C$ being one of the following. 

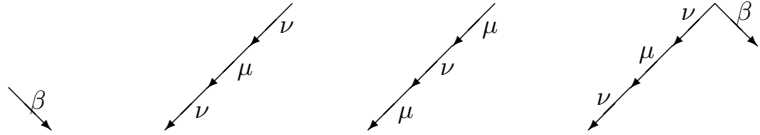
\setlength{\unitlength}{16pt}
\begin{picture}(15,5)(1,0)
\put(2.8,1.5){$\beta$}
\put(2.3,2){\vector(1,-1){1}}

\put(6.7,1.3){$\nu$}
\put(7,2){\vector(-1,-1){1}}
\put(7.7,2.3){$\mu$}
\put(8,3){\vector(-1,-1){1}}
\put(8.7,3.3){$\nu$}
\put(9,4){\vector(-1,-1){1}}

\put(11.5,1.3){$\mu$}
\put(11.8,2){\vector(-1,-1){1}}
\put(12.5,2.3){$\nu$}
\put(12.8,3){\vector(-1,-1){1}}
\put(13.5,3.3){$\mu$}
\put(13.8,4){\vector(-1,-1){1}}

\put(16.2,1.6){$\nu$}
\put(17,2){\vector(-1,-1){1}}
\put(17.2,2.7){$\mu$}
\put(18,3){\vector(-1,-1){1}}
\put(18.2,3.6){$\nu$}
\put(19,4){\vector(-1,-1){1}}
\put(19.5,3.6){$\beta$}
\put(19,4){\vector(1,-1){1}}
\end{picture}

\noindent
Thus, by Proposition \ref{AR translate}(2)(4), $\tau(S)$ and $\tau^2(S)$ are 
given by the following strings. 

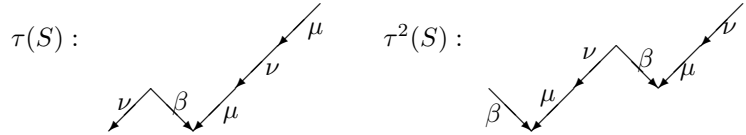
\setlength{\unitlength}{16pt}
\begin{picture}(15,5)(0,-0.5)
\put(0.7,2){$\tau(S):$}
\put(3.2,0.5){$\nu$}
\put(4,1){\vector(-1,-1){1}}
\put(4.5,0.5){$\beta$}
\put(4,1){\vector(1,-1){1}}
\put(5.7,0.3){$\mu$}
\put(6,1){\vector(-1,-1){1}}
\put(6.7,1.3){$\nu$}
\put(7,2){\vector(-1,-1){1}}
\put(7.7,2.3){$\mu$}
\put(8,3){\vector(-1,-1){1}}

\put(9.5,2){$\tau^2(S):$}
\put(11.9,0.2){$\beta$}
\put(12,1){\vector(1,-1){1}}
\put(13.2,0.7){$\mu$}
\put(14,1){\vector(-1,-1){1}}
\put(14.2,1.7){$\nu$}
\put(15,2){\vector(-1,-1){1}}
\put(15.5,1.5){$\beta$}
\put(15,2){\vector(1,-1){1}}
\put(16.5,1.3){$\mu$}
\put(17,2){\vector(-1,-1){1}}
\put(17.5,2.3){$\nu$}
\put(18,3){\vector(-1,-1){1}}
\end{picture}

\medskip
\noindent
It is clear that the same pattern repeats in 
$\tau^{2n+1}(S)$ and $\tau^{2n+2}(S)$, for $n>0$. 
\end{proof}

\begin{lem}
\label{complexity=2(2)}
Let $A$ be a symmetric special biserial algebra. 
Suppose that each node of the Gabriel quiver of $A$ 
has two outgoing arrows and two 
incoming arrows. Then, any simple $A$--module $S$ has 
the complexity $c_A(S)=2$. 
\end{lem}
\begin{proof}
As we mentioned before, $A$ is a biserial algebra. Thus, 
we have a description of the radical series of 
indecomposable projective 
$A$--modules. It follows that $A/S$ is a string 
algebra. As $A$ is self--injective, 
the Auslander--Reiten quiver of $\overline A$ coincides with 
the stable Auslander--Reiten quiver of $A$. 
By repeated use of Proposition \ref{AR translate}(4), 
we have the inequality 
\[
\operatorname{dim}_F\tau^n(S)\ge 2n+1\quad(n\in\mathbb N).
\]
Since $A$ is a symmetric algebra, 
$\tau=\Omega^2$ where $\Omega$ is the Heller 
loop operator. Thus, this inequality implies that $c_A(S)\ge 2$. 
As $c_A(S)\le 2$ by Theorem \ref{complexity}(3) 
and Theorem \ref{special biserial}, the result follows. 
\end{proof}

Let $A$ be a finite dimensional $F$--algebra, $\{P_1,\dots,P_s\}$ 
a complete set of the isomorphism classes of indecomposable projective 
$A$--modules. Then $A$ is symmetric if and only if 
$\End_A(P_1\oplus\cdots\oplus P_s)$ is symmetric. 
See \cite[Lemma I.3.3]{Er} for example. 
This fact will be used without further notice. 

\section{The case of $\H_n(q,Q)$}

\subsection{}
Let $F$ be an algebraically closed field 
and let $q, Q\in F^\times$. 
$\H_n(q,Q)$ is the $F$--algebra defined by generators 
$T_0,\dots,T_{n-1}$ and relations
\begin{gather*}
(T_0-Q)(T_0+1)=0, \quad
(T_i-q)(T_i+1)=0 \;\;(1\le i\le n-1),\\
T_0T_1T_0T_1=T_1T_0T_1T_0,\quad T_iT_j=T_jT_i \;\;(j\ge i+2),\\
T_iT_{i+1}T_i=T_{i+1}T_iT_{i+1}\;\;(i\ge1).
\end{gather*}

There are two cases to consider. The first case is the case 
where $-Q\not\in q^{\mathbb Z}$. 

\begin{thm}[{\cite[Theorem 4.17]{DJ2}}]
\label{Morita theorem}
Suppose that $Q\ne -q^f$ for any $f\in \mathbb Z$. 
Then $\H_n(q,Q)$ is Morita--equivalent to 
\[
\bigoplus_{k=0}^n\H^A_k(q)\otimes\H^A_{n-k}(q).
\]
\end{thm}

Further, we already know representation type for all of the 
block algebras of $\H_n^A(q)$ by \cite{EN}. 
Recall that block algebras of $\H^A_n(q)$ are labelled 
by $e$--cores $\kappa$ which satisfy the condition that 
$w_\kappa=\frac{n-|\kappa|}{e}$ is a non--negative integer 
\cite{DJ1}. 
Let $B_\kappa$ be the block algebra of $\H^A_n(q)$ labelled by 
an $e$--core $\kappa$. The integer 
$w_\kappa$ is called the {\sf weight} of $B_\kappa$. 
See \cite[2.7]{JK} for the terminology. 
The Erdmann--Nakano theorem is as follows. 

\begin{thm}[{\cite[Theorem 1.2]{EN}}]
\label{Erdmann-Nakano theorem}
Let $e$ be the multiplicative order of $q\ne1$ and 
$\kappa$ an $e$--core with 
$w_\kappa=\frac{n-|\kappa|}{e}\in{\mathbb Z}_{\ge0}$. Then, 
\begin{itemize}
\item[(1)]
$B_\kappa$ is semisimple if and only if 
$w_\kappa=0$. 
\item[(2)]
$B_\kappa$ is 
\begin{itemize}
\item[--]
finite if and only if $w_\kappa\le 1$, 
\item[--]
tame if and only if $e=2$ and $w_\kappa=2$, 
\item[--]
wild otherwise.
\end{itemize}
\end{itemize}
\end{thm}

Using this, we can prove the following. 

\begin{thm}
\label{separated parameter case}
Assume that $-Q\not\in q^{\mathbb Z}$ and that 
$q$ has the multiplicative order $e\ge 2$. 
Then, we have the following. 

\begin{itemize}
\item[(1)]
If $e\ge 3$ then $\H_n(q,Q)$ is 
\begin{itemize}
\item[--]
finite if and only if $n<2e$, 
\item[--]
wild otherwise. 
\end{itemize}
\item[(2)]
If $e=2$ then $\H_n(q,Q)$ is 
\begin{itemize}
\item[--]
finite if and only if $n<4$, 
\item[--]
tame if and only if $n=4$ or $5$, 
\item[--]
wild otherwise.
\end{itemize}
\end{itemize}
\end{thm}
\begin{proof}
(1) If $n<2e$ then one of $k$ and $n-k$ in Theorem 
\ref{Morita theorem} is strictly smaller than $e$. Thus, using 
Theorem \ref{Erdmann-Nakano theorem}(1), we know that one of 
$\H_k^A(q)$ and $\H_{n-k}^A(q)$ is a semisimple algebra, which is a direct 
sum of matrix algebras because $F$ is algebraically closed. 
Using Theorem \ref{Erdmann-Nakano theorem}(2) we also know that 
the other is finite. Hence $\H_n(q,Q)$ is finite in this case. 
To show that $n\ge 2e$ implies that $\H_n(q,Q)$ is wild, it is 
enough to prove that $\H_{2e}(q,Q)$ is wild by Corollary 
\ref{critical rank}. 
Theorem \ref{Erdmann-Nakano theorem}(2) implies that 
$\H_{2e}^A(q)$ is wild, which implies the result 
by Theorem \ref{Morita theorem}. 

\noindent
(2) By the assumption that $e=2$, 
the characteristic of $F$ is odd. First of all, 
Theorem \ref{Erdmann-Nakano theorem}(2) implies that 
$\H_2^A(q)$ and $\H_3^A(q)$ are finite. 
As $\H_0^A(q)=F$ and $\H_1^A(q)=F$ by definition, 
this implies 
that $\H_k^A(q)\otimes H_{2-k}^A(q)$, for $0\le k\le 2$, and 
$\H_k^A(q)\otimes H_{3-k}^A(q)$, for $0\le k\le 3$, 
are all finite, proving that $\H_2(q,Q)$ and $\H_3(q,Q)$ are finite. 

Next observe that $\H_4^A(q)$ and $\H_5^A(q)$ are tame. 
In fact, if $n=4$ then there is only one $2$--core, which has 
weight $2$, and if $n=5$ then there are two $2$--cores, one of 
which has weight $1$, the other of which has weight $2$. 
Thus, $\H_4^A(q)$ and $\H_5^A(q)$ are tame 
by Theorem \ref{Erdmann-Nakano theorem}(2). 
This implies that $\H_k^A(q)\otimes H_{4-k}^A(q)$, for 
$k=0,1,3,4$, and $\H_k^A(q)\otimes H_{5-k}^A(q)$, for 
$k=0,1,4,5$, are tame. As $\H_2^A(q)=F[T_1]/(T_1+1)^2$ implies that 
$\H_2^A(q)\otimes\H_2^A(q)$ is tame because 
\[
\H_2^A(q)\otimes\H_2^A(q)\simeq F[X,Y]/(X^2,Y^2), 
\]
we have proved that $\H_4(q,Q)$ is tame. 
On the other hand, $\H_3^A(q)$ is sum of two block algebras. 
One is isomorphic to the matrix algebra 
$\End_F(S^{(2,1)})$, where 
$S^{(2,1)}$ is the Specht module labelled by $(2,1)$, which is 
a projective $\H_3^A(q)$--module 
in the present case. The other is isomorphic to $F[X]/(X^2)$. 
Hence, $\H_k^A(q)\otimes \H_{5-k}^A(q)$, for $k=2$ and $k=3$, 
is Morita--equivalent 
to $F[X,Y]/(X^2,Y^2)\oplus F[Z]/(Z^2)$, proving that $\H_5(q,Q)$ is tame. 

As the block algebra of $\H_6^A(q)$ labelled by the empty $2$--core 
is wild, 
$\H_n^A(q)$ with $n\ge 6$ is wild by Corollary 
\ref{critical rank}. Thus, $\H_n(q,Q)$ is wild by 
Theorem \ref{Morita theorem}. 
\end{proof}

In the proof above, we heavily rely on Theorem \ref{Erdmann-Nakano theorem}. 
However, if we use the Fock space theory, which we will explain in the 
next subsection, we can give another proof which does not use 
Theorem \ref{Erdmann-Nakano theorem}. 

\subsection{}
The second case is the case where $-Q=q^f$ for some $f\in\mathbb Z$. 
Note that $\H_n(1,-1)$ is isomorphic to the semidirect product of 
$F[x_1,\dots, x_n]/(x_1^2,\dots, x_n^2)$ with the group algebra of 
the symmetric group $S_n$, where $S_n$ acts on $\{x_1,\dots, x_n\}$ in 
the natural way. Let $\ell$ be the characteristic of $F$, 
then.

\begin{prop}
$\H_n(1,-1)$ is 
\begin{itemize}
\item[--]
finite if $n=1$. 
\item[--]
tame if $n=2$ and $\ell\ne 2$.
\item[--]
wild otherwise.
\end{itemize}
\end{prop}
\begin{proof}
The case $n=1$ is obvious. As 
$F[x_1,\dots, x_n]/(x_1^2,\dots, x_n^2)$ is wild when $n\ge 3$, 
which follows from Theorem \ref{complexity} or 
Theorem \ref{stable equivalence}, 
Proposition \ref{reduction to critical rank} implies that 
$\H_n(1,-1)$ is wild in this case. 
The remaining cases are for $n=2$. Let $\sigma$ 
be the unique transposition of $S_2$. 

Suppose that $\ell\ne 2$. Then we have two simple modules 
$D^{\pm}$ defined by $x_i=0$, for $i=1, 2$, and 
$\sigma=\pm 1$. As the radical of $\H_2(1,-1)$ is 
$F[x_1, x_2]/(x_1^2, x_2^2)$, an explicit computation of the 
regular representation implies that 
the radical series of the projective cover of $D^{\pm}$ 
is as follows. 
\begin{equation*}
\begin{split}
&D^{\pm}\\
D^+&\oplus D^-\\
&D^{\pm}
\end{split}
\end{equation*}

Let $Q$ be the directed graph defined by the adjacency matrix 
$\binom{1\; 1}{1\; 1}$, and we denote the loops on nodes $1$, $2$ 
by $\alpha$ and $\beta$, the edges $1\rightarrow 2$, 
$1\leftarrow 2$ by $\mu$ and $\nu$. Then the Gabriel quiver of 
$\H_2(1,-1)$ is $Q$ and the relations include 
\[
\mu\alpha=\alpha\nu=0,\;\;\nu\beta=\beta\mu=0.
\]
Hence $\H_2(1,-1)$ is a special biserial algebra and 
Theorem \ref{special biserial} implies that $\H_2(1,-1)$ is 
finite or tame. On the other hand, Lemma \ref{complexity=2(2)} 
and Theorem \ref{complexity}(2) imply that $\H_2(1,-1)$ cannot 
be finite. Hence the result. 

Suppose that $\ell=2$. Then $\H_2(1,-1)$ is the projective cover 
of the unique simple module, and $\Rad \H_2(1,-1)/\Rad^2 \H_2(1,-1)$ 
constitutes of three simple modules. Hence $\H_2(1,-1)$ is wild by 
Theorem \ref{stable equivalence}.
\end{proof}

The theorem we are going to prove is the following. 

\begin{thm}
\label{two parameter case}
Assume that $q$ has the multiplicative order $e\ge 2$, and 
$Q=-q^f$, for some $0\le f\le e-1$. Then we have the following. 

\begin{itemize}
\item[(1)]
Suppose that $e\ge 3$. Then $\H_n(q,Q)$ is 
\begin{itemize}
\item[--]
finite if $n<\operatorname{min}\{e,\;2f+4,\;2e-2f+4\}$.
\item[--]
tame if $f=0$ and $4\le n<\operatorname{min}\{e,9\}$.
\item[--]
wild otherwise.
\end{itemize}
\item[(2)]
Suppose that $e=2$. Then $\H_n(q,Q)$ is 
\begin{itemize}
\item[--]
finite if $n=1$. 
\item[--]
tame if $n=2$ or $n=3$ and $f=1$.
\item[--]
wild otherwise.
\end{itemize}
\end{itemize}
\end{thm}

We have already proved in \cite[Theorem 1.4]{AM2} 
that $\H_n(q,-q^f)$ is finite if and only if 
$n<\operatorname{min}\{e,\;2f+4,\;2e-2f+4\}$, under the 
assumption that $e\ge 3$. The case $e=2$ is easy to handle, and 
the proof is given in \cite{AM3}. We remark that 
the radical series given in $(case\;1)$ in the proof of 
\cite[Theorem 5.25]{AM2} needs correction as in 
(case 5a) below. 

Hence, it is enough to prove the statements for tameness 
and wildness. 
To prove these, we need the Specht module theory for $\H_n(q,Q)$, 
developed by Dipper, James and Murphy \cite{DJM}, and 
the Fock space theory developed by the author \cite{A1}, \cite{A2} 
and \cite{A3}, as for the proof of finiteness in \cite[Theorem 1.4]{AM2}. 
We redefine $T_0$ by $T_0^{new}=-T_0^{old}$ if $0\le f\le\frac{e}{2}$, and 
$T_0^{new}=-q^{-f}T_0^{old}$ if $e$ is finite and $\frac{e}{2}<f\le e-1$. 
Thus, we may and do assume $0\le f\le\frac{e}{2}$ without loss of 
generality, and we assume that $T_0$ satisfies $(T_0-1)(T_0-q^f)=0$ 
in the rest of the paper. 

\begin{rem}
The Specht module theory was generalized to the Hecke algebras 
of type $(m,1,n)$ by Dipper, James and Mathas \cite{DJM'}. 
This theory is now viewed as an example of the 
cell module theory for cellular algebras in the sense of 
Graham and Lehrer. Also, if we restrict ourselves to type $A$ case, 
the Fock space theory was generalized to 
$q$--Schur algebras by Varagnolo and Vasserot \cite{VV}. 
\end{rem}

Let us begin by explaining the Specht module theory for $\H_n(q,Q)$. 
Let $\mathcal{BP}$ be the set of bipartitions 
$\lambda=(\lambda^{(1)},\lambda^{(2)})$. The size 
of $\lambda$ is denoted by $|\lambda|$. If $|\lambda|=n$ then 
we write $\lambda\vdash n$. Let $\mathcal{BP}(n)$ be the set 
of bipartitions with $\lambda\vdash n$. 
Then $\mathcal{BP}(n)$ is a poset with a partial order 
$\trianglelefteq$, called 
the {\sf dominance ordering}; we say that 
$\lambda\trianglelefteq\mu$ if
\begin{gather*}
\sum_{i=1}^k\lambda^{(1)}_i\le\sum_{i=1}^k\mu^{(1)}_i, \;
\mbox{for all $k$, and}\\
|\lambda^{(1)}|+\sum_{j=1}^k\lambda^{(2)}_j\le
|\mu^{(1)}|+\sum_{j=1}^k\mu^{(2)}_j, \;\mbox{for all $k$}. 
\end{gather*}
If $\lambda\trianglelefteq\mu$ and $\lambda\ne\mu$ 
we write $\lambda\triangleleft\mu$.

In \cite{DJM} it was shown that 
there exists a family $\{S^\lambda|\lambda\vdash n\}$ of 
$\H_n(q,Q)$--modules such that each of which is equipped with 
a symmetric bilinear form $\langle \;,\;\rangle$ satisfying 
$\langle T_iu,v\rangle=\langle u,T_iv\rangle$, for 
$0\le i\le n-1$ and for all $u,v\in S^\lambda$.
We define an $\H_n(q,Q)$--module $D^\lambda$ by 
$D^\lambda=S^\lambda/\operatorname{rad}_{\langle\;,\;\rangle}S^\lambda$. 
If there is a need to specify the base field $F$, we write 
$D^\lambda_F$. 
Then $D^\lambda$ is absolutely irreducible or $D^\lambda=0$. 
If $D^\lambda\ne0$ then its projective cover 
is denoted by $P^\lambda$ or $P^\lambda_F$. 

The following results are fundamental in the Specht module theory. 

\begin{thm}
\begin{itemize}
\item[(1)]
Each $D^\lambda$ is self--dual. 
\item[(2)]
$\{D^\lambda\,|\, \lambda\vdash n, D^\lambda\ne 0\}$ is a comlete set 
of the isomorphism classes of simple $\H_n(q,Q)$--modules. 
\item[(3)]
If $D^\lambda\ne 0$ then $S^\lambda$ is indecomposable with 
$\Top S^\lambda=D^\lambda$. 
\item[(4)]
For each $\mu\vdash n$ with $D^\mu\ne 0$, 
$P^\mu$ has a Specht filtration 
\[
P^\mu=F_0\supset F_1\supset\cdots \supset F_N=0
\]
such that $S^\lambda$ appears 
$[S^\lambda:D^\mu]$ times in $\{ F_i/F_{i+1} | 0\le i<N\}$, 
for each $\lambda$. Further, if we denote 
the maximal element of $\{ \lambda\,|\,[S^\lambda:D^\mu]\ne 0\}$ 
in the dominance ordering by $\lambda_{\rm max}$, then 
$S^{\lambda_{\rm max}}$ appears as a submodule of $P^\mu$. 
\end{itemize}
\end{thm}

The statement (1) implies that if there is an arrow in the 
Gabriel quiver of $\H_n(q,Q)$ then there is always an arrow 
with the same endpoints and the opposite direction in the 
Gabriel quiver. 

Another important remark about the Specht module theory is that 
if we have enough knowledge on 
the decomposition numbers $[S^\lambda:D^\mu]$, 
then, by using seminormal representations, 
we can obtain an explicit matrix representation of $D^\mu$. 

We denote the set of bipartitions $\lambda\in\mathcal{BP}(n)$ with 
$D^\lambda\ne0$ by $\mathcal{KBP}(n)$, and 
the disjoint union $\sqcup_{n\in\mathbb N}\mathcal{KBP}(n)$ 
by $\mathcal{KBP}$. Then we can describe $\mathcal{KBP}$ by using 
the crystal graph theory. Further, the decomposition numbers 
$[S^\lambda:D^\mu]$ can be described in terms of the canonical basis 
in a combinatorial Fock space. This is the theory which I called 
the Fock space theory. To explain this, 
we review the author's previous work \cite{A2}, which sprang from 
an important observation and results of Lascoux, Leclerc and Thibon 
\cite{LLT}. See chapters 12--14 of \cite{A1} for details. 
Here, we state the results only for $\H_n(q,Q)$. 

Assume as before that $q\in F^\times$, which appeared as one of the parameters 
$q$ and $Q$ of $\H_n(q,Q)$, is a primitive $e^{th}$ root of 
unity with $e\ge 2$. 
Then we consider the Kac--Moody Lie algebra 
$\mathfrak g$ of type $A^{(1)}_{e-1}$ and its quantized enveloping 
algebra $U_v(\mathfrak g)$, where $v$ is an indeterminate. 
We denote by $V_v(\Lambda)$ the irreducible 
integrable highest weight $U_v(\mathfrak g)$--module with 
highest weight $\Lambda=\Lambda_0+\Lambda_f$. 
Recall that $f$ is such that $Q=-q^f$ and $0\le f\le\frac{e}{2}$. 

The {\sf $v$--deformed combinatorial Fock space} $\mathcal F_v(\Lambda)$ 
with hightest weight $\Lambda$ is the infinite dimensional vector space 
\[
\mathcal F_v(\Lambda)=\bigoplus_{\lambda\in\mathcal{BP}}\mathbb Q(v)\lambda,
\]
which is a $U_v(\mathfrak g)$--module via 
\begin{gather*}
e_i\lambda=\sum_{\mu\in{\mathcal{BP}}:\lambda/\mu=\fbox{\it i}}
v^{-N^a_i(\lambda/\mu)}\mu, \\
f_i\lambda=\sum_{\mu\in{\mathcal{BP}}:\mu/\lambda=\fbox{\it i}}
v^{N^b_i(\mu/\lambda)}\mu, \\
t_i\lambda=v^{N_i(\lambda)}\lambda,\quad 
v^d\lambda=v^{-W_0(\lambda)}\lambda.
\end{gather*}
The definitions of $N^a_i(\lambda/\mu)$, $N^b_i(\mu/\lambda)$, 
$N_i(\lambda)$ and $W_0(\lambda)$ in the formulas above 
are as in \cite[Definition 10.8]{A1} and we omit them. 
The proof that this is a $U_v(\mathfrak g)$--module is given 
in \cite[Theorem 10.10]{A1}. 

The following is called the Misra--Miwa theorem. See 
\cite[Theorem 11.11]{A1}. 

\begin{thm}
Define $\mathcal L(\Lambda)$ and $\mathcal B(\Lambda)$ by 
\[
\mathcal L(\Lambda)=\bigoplus_{\lambda\in\mathcal{BP}}
\mathbb Q[v]_{(v)}\lambda, \quad
\mathcal B(\Lambda)=
\{\lambda+v\mathcal L(\Lambda)\,|\,\lambda\in\mathcal{BP}\}.
\]
Then $(\mathcal L(\Lambda),\mathcal B(\Lambda))$ is a (lower) crystal 
basis of $\mathcal F_v(\Lambda)$ in the sense of Kashiwara. 
\end{thm}
We identify the crystal graph of $(\mathcal L(\Lambda),\mathcal B(\Lambda))$ 
with $\mathcal{BP}$. 
Now we identify the $U_v(\mathfrak g)$--submodule of 
$\mathcal F_v(\Lambda)$ generated by 
the empty bipartition $\emptyset$ with $V_v(\Lambda)$. 
Then $(\mathcal L(\Lambda),\mathcal B(\Lambda))$ 
defines a crystal basis of $V_v(\Lambda)$, which 
we denote by $(L(\Lambda), B(\Lambda))$. Note that 
$V_v(\Lambda)$ is a $U_v(\mathfrak g)$--submodule of $\mathcal F_v(\Lambda)$ 
and $B(\Lambda)$ is a subset of $\mathcal{BP}$ in our definition. 
The highest weight vector $v_{\Lambda}$ is identified with the empty 
bipartition $\emptyset\in\mathcal F_v(\Lambda)$. By a fundamental theorem 
of Kashiwara and Lusztig, the crystal basis $(L(\Lambda), B(\Lambda))$ 
can be lifted to the canonical basis of $V_v(\Lambda)$. 

The following collects some of the main theorems of the Fock space theory. 
The second part was proved in \cite[Theorem 4.4]{A2}; the proof 
requires plenty of results. I recommend reading 
the proof of \cite[Theorem 12.5]{A1}, which not only proves the second part 
but also includes background materials necessary for the proof. 
The remaining parts were proved in \cite[Theorem 4.2]{A3} and 
\cite[Corollary 3.16]{AM2}. 

We say that $(K,\mathcal O,F)$ is a {\sf modular system with parameters} if 
$\mathcal O$ is a discrete valuation ring and $K$ the fraction field, 
$F$ the residue field, such that there are elements 
$\hat q, \hat Q\in\mathcal O^\times$ which are lift 
of the parameters $q, Q\in F^\times$. 

\begin{thm}
\label{Fock space theory}
\begin{itemize}
\item[(1)]
$\mathcal B(\Lambda)=\mathcal{KBP}$. 
\item[(2)]
Assume that the characteristic of $F$ is zero. 
Let $\mu\in\mathcal{KBP}$ and denote 
the corresponding canonical basis element by $G(\mu)$. 
If we write 
\[
G(\mu)=\sum_{\lambda\in\mathcal{BP}}d_{\lambda\mu}(v)\lambda
\]
in $\mathcal F_v(\Lambda)$ then we have 
\[
d_{\lambda\mu}(1)=[S^\lambda:D^\mu].
\]
\item[(3)]
Let $(K,\mathcal O,F)$ be a modular system with parameters and 
suppose that $\mu$ is a Kleshchev bipartition such that 
\[
G(\mu)=f_{i_1}^{(m_1)}\dots f_{i_n}^{(m_n)}v_\Lambda
\]
for some $m_1,\dots,m_n\in\mathbb N$ and 
$i_1,\dots,i_n\in\mathbb Z/e\mathbb Z$. Then the decomposition
map sends $[P^\mu_K]$ to $[P^\mu_F]$. In particular, the decomposition 
numbers in the column labelled by $\mu$ do not depend on the 
characteristic of $F$. 
\end{itemize}
\end{thm}

\subsection{}
Now we start the proof of Theorem \ref{two parameter case}(1). 
As we assume that $e\ge 3$ and $0\le f\le\frac{e}{2}$, we have 
$e-f\ge 2$. We prove the theorem by case--by--case analysis. The cases 
we consider are as follows. 

\begin{itemize}
\item[(case 1a)]
$n=e$, $e-f\ge 3$ and $1\le f\le \frac{e}{2}$. 
\item[(case 2a)]
$n=e$, $e-f=2$ and $1\le f\le \frac{e}{2}$. 
\item[(case 3a)]
$n=e$ and $f=0$. 
\item[(case 4a)]
$n=2f+4$, $e>2f+4$ and $1\le f\le \frac{e}{2}$. 
\item[(case 5a)]
$4\le n<\operatorname{min}\{e, 9\}$ and $f=0$.
\item[(case 6a)]
$n=9<e$ and $f=0$.
\end{itemize}

Our aim is to show that $\H_n(q,Q)$ is wild in all the cases but (case 5a). 
In (case 5a) we show that $\H_n(q,Q)$ is tame and has an indecomposable module 
with complexity $2$. 

\subsection{}
The first three cases are for $n=e$. 
Then, following \cite{AM2}, we define bipartitions 
\begin{align*}
  \lambda_k&=\left((0),(k,1^{e-k})\right)\quad(1\le k\le e),\\
  \mu_k&=\left((k,1^{e-k}),(0)\right)\quad(1\le k\le e),\\
  \lambda_{k,l}&=\left((f-l,1^{e-f-k}),(k,1^l)\right)\quad
           (1\le k\le e-f\;\text{and}\; 0\le l<f).
\end{align*}
These bipartitions belong to the \lq\lq principal block\rq\rq, 
which we denote by $B$. 

\begin{prop}[{\cite[Proposition 4.22]{AM2}}]
\label{f>0}
We have the following. 

\begin{enumerate}
\item[(1)]
The complete set of Kleshchev bipartitions in $B$ is
\begin{equation*}
\{\lambda_k|1\le k<e\}
  \cup\{\lambda_{k,l}|1\le k\le e-f\And 0\le l<f\}.
\end{equation*}
\item[(2)]
For $1\le k<e$ we have
\begin{equation*}
[P^{\lambda_k}]=[S^{\lambda_k}]+[S^{\lambda_{k+1}}]
    +\begin{cases}%
 [S^{\lambda_{k,f-1}}]+[S^{\lambda_{k+1,f-1}}]\quad&(k<e-f)\\[2pt]
 [S^{\lambda_{e-f,f-1}}] &(k=e-f)\\[2pt]
 [S^{\lambda_{e-f,e-k}}]+[S^{\lambda_{e-f,e-k-1}}]&(k>e-f)
\end{cases}.
\end{equation*}
\item[(3)]
For $1\le k\le e-f$ and $0\le l<f$ we have
\begin{equation*}
[P^{\lambda_{k,l}}]=[S^{\lambda_{k,l}}]+
\begin{cases} [S^{\lambda_{k-1,0}}]+[S^{\mu_{f+k-1}}]+[S^{\mu_{f+k}}]\quad
         &(k\ne1\And l=0)\\[2pt]
[S^{\mu_f}]+[S^{\mu_{f+1}}]
         &(k=1\And l=0)\\[2pt]
[S^{\lambda_{k,l-1}}]+[S^{\lambda_{k-1,l}}]+[S^{\lambda_{k-1,l-1}}]
         &(k\ne1\And l\ne0)\\[2pt]
[S^{\lambda_{1,l-1}}]+[S^{\mu_{f-l}}]+[S^{\mu_{f-l+1}}]&
    (k=1\And l\ne0)
\end{cases}.
\end{equation*}
\end{enumerate}
\end{prop}

\bigskip
Assume that we are in (case 1a). The argument is very similar to 
that in \cite{A5}, which was for $\H_n^B(q)$. By Proposition \ref{f>0}, 
\begin{align*}
[P^{\lambda_1}]&=
[S^{\lambda_1}]+[S^{\lambda_2}]+[S^{\lambda_{2,f-1}}]+[S^{\lambda_{1,f-1}}], 
\\
[P^{\lambda_2}]&=
[S^{\lambda_2}]+[S^{\lambda_3}]+[S^{\lambda_{3,f-1}}]+[S^{\lambda_{2,f-1}}], 
\end{align*}
and
\begin{align*}
[S^{\lambda_{1,f-1}}]&=
[D^{\lambda_{1,f-1}}]+[D^{\lambda_{2,f-1}}]+[D^{\lambda_1}], 
\\
[S^{\lambda_{2,f-1}}]&=
[D^{\lambda_{2,f-1}}]+[D^{\lambda_{3,f-1}}]+[D^{\lambda_2}]
+[D^{\lambda_1}].
\end{align*}
Then, that $S^{\lambda_{1,f-1}}$ is a submodule of $P^{\lambda_1}$ 
implies that the radical series of $S^{\lambda_{1,f-1}}$ is as follows. 
\begin{equation*}
\begin{split}
&D^{\lambda_{1,f-1}}\\
&D^{\lambda_{2,f-1}}\\
&D^{\lambda_1}
\end{split}
\end{equation*}
On the other hand, 
$S^{\lambda_2}$ has the radical series of the following form. 
\begin{equation*}
\begin{split}
&D^{\lambda_2}\\
&D^{\lambda_1}
\end{split}
\end{equation*}
Then, as in the proof of (case 1) of \cite[Theorem 4.21]{AM2}, 
we have 
\[
\Rad P^{\lambda_1}/\Rad^2 P^{\lambda_1}=D^{\lambda_2}\oplus 
D^{\lambda_{2,f-1}}. 
\]
Note that $D^{\lambda_{1,f-1}}$ cannot appear in 
$\Rad P^{\lambda_1}/\Rad^2 P^{\lambda_1}$ because 
if otherwise then $D^{\lambda_{1,f-1}}$ would appear in 
$\Soc^2 P^{\lambda_1}/\Soc P^{\lambda_1}$. This contradicts 
to $[P^{\lambda_1}:D^{\lambda_{1,f-1}}]=1$ as $D^{\lambda_{1,f-1}}$ 
has already appeared in $\Soc^3 P^{\lambda_1}/\Soc^2 P^{\lambda_1}$. 

Another property we need is $\Soc S^{\lambda_{2,f-1}}=D^{\lambda_2}$. 
This follows from the fact 
that $S^{\lambda_{2,f-1}}$ is a submodule of $P^{\lambda_2}$. 

Now, we determine the radical series of $S^{\lambda_{2,f-1}}$. 
If $\Rad^2 S^{\lambda_{2,f-1}}=0$ then $\Soc S^{\lambda_{2,f-1}}$ is 
the direct sum 
$D^{\lambda_{3,f-1}}\oplus D^{\lambda_1}\oplus D^{\lambda_2}$, 
a contradiction. 

If $\Rad^3 S^{\lambda_{2,f-1}}\ne 0$ then $S^{\lambda_{2,f-1}}$ is uniserial 
whose top is $D^{\lambda_{2,f-1}}$ and whose socle is $D^{\lambda_2}$. 
This implies that there exists a uniserial module of the following form. 
\begin{equation*}
\begin{split}
&D^{\lambda_1}\\
&D^{\lambda_{3,f-1}}
\end{split}
\end{equation*}
This contradicts to 
$\Rad P^{\lambda_1}/\Rad^2 P^{\lambda_1}=D^{\lambda_2}\oplus 
D^{\lambda_{2,f-1}}$. 
Hence, the radical length of $S^{\lambda_{2,f-1}}$ is $3$ and 
the radical series has the following form. 
\begin{equation*}
\begin{split}
& D^{\lambda_{2,f-1}}\\
D^{\lambda_1}&\oplus D^{\lambda_{3,f-1}}\\
& D^{\lambda_2}
\end{split}
\end{equation*}
As a conclusion, the Gabriel quiver of $\H_n(q,Q)$ contains the following 
quiver as a subquiver. 

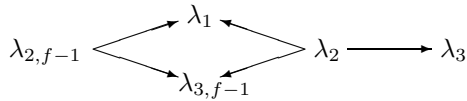
\setlength{\unitlength}{16pt}
\begin{picture}(15,3.5)(-3,-1)
\put(1,0.8){$\lambda_{2,f-1}$}
\put(3,1){\vector(3,1){2}}
\put(5.2,1.6){$\lambda_1$}
\put(3,1){\vector(3,-1){2}}
\put(5,0){$\lambda_{3,f-1}$}
\put(8,1){\vector(-3,1){2}}
\put(8,1){\vector(-3,-1){2}}
\put(8.2,0.8){$\lambda_2$}
\put(9,1){\vector(1,0){2}}
\put(11.2,0.8){$\lambda_3$}
\end{picture}

\noindent
Therefore, Lemma \ref{criterion for wildness-1} shows that 
$\H_n(q,Q)$ is wild in this case. 

\bigskip
Next assume that we are in (case 2a). Then we have either 
$e=3$ and $f=1$ or $e=4$ and $f=2$. 

First we consider the case where $e=3$ and $f=1$. 
Then, by Proposition \ref{f>0}, 
the decomposition matrix of the block algebra $B$ is 
as follows. 

\medskip
\begin{center}
\begin{tabular}{c|cccc}
 & $\lambda_1$ & $\lambda_2$ & $\lambda_{2,0}$ & $\lambda_{1,0}$ \\\hline
$\lambda_1$      & 1 & 0 & 0 & 0 \\
$\lambda_2$      & 1 & 1 & 0 & 0 \\
$\lambda_3$      & 0 & 1 & 0 & 0 \\
$\lambda_{2,0}$  & 1 & 1 & 1 & 0 \\
$\lambda_{1,0}$  & 1 & 0 & 1 & 1 \\
$\mu_1$          & 0 & 0 & 0 & 1 \\
$\mu_2$          & 0 & 0 & 1 & 1 \\
$\mu_3$          & 0 & 0 & 1 & 0 \\
\end{tabular}
\end{center}

\medskip
Thus, simple $B$--modules are all one dimensional and they are 
given by
\begin{alignat*}{2}
D^{\lambda_1}:\;& T_0\mapsto q, \quad& T_1,\, T_2 \mapsto -1, \\
D^{\lambda_2}:\;& T_0\mapsto q, \quad& T_1,\, T_2 \mapsto q, \\
D^{\lambda_{2,0}}:\;& T_0\mapsto 1, \quad& T_1,\, T_2 \mapsto q, \\
D^{\lambda_{1,0}}:\;& T_0\mapsto 1, \quad& T_1,\, T_2 \mapsto -1. 
\end{alignat*}
This is easy to see: observe that $D^{\lambda_2}=S^{\lambda_3}$, 
$D^{\lambda_{2,0}}=S^{\mu_3}$ and $D^{\lambda_{1,0}}=S^{\mu_1}$. 

The block algebra $B$ admits a symmetry. To see this, 
we define an algebra automorphism $\omega$ of $\H_n(q,Q)$ as follows. 
\begin{equation*}
\omega\;:\;\begin{cases} T_0\;\; \mapsto &1+q-T_0\;(=qT_0^{-1})\\
T_i\;\; \mapsto & q-1-T_i\;(=-qT_i^{-1})\quad\text{(for $i=1,2$)}
\end{cases}
\end{equation*}

Then, $\omega$ interchanges $D^{\lambda_1}$ and $D^{\lambda_{2,0}}$, 
$D^{\lambda_2}$ and $D^{\lambda_{1,0}}$, respectively. 

\begin{lem}
\label{ext}
$\Ext_B^1(D^{\lambda_a},D^{\lambda_b})\ne 0$ if and 
only if $\{\lambda_a,\lambda_b\}$ is one of 
\[
\{\lambda_1,\lambda_2\},\;\;\{\lambda_1,\lambda_{2,0}\},\;\;
\{\lambda_{1,0},\lambda_{2,0}\}.
\]
\end{lem}
\begin{proof}
By the self--duality of simple modules, 
$\Ext_B^1(D^{\lambda_a},D^{\lambda_b})\ne 0$ implies 
$\Ext_B^1(D^{\lambda_b},D^{\lambda_a})\ne 0$, for any pair 
$\{\lambda_a,\lambda_b\}$. 

As $S^{\lambda_2}$ is an indecomposable module 
with $\Top S^{\lambda_2}=D^{\lambda_2}$ 
and $\Soc S^{\lambda_2}=D^{\lambda_1}$, we have 
$\Ext_B^1(D^{\lambda_1},D^{\lambda_2})\ne 0$ and 
$\Ext_B^1(D^{\lambda_2},D^{\lambda_1})\ne 0$. 
Applying $\omega$, we also have 
$\Ext_B^1(D^{\lambda_{1,0}},D^{\lambda_{2,0}})\ne 0$ and 
$\Ext_B^1(D^{\lambda_{2,0}},D^{\lambda_{1,0}})\ne 0$. 

To prove that $\Ext_B^1(D^{\lambda_1},D^{\lambda_{2,0}})\ne 0$, 
we consider 
\[
T_0 \mapsto \begin{pmatrix} 1 & 1-q \\ 0 & q \end{pmatrix},\;\;
T_1 \mapsto \begin{pmatrix} q & 0 \\ 0 & -1 \end{pmatrix},\;\;
T_2 \mapsto \begin{pmatrix} q & 1+q \\ 0 & -1 \end{pmatrix}.
\]
Then, this defines an indecomposable representation of 
$\H_3(q,Q)$. Hence, the result follows. 
To prove that $\Ext_B^1(D^{\lambda_1},D^{\lambda_{1,0}})=0$, 
we consider 
\[
T_0 \mapsto \begin{pmatrix} 1 & \alpha \\ 0 & q \end{pmatrix},\;\;
T_1 \mapsto \begin{pmatrix} -1 & \beta \\ 0 & -1 \end{pmatrix},\;\;
T_2 \mapsto \begin{pmatrix} -1 & \gamma \\ 0 & -1 \end{pmatrix},
\]
and require that they satisfy the defining relations. 
Since $T_i-q$, for $i=1,2$, are invertible in this case, 
$\beta=0$ and $\gamma=0$ follow. Then, as $q\ne 1$, we can diagonalize 
$T_0$. This implies that a short exact sequence 
\[
0 \longrightarrow D^{\lambda_{1,0}} \longrightarrow M 
\longrightarrow D^{\lambda_1} \longrightarrow 0
\]
always splits, proving the result. Applying $\omega$ we also get 
$\Ext_B^1(D^{\lambda_{2,0}},D^{\lambda_2})=0$. 

Similarly, to prove that 
$\Ext_B^1(D^{\lambda_2},D^{\lambda_{1,0}})=0$, 
we consider 
\[
T_0 \mapsto \begin{pmatrix} 1 & \alpha \\ 0 & q \end{pmatrix},\;\;
T_1 \mapsto \begin{pmatrix} -1 & \beta \\ 0 & q \end{pmatrix},\;\;
T_2 \mapsto \begin{pmatrix} -1 & \gamma \\ 0 & q \end{pmatrix}.
\]
We may assume $\beta=0$ because $q\ne -1$. 
Then, by requiring the defining relations, 
we get $\alpha=0$ by $(T_0T_1)^2=(T_1T_0)^2$ and $\gamma=0$ by 
$T_0T_2=T_2T_0$. 

By the same argument, $q\ne\pm1$ implies that self--extensions are all zero. 
\end{proof}

Let $P=P^{\lambda_1}\oplus P^{\lambda_2}\oplus P^{\lambda_{1,0}}
\oplus P^{\lambda_{2,0}}$, and define an $F$--algebra $A$ by 
\[
A=\End_B(P/\Rad^3 P).
\]
Our aim is to determine the algebra structure of $A$. To do this, we consider 
the radical series of each indecomposable projective module. 
Observe that $S^{\lambda_{i,0}}$, for 
$i=1,2$, is a submodule of $P^{\lambda_i}$ respectively. 
Thus, the radical series of $S^{\lambda_{1,0}}$ and 
$S^{\lambda_{2,0}}$ are as follows. 
\begin{equation*}
\begin{split}
&D^{\lambda_{1,0}}\\
&D^{\lambda_{2,0}}\\
&D^{\lambda_1}
\end{split}
\qquad\text{and}\qquad
\begin{split}
&D^{\lambda_{2,0}}\\
&D^{\lambda_1}\\
&D^{\lambda_2}
\end{split}
\end{equation*}

We start with $P^{\lambda_1}$. 
By Lemma \ref{ext}, $\Rad P^{\lambda_1}/\Rad^2 P^{\lambda_1}=
D^{\lambda_2}\oplus D^{\lambda_{2,0}}$. As 
$\Soc^3 P^{\lambda_1}/\Soc^2 P^{\lambda_1}$ contains 
$D^{\lambda_{1,0}}$, $D^{\lambda_{1,0}}$ appears in 
$\Rad^2 P^{\lambda_1}/\Rad^3 P^{\lambda_1}$. By the Specht filtration, 
$P^{\lambda_1}$ has a submodule $W$ with 
\[
0\longrightarrow S^{\lambda_{1,0}}\longrightarrow W 
\longrightarrow S^{\lambda_{2,0}}\longrightarrow 0.
\]
We denote the pullback of $\Soc S^{\lambda_{2,0}}$ to $W$ by 
$\Soc S^{\lambda_{2,0}}+S^{\lambda_{1,0}}$ and define a $B$--module $V$ by 
\[
V=(\Rad P^{\lambda_1})/(\Soc S^{\lambda_{2,0}}+S^{\lambda_{1,0}}).
\]
Then we have a short exact sequence 
\[
0\;\;\longrightarrow \quad\begin{split}
&D^{\lambda_{2,0}}\\
&D^{\lambda_1}\end{split}\quad
\longrightarrow\;\; V \;\;\longrightarrow\quad 
\begin{split}
&D^{\lambda_2}\\
&D^{\lambda_1}
\end{split}\quad
\longrightarrow\;\; 0
\]
and $V/\Rad V=D^{\lambda_2}\oplus D^{\lambda_{2,0}}$. Thus, 
$\Rad V=D^{\lambda_1}\oplus D^{\lambda_1}$ by 
$\Ext_B^1(D^{\lambda_1},D^{\lambda_1})=0$ and 
$D^{\lambda_1}\oplus D^{\lambda_1}$ appears in 
$\Rad^2 P^{\lambda_1}/\Rad^3 P^{\lambda_1}$. Taking 
the Specht filtration into consideration, 
these imply that 
the radical series of $P^{\lambda_1}/\Rad^3 P^{\lambda_1}$ 
is as follows. 
\begin{equation*}
\begin{split}
&D^{\lambda_1}\\
D^{\lambda_2}&\oplus D^{\lambda_{2,0}}\\
D^{\lambda_1}\;\oplus\; &D^{\lambda_1}\oplus D^{\lambda_{1,0}}
\end{split}
\end{equation*}

Next we consider $P^{\lambda_2}$. By Lemma \ref{ext}, 
$\Rad P^{\lambda_2}/\Rad^2 P^{\lambda_2}=
D^{\lambda_1}$. Thus, by taking 
the Specht filtration into consideration again,
$P^{\lambda_2}/\Rad^2 P^{\lambda_2}=S^{\lambda_2}$ and 
$\Rad^2 P^{\lambda_2}$ has a submodule $F$ such that 
\[
(\Rad^2 P^{\lambda_2})/F= S^{\lambda_3}\simeq D^{\lambda_2},\quad
F=S^{\lambda_{2,0}}. 
\]
If the radical length of $\Rad^2 P^{\lambda_2}$ was $4$ then 
$P^{\lambda_2}$ would be a uniserial $B$--module and 
$\Rad^3 P^{\lambda_2}=S^{\lambda_{2,0}}$. 
This contradicts to the self--duality of $P^{\lambda_2}$. 
Thus, the radical length of $\Rad^2 P^{\lambda_2}$ is $3$. 
If $(\Rad^2 P^{\lambda_2})/F\simeq D^{\lambda_2}$ 
appeared in $\Rad^3 P^{\lambda_2}$ then 
$D^{\lambda_2}$ would appear as a submodule of 
$\Soc^2 P^{\lambda_2}/\Soc P^{\lambda_2}$, which contradicts to 
$\Ext_B^1(D^{\lambda_2},D^{\lambda_2})=0$. 
Thus, $D^{\lambda_2}$ appears in 
$\Rad^2 P^{\lambda_2}/\Rad^3 P^{\lambda_2}$. 
On the other hand, as $\Top F=D^{\lambda_{2,0}}$ appears in 
$\Soc^3 P^{\lambda_2}/\Soc^2 P^{\lambda_2}$, 
$D^{\lambda_{2,0}}$ appears in $\Rad^2 P^{\lambda_2}/\Rad^3 P^{\lambda_2}$. 
As a conclusion, 
the radical series of $P^{\lambda_2}/\Rad^3 P^{\lambda_2}$ 
is as follows. 
\begin{equation*}
\begin{split}
&D^{\lambda_2}\\
&D^{\lambda_1}\\
D^{\lambda_2}&\oplus D^{\lambda_{2,0}}
\end{split}
\end{equation*}

Applying $\omega$ to the radical structure of 
$P^{\lambda_1}/\Rad^3 P^{\lambda_1}$,  
we know that the radical series of 
$P^{\lambda_{2,0}}/\Rad^3 P^{\lambda_{2,0}}$ 
is as follows. 
\begin{equation*}
\begin{split}
&D^{\lambda_{2,0}}\\
D^{\lambda_1}&\oplus D^{\lambda_{1,0}}\\
D^{\lambda_2}\;\oplus\; &D^{\lambda_{2,0}}\oplus D^{\lambda_{2,0}}
\end{split}
\end{equation*}

Similarly, the radical series of $P^{\lambda_{2,0}}/\Rad^3 P^{\lambda_{2,0}}$ 
is as follows. 
\begin{equation*}
\begin{split}
&D^{\lambda_{1,0}}\\
&D^{\lambda_{2,0}}\\
D^{\lambda_1}&\oplus D^{\lambda_{1,0}}
\end{split}
\end{equation*}

Set $D_1=D^{\lambda_2}$, $D_2=D^{\lambda_1}$, 
$D_3=D^{\lambda_{2,0}}$ and $D_4=D^{\lambda_{1,0}}$. Then, 
the Gabriel quiver of the block algebra $B$ is 
the directed graph $Q$ in Lemma \ref{criterion for wildness-2}. 
We have a surjective algebra homomorphism 
\[
FQ/J^3 \longrightarrow A^{opp},
\]
where $J^3$ is the ideal of $FQ$ generated by paths of length $3$. 
Let us compare the dimensions of $FQ/J^3$ and $A$. 
The computation 
above shows $\operatorname{dim}_F A=20$. Hence, 
$\operatorname{dim}_F FQ/J^3=20$ implies that 
$FQ/J^3\simeq A^{opp}$. Therefore, 
we have determined the algebra structure of $A$ 
and Lemma \ref{criterion for wildness-2} 
implies that $A$ is wild. Since $A$ is a factor algebra of 
$\End_B(P)$ and $\End_B(P)^{opp}$ 
is Morita--equivalent to $B$, 
we have shown that $\H_n(q,Q)$ with $n=e=3, f=1$ is wild. 

Next we consider the case where $e=4$ and $f=2$. 
We have, by Proposition \ref{f>0}, that 
\begin{align*}
[P^{\lambda_1}]&=
[S^{\lambda_1}]+[S^{\lambda_2}]+[S^{\lambda_{2,1}}]+[S^{\lambda_{1,1}}], 
\\
[P^{\lambda_2}]&=
[S^{\lambda_2}]+[S^{\lambda_3}]+[S^{\lambda_{2,1}}], 
\\
[P^{\lambda_3}]&=
[S^{\lambda_3}]+[S^{\lambda_4}]+[S^{\lambda_{2,1}}]
+[S^{\lambda_{2,0}}], 
\end{align*}
and
\begin{align*}
[S^{\lambda_{1,1}}]&=
[D^{\lambda_{1,1}}]+[D^{\lambda_{2,1}}]+[D^{\lambda_1}], 
\\
[S^{\lambda_{2,1}}]&=
[D^{\lambda_{2,1}}]+[D^{\lambda_3}]+[D^{\lambda_2}]
+[D^{\lambda_1}], 
\\
[S^{\lambda_{2,0}}]&=
[D^{\lambda_{2,0}}]+[D^{\lambda_{2,1}}]+[D^{\lambda_3}].
\end{align*}

We argue as in (case 1a). Then the radical series of $S^{\lambda_{2,1}}$ has 
the following form. 
\begin{equation*}
\begin{split}
& D^{\lambda_{2,1}}\\
D^{\lambda_1}&\oplus D^{\lambda_3}\\
& D^{\lambda_2}
\end{split}
\end{equation*}

As $S^{\lambda_{2,0}}$ is a submodule of $P^{\lambda_3}$, 
$\Soc S^{\lambda_{2,0}}=D^{\lambda_3}$. Thus, $S^{\lambda_{2,0}}$ 
is uniserial of the following form. 
\begin{equation*}
\begin{split}
&D^{\lambda_{2,0}}\\
&D^{\lambda_{2,1}}\\
&D^{\lambda_3}
\end{split}
\end{equation*}
As a conclusion, the Gabriel quiver contains the following 
quiver as a subquiver. 

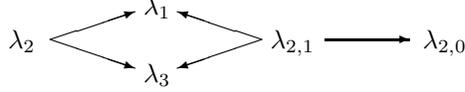
\setlength{\unitlength}{16pt}
\begin{picture}(15,4)(-3,-1)
\put(2,0.8){$\lambda_2$}
\put(3,1){\vector(3,1){2}}
\put(5.2,1.6){$\lambda_1$}
\put(3,1){\vector(3,-1){2}}
\put(5.2,0){$\lambda_3$}
\put(8,1){\vector(-3,1){2}}
\put(8,1){\vector(-3,-1){2}}
\put(8.2,0.8){$\lambda_{2,1}$}
\put(9.5,1){\vector(1,0){2}}
\put(11.8,0.8){$\lambda_{2,0}$}
\end{picture}

\noindent
Thus, Lemma \ref{criterion for wildness-1} implies that 
$\H_n(q,Q)$ with $n=e=4, f=2$ is wild. 

\bigskip
Now, suppose that we are in (case 3a). We consider the case 
where $e\ge 4$ first. 
Then, we have simple modules 
$D^{\lambda_i}$, for $i=1,2,3$. Considering $S^{\lambda_2}$ and 
$S^{\lambda_3}$, we obtain that 
$\Ext_B^1(D^{\lambda_i},D^{\lambda_j})\ne 0$, 
for $\{i,j\}=\{1,2\}$ and $\{2,3\}$. In \cite{AM2}, 
we introduced modules $M^{\lambda_i}$ after Proposition 4.24, 
and showed that 
the self--extensions of the simple modules $D^{\lambda_i}$ 
are non--zero. Thus, 
$\Ext_B^1(D^{\lambda_i},D^{\lambda_i})\ne 0$, 
for $i=2,3$. Therefore, 
we can apply Lemma \ref{criterion for wildness-3}, proving that 
$\H_n(q,Q)$ is wild in this case. 

Next suppose that $e=3$. We compute the decompostion numbers 
in the columns labelled by $\lambda_1$ and $\lambda_2$ as before. 
Then, $D^{\lambda_1}=S^{\lambda_1}$ and $D^{\lambda_2}=S^{\lambda_3}$. 
So, the Specht module theory implies that 
simple $B$--modules $D^{\lambda_1}$ and $D^{\lambda_2}$ 
are given by

\begin{alignat*}{2}
D^{\lambda_1}:\;& T_0\mapsto 1, \quad& T_1,\, T_2 \mapsto -1, \\
D^{\lambda_2}:\;& T_0\mapsto 1, \quad& T_1,\, T_2 \mapsto q. \\
\end{alignat*}

Note that $T_0\mapsto 2-T_0$, $T_i\mapsto q-1-T_i$, for $i=1,2$, 
is an algebra automorphism of $\H_n(q,Q)$ in this case. 
We denote this automorphism by $\omega$ again. Then, $\omega$ 
interchanges $D^{\lambda_1}$ and $D^{\lambda_2}$. Let us consider 
the following representation. 


\begin{gather*}
T_0\mapsto 
\begin{pmatrix}
1 & 0 & -q^2 & 0 & 1 \\
0 & 1 & 0 & 0 & 0 \\
0 & 0 & 1 & 0 & 0 \\
0 & 0 & 0 & 1 & 1 \\
0 & 0 & 0 & 0 & 1 \end{pmatrix}, \quad
T_1\mapsto 
\begin{pmatrix}
q & 0 & 0 & q^2 & 0 \\
0 & -1 & 0 & 0 & 0 \\
0 & 0 & q & 0 & 0 \\
0 & 0 & 0 & -1 & 0 \\
0 & 0 & 0 & 0 & -1 \end{pmatrix}, \\[5pt]
T_2\mapsto 
\begin{pmatrix}
q & 0 & 0 & 0 & 0 \\
0 & -1 & -q^2 & 0 & 1 \\
0 & 0 & q & 0 & 1 \\
0 & 0 & 0 & -1 & 0 \\
0 & 0 & 0 & 0 & -1 \end{pmatrix}.
\end{gather*}

We denote the corresponding $B$--module by $M$. 
The underlying space is $F^5$, the column vectors of dimension $5$. 
If we write $F^4$, this means the subspace consisting of 
the vectors of $F^5$ whose final entry is $0$. 

\begin{lem}
$\Top M=D^{\lambda_1}\;\;\text{and}\;\;\Rad^2 M=\Soc M
=D^{\lambda_1}\oplus D^{\lambda_2}$. 
\end{lem}
\begin{proof}
We can compute $\Top M$ and $\Soc M$ explicitly. 
In the computation of $\Top M=D^{\lambda_1}$, 
we obtain $\Rad M=F^4$. We can also show that 
$\Rad^2 M=\Soc M$ by computing $\Top(\Rad M)$ explicitly. 
\end{proof}

\begin{rem}
$\Rad M$ is indecomposable because 
\[
\End_B(\Rad M)=
FI_4\oplus FE_{13}\oplus FE_{24}\subset M(4,4,F),
\]
is a local $F$--algebra. Here, $I_4$ is the identity matrix 
and $E_{ij}$ are matrix units. 
\end{rem}

As $\Top M=D^{\lambda_1}$, $M$ is a factor module of $P^{\lambda_1}$. 
By twisting $M$ by $\omega$, we know that 
$M^\omega$ is a factor module of $P^{\lambda_2}$. 
Define 
\[
A=\End_B(M\oplus M^\omega).
\]
Then $A$ is a factor algebra of 
$\End_B(P^{\lambda_1}\oplus P^{\lambda_2})$. 

We shall show that $A$ is wild. As $B$ is Morita--equivalent to 
$\End_B(P^{\lambda_1}\oplus P^{\lambda_2})^{opp}$, this 
implies that $B$ is wild. Define
\[
\alpha\in\End_B(M)\subset M(5,5,F)\;\;\text{and}\;\;
\beta\in\End_B(M^\omega)\subset M(5,5,F)
\]
by the matrix $-q^2E_{13}+E_{15}+E_{45}$, and define 
\[
\mu\in\Hom_B(M,M^\omega)\subset M(5,5,F)\;\;\text{and}\;\;
\nu\in\Hom_B(M^\omega,M)\subset M(5,5,F)
\]
by the matrix $E_{14}+E_{23}+qE_{35}$. The embeddings into $M(5,5,F)$ 
are given by the natural basis of $F^5$ which is the underlying space 
of $M$ and $M^\omega$. That they commute with $B$ 
is proved by explicit computations. Finally, 
define $e_1$ and $e_2$ to be the projectors to $M$ and $M^\omega$ 
respectively. 

\begin{prop}
\begin{itemize}
\item[(1)]
$A$ has the basis 
\[
\{e_1,\; e_2,\; \alpha,\; \beta,\; \mu,\; 
\nu,\; \mu\nu,\; \nu\mu,\; \mu\alpha,\; \nu\beta\}.
\]
\item[(2)]
The following relations hold in $A$. 
\begin{gather*}
\alpha^2=0,\;\;\beta^2=0,\;\;\alpha\nu=-\nu\beta,\;\;
\beta\mu=-\mu\alpha,\\
\nu\beta\mu=0,\;\;\mu\alpha\nu=0,\;\;
\mu\nu\mu=0,\;\;\nu\mu\nu=0.
\end{gather*}
\item[(3)]
Let $Q$ be the directed graph in Lemma \ref{criterion for wildness-4} 
then $A\simeq FQ/I$ where the relations are as in (2). 
\end{itemize}
\end{prop}
\begin{proof}
These are proved by explicit computations. 
\end{proof}

Let us consider the factor algebra $\overline A$ 
of $A$ which is defined by requiring one more relation $\nu\beta=0$. 
Then, $\overline A$ has the defining relations 
\begin{gather*}
\alpha^2=0,\;\;\beta^2=0,\;\;\alpha\nu=0,\;\;
\nu\beta=0,\\
\beta\mu=-\mu\alpha,\;\;
\mu\nu\mu=0,\;\;\nu\mu\nu=0.
\end{gather*}
By changing the sign of $\alpha$, we know that 
$\overline A\simeq \Lambda_{32}$. 
Hence, $A$ is wild by Lemma \ref{criterion for wildness-4}. 
We have proved that $\H_n(q,Q)$ is wild in this case. 

\subsection{}
The remaining two cases are for $n=2f+4$. 

Assume that we are in (case 4a). 
Following \cite{AM2}, we define 
bipartitions $\lambda_1=((0),(2^{f+2}))$,
$\lambda_2=((1),(2^{f+1},1))$, $\lambda_3=((1^2),(2^{f+1}))$,
$\lambda_4=((2),(2^f,1^2))$ and $\lambda_5=((2,1),(2^f,1))$. 
They belong to the same block, which we denote by $B$. 

\begin{lem}[{\cite[Lemma 5.6]{AM2}}]
The bipartitions $\lambda_1,\lambda_2,\lambda_3,\lambda_4$
and $\lambda_5$ are all Kleshchev. Furthermore, 
$[P^{\lambda_2}]=[S^{\lambda_2}]+[S^{\lambda_3}]
                            +[S^{\lambda_4}]+[S^{\lambda_5}]$
and the first five rows of the columns of 
the decomposition matrix of $B$ labelled by 
$\lambda_i$, for $1\le i\le 5$, are as
follows. 
$$\begin{array}{l|*5c}
& \lambda_1 & \lambda_2 &\lambda_3 & \lambda_4 & \lambda_5
\\\hline
\lambda_1 & 1 & 0 & 0 & 0 & 0\\
\lambda_2 & 1 & 1 & 0 & 0 & 0\\
\lambda_3 & 0 & 1 & 1 & 0 & 0\\
\lambda_4 & 0 & 1 & 0 & 1 & 0\\
\lambda_5 & 1 & 1 & 1 & 1 & 1
\end{array}$$
\end{lem}

Note that $\Soc S^{\lambda_5}=D^{\lambda_2}$. 
Our aim is to show that 
the radical series of $S^{\lambda_5}$ has the following form. 
\begin{equation*}
\begin{split}
&D^{\lambda_5}\\
D^{\lambda_1}\oplus\; &D^{\lambda_3}\oplus D^{\lambda_4}\\
&D^{\lambda_2}
\end{split}
\end{equation*}
Once this is proved, then the Gabriel quiver contains 
the following quiver as a subquiver. 

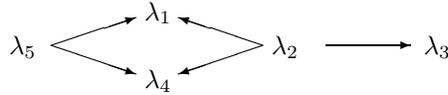
\setlength{\unitlength}{16pt}
\begin{picture}(15,3.5)(-3,-1)
\put(2,0.8){$\lambda_5$}
\put(3,1){\vector(3,1){2}}
\put(5.2,1.6){$\lambda_1$}
\put(3,1){\vector(3,-1){2}}
\put(5.2,0){$\lambda_4$}
\put(8,1){\vector(-3,1){2}}
\put(8,1){\vector(-3,-1){2}}
\put(8.2,0.8){$\lambda_2$}
\put(9.5,1){\vector(1,0){2}}
\put(11.8,0.8){$\lambda_3$}
\end{picture}

\noindent
So, Lemma \ref{criterion for wildness-1} implies that 
$\H_n(q,Q)$ is wild in this case. 

Now we determine the radical structure of $S^{\lambda_5}$. 
In (case 2) of the proof of \cite[Theorem 5.25]{AM2}, we showed 
\[
\Rad P^{\lambda_2}/\Rad^2 P^{\lambda_2}=
D^{\lambda_1}\oplus D^{\lambda_3}\oplus D^{\lambda_4}.
\]
This implies that 
\[
\Soc^2 S^{\lambda_5}/\Soc S^{\lambda_5}\subset 
\Soc^2 P^{\lambda_2}/\Soc P^{\lambda_2}=
D^{\lambda_1}\oplus D^{\lambda_3}\oplus D^{\lambda_4}.
\]
If the inclusion is strict, then one of 
$D^{\lambda_1}$, $D^{\lambda_3}$, $D^{\lambda_4}$ appears 
in $S^{\lambda_2}$, $S^{\lambda_3}$ or $S^{\lambda_4}$. 
Note that these Specht modules are uniserial of the following form. 
\begin{equation*}
S^{\lambda_2}=\;\;
\begin{split}
&D^{\lambda_2}\\
&D^{\lambda_1}
\end{split}\;,\quad
S^{\lambda_3}=\;\;
\begin{split}
&D^{\lambda_3}\\
&D^{\lambda_2}
\end{split}\;,\quad
S^{\lambda_4}=\;\;
\begin{split}
&D^{\lambda_4}\\
&D^{\lambda_2}
\end{split}\;.
\end{equation*}

Assume that $D^{\lambda_1}$ does not appear in 
$\Soc^2 S^{\lambda_5}/\Soc S^{\lambda_5}$. Then 
the $D^{\lambda_1}$ which appears in 
$\Soc^2 P^{\lambda_2}/\Soc P^{\lambda_2}$ must come from 
$\Soc S^{\lambda_2}$. So the $D^{\lambda_1}$ also appears in 
$\Rad P^{\lambda_2}/\Rad^2 P^{\lambda_2}$. This implies that 
the heart $H(P^{\lambda_2})$ of $P^{\lambda_2}$ has the form 
\[
H(P^{\lambda_2})=D^{\lambda_1}\oplus M,
\]
where $M$ is some $B$--module. Observe that 
\begin{itemize}
\item[(a)]
$M$ is self--dual because $H(P^{\lambda_2})$ is self--dual,
\item[(b)]
$\Top M=D^{\lambda_3}\oplus D^{\lambda_4}$,
\item[(c)]
$M$ has a filtration whose succesive quotients are 
$S^{\lambda_3}$, $S^{\lambda_4}$ and 
$S^{\lambda_5}/\Soc S^{\lambda_5}$.
\end{itemize}
As $\Soc M=D^{\lambda_3}\oplus D^{\lambda_4}$, both 
$\Soc S^{\lambda_3}=D^{\lambda_2}$ and $\Soc S^{\lambda_4}=D^{\lambda_2}$ 
cannot appear in $\Soc M$. Thus, 
$\Soc(S^{\lambda_5}/\Soc S^{\lambda_5})=D^{\lambda_3}\oplus D^{\lambda_4}$. 
As $S^{\lambda_5}$ has the unique maximal submodule $\Rad S^{\lambda_5}$, 
the radical series of $S^{\lambda_5}/\Soc S^{\lambda_5}$ has the 
following form. 
\begin{equation*}
\begin{split}
&D^{\lambda_5}\\
&D^{\lambda_1}\\
D^{\lambda_3}&\oplus D^{\lambda_4}
\end{split}
\end{equation*}
In particular, $H(M)$ contains a uniserial module of length $2$ whose top is 
$D^{\lambda_5}$ and whose socle is $D^{\lambda_1}$. As $H(M)$ is 
also self--dual, $H(M)$ must have a uniserial module of length $2$ 
whose top is $D^{\lambda_1}$ and whose socle is $D^{\lambda_5}$ 
as a factor module. However, this is impossible because 
$[H(M):D^{\lambda_1}]=1$ and $[H(M):D^{\lambda_5}]=1$. 

Next assume that $D^{\lambda_3}$ does not appear in 
$\Soc^2 S^{\lambda_5}/\Soc S^{\lambda_5}$. Then the $D^{\lambda_3}$ 
which appears in $\Soc^2 P^{\lambda_2}/\Soc P^{\lambda_2}$ 
must come from $S^{\lambda_3}$. Then, 
$S^{\lambda_3}\subset \Soc^2 P^{\lambda_2}$. Thus, both 
$\Soc S^{\lambda_3}$ and $\Soc S^{\lambda_5}$ appear in 
$\Soc P^{\lambda_2}$, that is, 
$D^{\lambda_2}\oplus D^{\lambda_2}\subset \Soc P^{\lambda_2}$. 
This is a contradiction. 

Assume that $D^{\lambda_4}$ does not appear in 
$\Soc^2 S^{\lambda_5}/\Soc S^{\lambda_5}$. Then, by 
the similar reason as in the previous case, we reach 
a contradiction. 

$\Soc^2 S^{\lambda_5}/\Soc S^{\lambda_5}=
D^{\lambda_1}\oplus D^{\lambda_3}\oplus D^{\lambda_4}$ 
implies that the radical series of $S^{\lambda_5}$ 
is in the desired form, since $\Top S^{\lambda_5}=D^{\lambda_5}$ and 
$\Soc S^{\lambda_5}=D^{\lambda_2}$. 

\bigskip
Assume that we are in (case 5a). In \cite[{section 6}]{AM2}, we 
introduced the notion of path sequences for bipartitions. 
Let $A=\binom{0}{1}, B=\binom{1}{0}, C=\binom{0}{1}, D=\binom{1}{1}$ 
and define $x_r$ and $x_l$ for $x\in\{a, b, c, d\}$ 
as in [{{\it loc. cit.}}]. 
As $f=0$, we do not have to define $x_m$ for $x\in\{a, b, c, d\}$. 
By \cite[{(6.28),(6.29)}]{AM2} we have 
\[
a_l+a_r=b_l+b_r,\quad n>(a_l+a_r)^2.
\]
Our assumption $n<9$ implies that the number of $A$'s is at most $2$. 
On the other hand, as in (case 1) in the proof 
of \cite[Theorem 6.30]{AM2}, $n<e$ implies that Specht modules which belongs 
to a block must have the same set of contents. Further, if 
$S^\lambda$ and $S^\mu$ belong to the block then the 
path sequence of $\mu$ is obtained from that of $\lambda$ 
by a sequence of interchanging $A$ and $B$. 

\begin{defn}
For a bipartition $\lambda=(\lambda^{(1)},\lambda^{(2)})$, 
we denote $(\lambda^{(2)},\lambda^{(1)})$ by $\lambda^\sharp$. 
\end{defn}

Note that $\lambda^\sharp$ has the same set of contents as $\lambda$ 
because of $f=0$. 

\begin{lem}
\label{a=0,1}
Assume that $4\le n<\operatorname{min}\{e, 9\}$. 
If $B$ is a block algebra of $\H_n(q,-1)$ then 
the number of Specht modules belonging to $B$ is $1$, $2$ or $6$. 
If the number is $1$ then $B$ is 
semisimple. If the number is $2$ then 
$B$ is Morita--equivalent to $F[X]/(X^2)$. 
\end{lem}
\begin{proof}
Choose a Specht module $S^\lambda$ which belongs to $B$ so that 
$\lambda^\sharp\triangleleft\lambda$ does not hold, and 
let $a=a_l+a_r$ be the number of $A$'s in the path 
sequence of $\lambda$. 

If $a=0$ then we cannot interchange 
$A$ and $B$, which implies that $S^\lambda$ is the unique Specht module which 
belongs to $B$. Thus, $P^\lambda=S^\lambda=D^\lambda$ and 
$B$ is semisimple. 

If $a=1$ then the number of Specht modules which belongs to 
$B$ is at most $2$. If it is $1$ then $B$ is semisimple. 
If it is $2$, then $S^{\lambda^\sharp}$ belongs to $B$. Observe that 
\[
\operatorname{dim}_F B=2(\operatorname{dim}_F S^\lambda)^2, \;\;
P^\lambda\ne S^\lambda=D^\lambda.
\]
Thus we can write $[P^\lambda]=[S^\lambda]+m[S^{\lambda^\sharp}]$, 
for some $m\ge 1$, and 
\[
(m+1)(\operatorname{dim}_F S^\lambda)^2=
\operatorname{dim}_F P^\lambda\cdot\operatorname{dim}_F D^\lambda
\le \operatorname{dim}_F B=2(\operatorname{dim}_F S^\lambda)^2
\]
implies that $m=1$. Hence, $S^{\lambda^\sharp}=D^{\lambda}$ by 
dimension counting, and we have that $D^\lambda$ is the unique 
simple $B$--module and that $P^\lambda$ is uniserial of length $2$. 
Hence $B$ is Morita--equivalent to $F[X]/(X^2)$. 

If $a\ge 2$ then $a=2$ and the number of Specht modules which 
belongs to $B$ is at most $6$. However, in Proposition \ref{a=2} below, 
we will list all of these possibilities and compute 
the decomposition numbers explicitly. 
From the decomposition matrix we know that these six Specht modules 
always consititute a single block. 
\end{proof}

To cover the case where $a=a_l+a_r=2$, we shall prove the following.

\begin{prop}
\label{a=2}
Assume that $4\le n<\operatorname{min}\{e, 9\}$. 
If $B$ is a block algebra of $\H_n(q,-1)$ such that the 
number of Specht modules belonging to $B$ is not $1$ nor $2$, 
then there are Kleshchev bipartitions $\lambda_1$ and $\lambda_2$, 
and a bipartition $\lambda_3$ which is not Kleshchev, 
such that the transpose of the 
decomposition matrix of $B$ has the following form.

\begin{center}
\begin{tabular}{c|cccccc}
 & $\lambda_1$ & $\lambda_2$ & $\lambda_3$ & $\lambda_3^\sharp$ 
& $\lambda_2^\sharp$ & $\lambda_1^\sharp$ \\\hline
$\lambda_1$      & 1 & 1 & 0 & 0 & 1 & 1 \\
$\lambda_2$      & 0 & 1 & 1 & 1 & 1 & 0 \\
\end{tabular}
\end{center}

Further, there are bipartitions $\mu_1, \mu_2$ of smaller ranks 
such that, for $i=1,\;2$, the following hold. 
\begin{itemize}
\item[(i)]
$D^{\mu_i}=S^{\mu_i}$, $P^{\mu_i}$ is uniserial of length $2$. 
\item[(ii)]
$P^{\lambda_i}=P^{\mu_i}\!\uparrow^B$.
\item[(iii)]
If we set $M_i=D^{\mu_i}\!\uparrow^B$ then
\begin{equation*}
\begin{split}
[M_1]&=[S^{\lambda_1}]+[S^{\lambda_2}]=2[D^{\lambda_1}]+[D^{\lambda_2}],\\
[M_2]&=[S^{\lambda_2}]+[S^{\lambda_3}]=[D^{\lambda_1}]+2[D^{\lambda_2}].
\end{split}
\end{equation*}
\item[(iv)]
Let $B_i$ be the block algebra which $D^{\mu_i}$ belongs, for $i=1,2$. 
Then
\begin{equation*}
\begin{split}
D^{\lambda_1}\!\downarrow_{B_1}=D^{\mu_1},\;\;&
D^{\lambda_1}\!\downarrow_{B_2}=0,\\
D^{\lambda_2}\!\downarrow_{B_2}=D^{\mu_2},\;\;&
D^{\lambda_2}\!\downarrow_{B_1}=0.
\end{split}
\end{equation*}
\end{itemize}
\end{prop}
\begin{proof}
It is enough to verify the decomposition matrix and 
the condition (iii) and (iv) for those Kleshchev bipartitions 
$\mu_1, \mu_2$ with $a_l+a_r=1$; 
if $a=1$ for $\mu_i$ and $\mu_i$ is Kleshchev 
then Lemma \ref{a=0,1} implies (i). Then 
$[P^{\mu_i}\!\uparrow^B]=2[M_i]=[P^{\lambda_i}]$ implies (ii). 

We start with the cases where no 
$S^{\lambda}$ with $\lambda^{(1)}=\emptyset$ belongs to $B$. 
Then, by explicit computation, we know that $n=7$ or $n=8$. 
The possibilities for $n=7$ are
\[
\lambda_1=((1),(3^2)),\;\lambda_2=((2),(3,2)),\;\lambda_3=((3),(2^2))
\]
and
\[
\lambda_1=((1),(2^3)),\;\lambda_2=((1^2),(2^2,1)),\;\lambda_3=((1^3),(2^2)).
\]
In the first case, 
\begin{equation*}
\begin{split}
G(\lambda_1)&=f_1f_2f_0f_1f_{-1}f_0^{(2)}((0),(0))=
\lambda_1+v\lambda_2+v\lambda_2^\sharp+v^2\lambda_1^\sharp,\\
G(\lambda_2)&=f_2f_0f_{-1}f_1^{(2)}f_0^{(2)}((0),(0))=
\lambda_2+v\lambda_3+v\lambda_3^\sharp+v^2\lambda_2^\sharp,
\end{split}
\end{equation*}
and we can check that $\lambda_1, \lambda_2$ are Kleshchev and 
the other four bipartitions are not. Useful criteria for this 
are
\begin{itemize}
\item[(x1)]
Assume that $f=0$. Then $\lambda$ with $\lambda^{(1)}\ne\emptyset$ 
and $\lambda^{(2)}=\emptyset$ is not Kleshchev. 
\item[(x2)]
Assume that $f=0$ and $e>n$ and that both $\lambda^{(1)}$ and 
$\lambda^{(2)}$ are hooks. Then $\lambda$ is Kleshchev if and 
only if $\lambda^{(1)}\subset \lambda^{(2)}$. 
\item[(x3)]
Assume that $f=0$ and $e>n=s+4$. Then 
$((2^2),(s))$ and $((2^2),(1^s))$, for $s\ge0$, are not Kleshchev, 
\item[(x4)]
Assume that $f=0$ and $e>n=s+4$. Then 
$((s),(2^2))$ and $((1^s),(2^2))$, for $s\ge3$, are not Kleshchev. 
\end{itemize}
Thus, the decomposition matrix follows. 

If we set $\mu_1=((1),(3,2))$ and $\mu_2=((2),(2^2))$ then 
we can verify the conditions for $\mu_i$, for $i=1,2$. Note that 
$[S^{\lambda_1}]=[D^{\lambda_1}]$ and 
$[S^{\lambda_2}]=[D^{\lambda_1}]+[D^{\lambda_2}]$.

In the second case, 
\begin{equation*}
\begin{split}
G(\lambda_1)&=f_{-1}f_0f_1f_{-2}f_{-1}f_0^{(2)}((0),(0))=
\lambda_1+v\lambda_2+v\lambda_2^\sharp+v^2\lambda_1^\sharp,\\
G(\lambda_2)&=f_0f_1f_{-2}f_1^{(2)}f_0^{(2)}((0),(0))=
\lambda_2+v\lambda_3+v\lambda_3^\sharp+v^2\lambda_2^\sharp,
\end{split}
\end{equation*}
and we can check that $\lambda_1, \lambda_2$ are Kleshchev and 
the other four bipartitions are not. We set 
$\mu_1=((1),(2^2,1))$ and $\mu_2=((1^2),(2^2))$. 

Next, the possibilities for $n=8$ are six cases given in the table below. 
We always have 
\begin{equation*}
\begin{split}
G(\lambda_1)&=
\lambda_1+v\lambda_2+v\lambda_2^\sharp+v^2\lambda_1^\sharp,\\
G(\lambda_2)&=
\lambda_2+v\lambda_3+v\lambda_3^\sharp+v^2\lambda_2^\sharp.
\end{split}
\end{equation*}
The verification of the conditions for $\mu_i$, for $i=1,\;2$, 
given in the table, are left to the reader. 

\begin{center}
\begin{tabular}{ccc|cc}
$\lambda_1$ & $\lambda_2$ & $\lambda_3$ & $\mu_1$ & $\mu_2$ \\ \hline
$((1),(4,3))$ & $((2),(4,2))$ & $((2^2),(4))$ & 
$((1),(4,2))$ & $((2),(4))$\\
$((1),(3^2,1))$ & $((2),(3,2,1))$ & $((3),(2^2,1))$ &
$((1),(3,2,1))$ & $((2),(2^2,1))$\\
$((1),(3,2^2))$ & $((1^2),(3,2,1))$ & $((1^3),(3,2))$ & 
$((1),(3,2,1))$ & $((1^2),(3,2))$\\
$((1),(2^3,1))$ & $((1^2),(2^2,1^2))$ & $((1^4),(2^2))$ & 
$((1),(2^2,1^2))$ & $((1^2),(2^2))$\\
$((2),(2^3))$ & $((2,1),(2^2,1))$ & $((2,1^2),(2^2))$ &
$((2),(2^2,1))$ & $((2,1),(2^2))$\\
$((1^2),(3^2))$ & $((2,1),(3,2))$ & $((2^2),(3,1))$ &
$((1^2),(3,2))$ & $((2,1),(3,1))$
\end{tabular}
\end{center}

We turn to the cases where $S^\lambda$ with $\lambda^{(1)}=\emptyset$ 
belongs to $B$. For $n=4$, there is only one possibility 
$\lambda^{(2)}=(2^2)$. For $n=5$, $\lambda^{(2)}=(3,2)$ or $(2^2,1)$. 
For $n=6$, $\lambda^{(2)}$ is one of
$(4,2),\;(3^2),\;(3,2,1),\;(2^3),\;(2^2,1^2)$.
For $n=7$, there are $8$ possibilities: $\lambda^{(2)}$ is one of 
$(5,2),\;(4,3),\;(4,2,1),\;(3^2,1)$ and their transposes. 
For $n=8$, there are $14$ possibilities: $\lambda^{(2)}$ is one of 
\[
(6,2),\;(5,3),\;(5,2,1),\;(4^2),\;(4,3,1),\;(4,2^2),\;
(4,2,1^2),\;(3^2,2)\]
and their transposes. Note that the last two are symmetric partitions. 

In these cases, $\lambda^{(2)}$ is of the form
\[
(j+1,i+2,2^{k-1},1^{l-k})\;\;(0\le i<j,\;1\le k\le l)
\]
That is, $\lambda^{(2)}$ is obtained from $(j+1,1^l)$ by adding 
the hook $(i+1,1^{k-1})$. Let $\lambda_1$ be the bipartition 
$\lambda=(\emptyset,\lambda^{(2)})$ and define $\lambda_2$ and $\lambda_3$ by
\[
\lambda_2=((i+1,1^{k-1}),(j+1,1^l)),\;\;
\lambda_3=((j+1,1^{k-1}),(i+1,1^l)).
\]
Then, 
\begin{equation*}
\begin{split}
G(\lambda_1)&=f_{-k+1}\cdots f_{-1}f_i\cdots f_0
f_{-l}\cdots f_{-1}f_j\cdots f_0((0),(0))\\
&=\lambda_1+v\lambda_2+v\lambda_2^\sharp+v^2\lambda_1^\sharp,\\
G(\lambda_2)&=f_{-l}\cdots f_{-k}f_j\cdots f_{i+1}
f_{-k+1}^{(2)}\cdots f_{-1}^{(2)}f_i^{(2)}\cdots f_0^{(2)}((0),(0))\\
&=\lambda_2+v\lambda_3+v\lambda_3^\sharp+v^2\lambda_2^\sharp.
\end{split}
\end{equation*}
We can check that the other four bipartitions are not Kleshchev. 

We define $\mu_1$ to be the 
bipartition obtained from $\lambda_1$ by deleting the hook 
$(i+1,1^{k-1})$, and $\mu_2$ to be the 
bipartition obtained from $\lambda_2$ by deleting 
$(j-i)$ nodes from the first row. 
It is easy to check the conditions for $\mu_i$. 
\end{proof}

Now, we are in a position to determine the radical structure of 
$P^{\lambda_1}$ and $P^{\lambda_2}$ for all of the blocks with 
six Specht modules. Let $B$ be such a block algebra of 
$\H_n(q,-1)$. By Proposition \ref{a=2}, $P^{\lambda_i}$ has 
a submodule which is isomorphic to $M_i$, which implies that 
$\Soc M_i=D^{\lambda_i}$. As $M_i$ is self--dual and 
$[M_i]=2[D^{\lambda_i}]+[D^{\lambda_{3-i}}]$, $M_i$ is uniserial 
of length $3$, whose top and bottom are $D^{\lambda_i}$. 

Let $n_i=|\mu_i|$, for $i=1, 2$. 
As $\H_n(q,-1)$ is projective as a right $\H_{n_i}(q,-1)$--module, 
the Eckmann--Shapiro lemma \cite[Corollary 2.8.4]{B2} implies that 
\[
\Ext^1_B(M_i,D^{\lambda_j})\simeq 
\Ext^1_{\H_{n_i}(q,-1)}(D^{\mu_i},D^{\lambda_j}\!\downarrow_{\H_{n_i}(q,-1)})
=\Ext^1_{B_i}(D^{\mu_i},D^{\lambda_j}\!\downarrow_{B_i}).
\]
Thus, $\operatorname{dim}_F \Ext^1_B(M_i,D^{\lambda_j})=\delta_{ij}$. 
Since $M_i$ and $D^{\lambda_i}$ is self--dual, this also implies that 
$\operatorname{dim}_F \Ext^1_B(D^{\lambda_i},M_j)=\delta_{ij}$. As 
\begin{multline*}
0=\Hom_B(D^{\lambda_2},M_1)\longrightarrow
\Hom_B(D^{\lambda_2},M_1/\Soc M_1)\\
\longrightarrow \Ext^1_B(D^{\lambda_2},D^{\lambda_1})
\longrightarrow \Ext^1_B(D^{\lambda_2},M_1)=0,
\end{multline*}
we conclude that $\Ext^1_B(D^{\lambda_1},D^{\lambda_2})\simeq 
\Ext^1_B(D^{\lambda_2},D^{\lambda_1})=F$. 

Nextly, we consider the long exact sequences 
\begin{multline*}
0=\Hom_B(D^{\lambda_i},M_i/\Soc M_i)\longrightarrow
\Ext^1_B(D^{\lambda_i},D^{\lambda_i})\\
\longrightarrow \Ext^1_B(D^{\lambda_i},M_i)
\stackrel{f_i}{\longrightarrow} 
\Ext^1_B(D^{\lambda_i},M_i/\Soc M_i),
\end{multline*}
for $i=1,2$. We show that $f_i=0$. Assume to the contrary that $f_i\ne 0$. 
Then, $\Ext_B^1(D^{\lambda_i},M_i)=F$ implies that 
$\Ext_B^1(D^{\lambda_i},D^{\lambda_i})=0$. Thus, 
$H(P^{\lambda_i})$ has the simple head and the simple socle, 
and the heart of $H(P^{\lambda_i})$ has composition factors 
$2[D^{\lambda_i}]$. As $\Ext_B^1(D^{\lambda_i},D^{\lambda_i})=0$, 
the heart of $H(P^{\lambda_i})$ must be semisimple. Hence, 
$\Rad P^{\lambda_i}/\Soc^2 P^{\lambda_i}$ has the 
following radical structure. 
\begin{equation*}
\begin{split}
&D^{\lambda_j}\\
D^{\lambda_i}&\oplus D^{\lambda_i}
\end{split}
\end{equation*}
where $j=3-i$. However, this contradicts to 
$\Ext_B^1(D^{\lambda_1},D^{\lambda_2})=F$. We have proved 
that $f_i=0$. This implies that 
$\Ext_B^1(D^{\lambda_i},D^{\lambda_i})=F$, for $i=1,2$. 
In particular, the bipartite graph associated with the Gabriel quiver of $B$ 
is not a Dynkin diagram and Theorem \ref{stable equivalence} implies that 
$B$ is not finite. 

We shall show 
\begin{equation*}
H(P^{\lambda_1})=
D^{\lambda_1}\bigoplus\;
\begin{split}
&D^{\lambda_2}\\
&D^{\lambda_1}\\
&D^{\lambda_2}
\end{split}\;,\quad
H(P^{\lambda_2})=
D^{\lambda_2}\bigoplus\; 
\begin{split}
&D^{\lambda_1}\\
&D^{\lambda_2}\\
&D^{\lambda_1}
\end{split}\;.
\end{equation*}

First suppose that $H(P^{\lambda_1})$ has radical length $2$. 
$H(P^{\lambda_1})$ has Specht filtration whose successive 
quotients are $S^{\lambda_2}$ and $S^{\lambda_2^\sharp}$. 
As $S^{\lambda_2^\sharp}$ is a submodule of $P^{\lambda_2}$, 
the decomposition matrix tells us that $S^{\lambda_2^\sharp}$ is 
uniserial of length $2$ with $\Top S^{\lambda_2^\sharp}=D^{\lambda_1}$ and 
$\Soc S^{\lambda_2^\sharp}=D^{\lambda_2}$, which is 
the dual of $S^{\lambda_2}$. Our assumption implies that 
\[
H(P^{\lambda_1})=S^{\lambda_2}\oplus S^{\lambda_2^\sharp}.
\]

Let $\mu\in\Hom_B(P^{\lambda_1},P^{\lambda_2})$ be such that 
$\Im \mu$ is equal to 
$D^{\lambda_1}\subset\Top(\Rad P^{\lambda_2})$ modulo 
$\Rad^2 P^{\lambda_2}$. 
$\Im \mu$ has the simple socle $\Soc P^{\lambda_2}$. Thus, 
$\mu$ factors through $P^{\lambda_1}/\Soc P^{\lambda_1}$ and 
$D^{\lambda_1}=\Soc S^{\lambda_2}\subset 
\Soc(P^{\lambda_1}/\Soc P^{\lambda_1})$ must vanish. Therefore, 
$\Im \mu$ is uniserial of the following form. 
\begin{equation*}
\Im \mu\simeq\; 
\begin{split}
&D^{\lambda_1}\\
&D^{\lambda_1}\\
&D^{\lambda_2}
\end{split}
\end{equation*}
As $H(P^{\lambda_2})$ contains the submodule $\Im \mu/\Soc P^{\lambda_2}$ and 
$H(P^{\lambda_2})$ is self--dual, this implies 
\begin{equation*}
0 \longrightarrow 
\begin{split}
&D^{\lambda_2}\\
&D^{\lambda_2}
\end{split} \longrightarrow H(P^{\lambda_2}) \longrightarrow 
\begin{split}
&D^{\lambda_1}\\
&D^{\lambda_1}
\end{split} \longrightarrow 0
\end{equation*}
because $\Ext^1_B(D^{\lambda_2},D^{\lambda_2})=F$ implies that 
$[\Soc(H(P^{\lambda_2})):D^{\lambda_2}]=1$. Since we have the submodule 
$\Im \mu/\Soc P^{\lambda_2}$, this exact sequence splits. This leads to 
the following contradiction. 
\begin{equation*}
M_2/\Soc M_2\simeq\; 
\begin{split}
&D^{\lambda_2}\\
&D^{\lambda_1}
\end{split}\; \subset H(P^{\lambda_2})\simeq\; 
\begin{split}
&D^{\lambda_1}\\
&D^{\lambda_1}
\end{split}\; \bigoplus\; 
\begin{split}
&D^{\lambda_2}\\
&D^{\lambda_2}
\end{split}
\end{equation*}
Hence, the radical length of $H(P^{\lambda_1})$ is greater than $2$. 

Recall that $\Rad P^{\lambda_1}/\Rad^2 P^{\lambda_1}
=D^{\lambda_1}\oplus D^{\lambda_2}$. 
Further, $D^{\lambda_1}$ must appear in 
$\Rad^2 P^{\lambda_1}/\Rad^3 P^{\lambda_1}$ 
because $S^{\lambda_2}$ is a factor module of $\Rad P^{\lambda_1}$. Hence 
the radical series of $P^{\lambda_1}$ is as follows. 
\begin{equation*}
\begin{split}
&D^{\lambda_1}\\
D^{\lambda_1} &\oplus D^{\lambda_2}\\
&D^{\lambda_1}\\
&D^{\lambda_2}\\
&D^{\lambda_1}
\end{split}
\end{equation*}

Let $\nu\in\Hom_B(P^{\lambda_2},P^{\lambda_1})$ be such that 
$\Im \nu$ is equal to 
$D^{\lambda_2}\subset\Top(\Rad P^{\lambda_1})$ modulo 
$\Rad^2 P^{\lambda_1}$. We also choose $\alpha\in\End_B(P^{\lambda_1})$ in 
such a way that $\Im \alpha$ is equal to 
$D^{\lambda_1}\subset\Top(\Rad P^{\lambda_1})$ modulo 
$\Rad^2 P^{\lambda_1}$. 

If $\Rad(\Im \alpha)$ contained $\Rad^2 P^{\lambda_1}$ then 
$\Rad^2(\Im \alpha)/\Rad^3(\Im \alpha)=D^{\lambda_2}$ would appear 
in $\Rad^2 P^{\lambda_1}/\Rad^3 P^{\lambda_1}$, a contradiction. 
Thus $\Rad(\Im \nu)$ contains $\Rad^2 P^{\lambda_1}$: 
\begin{equation*}
\Im \nu \;\simeq\; 
\begin{split}
&D^{\lambda_2}\\
&D^{\lambda_1}\\
&D^{\lambda_2}\\
&D^{\lambda_1}
\end{split}
\end{equation*}
We have the following exact sequence.
\begin{equation*}
0 \longrightarrow 
\Ker \nu=\begin{split} &D^{\lambda_2}\\&D^{\lambda_2}\end{split} 
\longrightarrow \Rad P^{\lambda_2} \longrightarrow 
\Rad(\Im \nu)=
\begin{split}&D^{\lambda_1}\\&D^{\lambda_2}\\&D^{\lambda_1}\end{split} 
\longrightarrow 0
\end{equation*}
As $\Top(\Rad P^{\lambda_2})=D^{\lambda_1}\oplus D^{\lambda_2}$, this implies 
that $H(P^{\lambda_2})=D^{\lambda_2}\oplus N_2$, for some submodule $N_2$. 
Now $N_2\simeq \Rad P^{\lambda_2}/\Ker \nu$ implies that 
$H(P^{\lambda_2})$ has the desired form. In particular, 
the radical series of $P^{\lambda_2}$ is as follows. 
\begin{equation*}
\begin{split}
&D^{\lambda_2}\\
D^{\lambda_1} &\oplus D^{\lambda_2}\\
&D^{\lambda_2}\\
&D^{\lambda_1}\\
&D^{\lambda_2}
\end{split}
\end{equation*}

Let $\beta\in\End_B(P^{\lambda_2})$ be 
such that $\Im \beta$ is equal to 
$D^{\lambda_2}\subset\Top(\Rad P^{\lambda_2})$ modulo 
$\Rad^2 P^{\lambda_2}$. More precisely, we choose $\beta$ as follows. 
\[
\Im \beta/\Soc P^{\lambda_2}=D^{\lambda_2}\subset 
D^{\lambda_2}\oplus N_2=H(P^{\lambda_2}).
\]
As $\Rad(\Im \beta)=\Rad^4 P^{\lambda_2}$, 
$\Rad(\Im \mu)$ contains $\Rad^2 P^{\lambda_2}$: 
\begin{equation*}
\Im \mu \;\simeq\; 
\begin{split}
&D^{\lambda_1}\\
&D^{\lambda_2}\\
&D^{\lambda_1}\\
&D^{\lambda_2}
\end{split}
\end{equation*}
This implies that we can choose $N_2=\Im \mu/\Soc P^{\lambda_2}$ and 
\[
H(P^{\lambda_2})=\Im \beta/\Soc P^{\lambda_2}\oplus 
\Im \mu/\Soc P^{\lambda_2}. 
\]

If we consider 
\begin{equation*}
0 \longrightarrow 
\Ker \mu=\begin{split} &D^{\lambda_1}\\&D^{\lambda_1}\end{split} 
\longrightarrow \Rad P^{\lambda_1} \longrightarrow 
\Rad(\Im \mu)=
\begin{split}&D^{\lambda_2}\\&D^{\lambda_1}\\&D^{\lambda_2}\end{split} 
\longrightarrow 0
\end{equation*}
and argue in the same way as above, we can write 
$H(P^{\lambda_1})=D^{\lambda_1}\oplus N_1$, for some $N_1$, and then 
we can choose $\alpha$ so that 
$\Im \alpha/\Soc P^{\lambda_1}=D^{\lambda_1}$ and 
$N_1=\Im \nu/\Soc P^{\lambda_1}$. Thus, 
$H(P^{\lambda_1})$ has the desired form and 
\[
H(P^{\lambda_1})=\Im \alpha/\Soc P^{\lambda_2}\oplus 
\Im \nu/\Soc P^{\lambda_1}.
\]

It is now straightforward to check that 
$\End_B(P^{\lambda_1}\oplus P^{\lambda_2})$ is generated by 
$\alpha, \beta, \mu, \nu$ subject to the relations 
\begin{gather*}
\mu\alpha=0,\;\alpha\nu=0,\;\beta\mu=0,\;\nu\beta=0,\\[5pt]
\alpha^2=(\nu\mu)^2,\;\beta^2=(\mu\nu)^2.
\end{gather*}
The basis is given by 
\[
\{e_1,\;e_2,\;\alpha,\;\beta,\;\mu,\;\nu,\;\alpha^2,\;\beta^2,\;
\mu\nu,\;\nu\mu,\;\mu\nu\mu,\;\nu\mu\nu\}.
\]
This is a special biserial algebra and we have already proved that 
it is not finite. Thus, Theorem \ref{special biserial} implies that 
$B$ is tame, and Lemma \ref{complexity=2(2)} says that there is an 
indecomposable $B$-module with complexity $2$. 

\bigskip
Finally, assume that we are in (case 6a). 
Define $\lambda_i$, for $1\le i\le 10$, by
\begin{gather*}
\lambda_1=((0),(3^3)),\;\lambda_2=((1),(3^2,2)),\;
\lambda_3=((1^2),(3^2,1)),\;
\lambda_4=((2),(3,2^2)),\\
\lambda_5=((1^3),(3^2)),\;\lambda_6=((2,1),(3,2,1)),\;\lambda_7=((3),(2^3)),\\
\lambda_8=((2,1^2),(3,2)),\;
\lambda_9=((2^2),(3,1^2)),\;\lambda_{10}=((3,1),(2^2,1)).
\end{gather*}
Then, $\lambda_1,\lambda_2,\lambda_3,\lambda_4$ and $\lambda_6$ are 
Kleshchev and the others are not. The canonical basis elements are 
as follows. 
\begin{equation*}
\begin{split}
G(\lambda_1)&=f_0f_1f_2f_{-1}f_0f_1f_{-2}f_{-1}f_0((0),(0))\\
&=\lambda_1+v\lambda_2+v\lambda_6+v^2\lambda_9
+v\lambda_9^\sharp+v^2\lambda_6^\sharp+v\lambda_1^\sharp\\
G(\lambda_2)&=f_1f_2f_{-1}f_0f_1f_{-2}f_{-1}f_0^{(2)}((0),(0))\\
&=\lambda_2+v\lambda_3+v\lambda_4+v^2\lambda_6
+v^2\lambda_6^\sharp+v^3\lambda_4^\sharp+v^3\lambda_3^\sharp
+v^4\lambda_2^\sharp\\
G(\lambda_3)&=f_1f_2f_0f_{-1}f_{-2}f_{-1}^{(2)}f_0^{(2)}((0),(0))\\
&=\lambda_3+v\lambda_5+v\lambda_6+v^2\lambda_8+
+v\lambda_8^\sharp+v^2\lambda_6^\sharp+v^2\lambda_5^\sharp
+v^3\lambda_3^\sharp\\
G(\lambda_4)&=f_2f_{-1}f_0f_{-2}f_{-1}f_1^{(2)}f_0^{(2)}((0),(0))\\
&=\lambda_4+v\lambda_6+v\lambda_7+v^2\lambda_{10}+
+v\lambda_{10}^\sharp+v^2\lambda_7^\sharp+v^2\lambda_6^\sharp
+v^3\lambda_4^\sharp\\
G(\lambda_6)&=f_2f_0f_{-2}f_1^{(2)}f_{-1}^{(2)}f_0^{(2)}((0),(0))\\
&=\lambda_6+v\lambda_8+v\lambda_9+v\lambda_{10}+
+v^2\lambda_{10}^\sharp+v^2\lambda_9^\sharp+v^2\lambda_8^\sharp
+v^3\lambda_6^\sharp
\end{split}
\end{equation*}

Thus, the bipartitions $\lambda_i$ and $\lambda_i^\sharp$, for $1\le i\le 10$, 
form a block algebra of $\H_9(q,-1)$, which we denote by $B$. 
We have
\begin{equation*}
S^{\lambda_2}=\begin{split}&D^{\lambda_2}\\&D^{\lambda_1}\end{split}\;,\;\;
S^{\lambda_3}=\begin{split}&D^{\lambda_3}\\&D^{\lambda_2}\end{split}\;,\;\;
S^{\lambda_4}=\begin{split}&D^{\lambda_4}\\&D^{\lambda_2}\end{split}\;.
\end{equation*}

Let $\mu_1=((1^2),(3^2))$ and $\mu_2=((2,1),(3,2))$. These bipartitions 
appeared in the last line of the table in (case 5a) as $\lambda_1$ and 
$\lambda_2$. Thus, $D^{\mu_1}$ and $D^{\mu_2}$ belong to a tame block, say 
$B'$, of $\H_8(q,-1)$ and we have the following.
\begin{equation*}
S^{\mu_1}=D^{\mu_1}, \quad
S^{\mu_2}=\begin{split}&D^{\mu_2}\\&D^{\mu_1}\end{split}
\end{equation*}
Hence, by $[S^{\mu_1}\!\uparrow^B]=[S^{\lambda_3}]+[S^{\lambda_5}]$ and 
$[S^{\mu_2}\!\uparrow^B]=[S^{\lambda_6}]+[S^{\lambda_8}]$, we have 
\[
[D^{\mu_2}\!\uparrow^B]=[D^{\lambda_1}]+[D^{\lambda_4}]+2[D^{\lambda_6}].
\]
Now $D^{\lambda_1}=S^{\lambda_1}$ and $D^{\lambda_4}=S^{\lambda_7}$ 
imply $D^{\lambda_1}\!\downarrow_{B'}=0$, $D^{\lambda_4}\!\downarrow_{B'}=0$, 
thus 
\[
\dim\;\Hom_B(D^{\mu_2}\!\uparrow^B, D^{\lambda_i})\le 
\dim\;\Hom_{B'}(D^{\mu_2}, D^{\lambda_i}\!\downarrow_{\H_8(q,-1)})=0,
\]
for $i=1,4$. Hence $D^{\mu_2}\!\uparrow^B$ is not semisimple. 
If the radical length was $2$ then 
$D^{\lambda_1}\oplus D^{\lambda_4}$ would appear in 
$\Soc(D^{\mu_2}\!\uparrow^B)$. However, as $D^{\mu_2}\!\uparrow^B$ 
is self--dual, 
$D^{\lambda_1}\oplus D^{\lambda_4}$ would also appear in 
$\Top(D^{\mu_2}\!\uparrow^B)$, a contradiction. If 
$D^{\lambda_6}$ appeared twice in $\Top(D^{\mu_2}\!\uparrow^B)$, 
$D^{\lambda_6}\oplus D^{\lambda_6}$ 
would be a direct summand of $D^{\mu_2}\!\uparrow^B$. 
Then $D^{\lambda_1}$ or $D^{\lambda_4}$ would appear in 
$\Top(D^{\mu_2}\!\uparrow^B)$, a contradiction again. 
Therefore, 
\[
\Top(D^{\mu_2}\!\uparrow^B)=D^{\lambda_6}=\Soc(D^{\mu_2}\!\uparrow^B)
\]
and the heart of $D^{\mu_2}\!\uparrow^B$ must be 
$D^{\lambda_1}\oplus D^{\lambda_4}$ as it is self--dual. We have proved 
that the radical series of $D^{\mu_2}\!\uparrow^B$ is as follows. 
\begin{equation*}
\begin{split}
&D^{\lambda_6}\\
D^{\lambda_1}&\oplus D^{\lambda_4}\\
&D^{\lambda_6}
\end{split}
\end{equation*}
To summarize, the Gabriel quiver contains the following 
quiver as a subquiver. 

\setlength{\unitlength}{16pt}
\begin{picture}(15,4)(-3,-1)
\put(2,0.8){$\lambda_6$}
\put(3,1){\vector(3,1){2}}
\put(5.2,1.6){$\lambda_1$}
\put(3,1){\vector(3,-1){2}}
\put(5.2,0){$\lambda_4$}
\put(8,1){\vector(-3,1){2}}
\put(8,1){\vector(-3,-1){2}}
\put(8.2,0.8){$\lambda_2$}
\put(9.2,1){\vector(1,0){2}}
\put(11.5,0.8){$\lambda_3$}
\end{picture}

\noindent
Thus, Lemma \ref{criterion for wildness-1} implies that 
$\H_9(q,-1)$ with $e>n=9, f=0$ is wild. 

\subsection{}
In this subsection we prove Theorem \ref{two parameter case}(2). 
Thus, we assume that $e=2$. This implies that $e\le 2f+4$. 
As we also assume that 
$0\le f\le\frac{e}{2}$, we have $f=0$ or $f=1$. 
The cases to consider are as follows.

\begin{itemize}
\item[(case 1b)]
$n=2$ and $f=0$. 
\item[(case 2b)]
$n=2$ and $f=1$. 
\item[(case 3b)]
$n=3$ and $f=0$.
\item[(case 4b)]
$n=3$ and $f=1$.
\item[(case 5b)]
$n=4$ and $f=1$.
\end{itemize}

Our aim is to show that $\H_n(q,Q)$ is tame in (case 1b), (case 2b) and 
(case 4b), and wild in (case 3b) and (case 5b). 
We also prove that, in the tame cases, the tame block algebras of 
$\H_n(q,Q)$ are special biserial algebras. We will also show, 
by using Lemma \ref{complexity=2(1)} and Lemma \ref{complexity=2(2)}, 
that there is an indecomposable $\H_n(q,Q)$--module with complexity $2$ 
in these cases. 

Before going into the case--by--case analysis, we write 
the first three layers of the crystal graphs 
$B(2\Lambda_0)$ and $B(\Lambda_0+\Lambda_1)$. As was explained before, 
the nodes of each of the layers of 
the crystal graphs parametrize simple $\H_n(-1,-1)$--modules 
and simple $\H_n(-1,1)$--modules, for $n=1,2,3$, respectively. 
Further, we can compute 
the full decomposition matrices in these cases. 

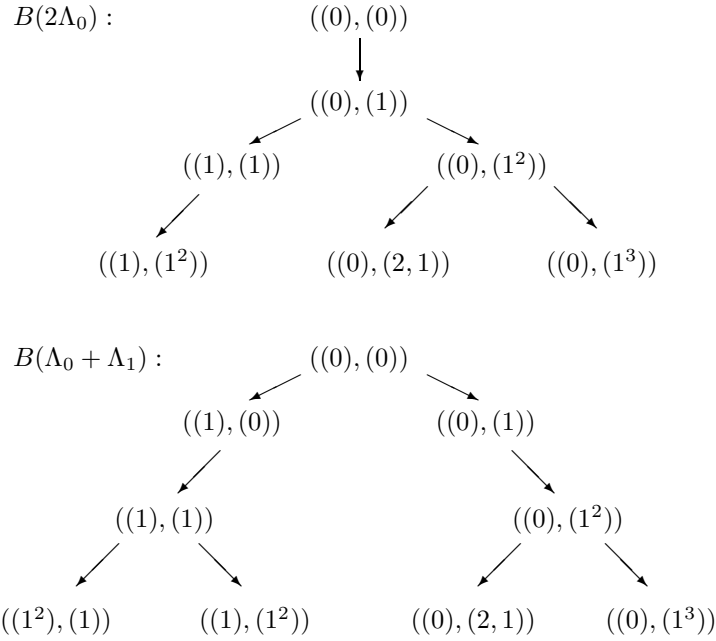
\setlength{\unitlength}{16pt}
\begin{picture}(15,8.5)(-3,0)
\put(-2,7){$B(2\Lambda_0):$}
\put(5,7){$((0),(0))$}
\put(6.2,6.7){\vector(0,-1){1}}
\put(5,5){$((0),(1))$}
\put(4.8,4.8){\vector(-2,-1){1.2}}
\put(2,3.5){$((1),(1))$}
\put(2.4,3){\vector(-1,-1){1}}
\put(0,1.2){$((1),(1^2))$}
\put(7.8,4.8){\vector(2,-1){1.2}}
\put(8,3.5){$((0),(1^2))$}
\put(7.8,3.2){\vector(-1,-1){1}}
\put(10.8,3.2){\vector(1,-1){1}}
\put(5.4,1.2){$((0),(2,1))$}
\put(10.6,1.2){$((0),(1^3))$}
\end{picture}

\setlength{\unitlength}{16pt}
\begin{picture}(15,8)(-3,0)
\put(-2,7){$B(\Lambda_0+\Lambda_1):$}
\put(5,7){$((0),(0))$}
\put(4.8,6.8){\vector(-2,-1){1.2}}
\put(2,5.5){$((1),(0))$}
\put(7.8,6.8){\vector(2,-1){1.2}}
\put(8,5.5){$((0),(1))$}

\put(2.9,5){\vector(-1,-1){1}}
\put(0.4,3.2){$((1),(1))$}
\put(9.8,5){\vector(1,-1){1}}
\put(9.8,3.2){$((0),(1^2))$}

\put(0.5,2.8){\vector(-1,-1){1}}
\put(2.4,2.8){\vector(1,-1){1}}
\put(-2.3,0.8){$((1^2),(1))$}
\put(2.4,0.8){$((1),(1^2))$}

\put(10,2.8){\vector(-1,-1){1}}
\put(12,2.8){\vector(1,-1){1}}
\put(7.4,0.8){$((0),(2,1))$}
\put(12,0.8){$((0),(1^3))$}
\end{picture}

First assume that we are in (case 1b). Then, $\H_2(q,Q)$ has two 
blocks. One is a semisimple block algebra $\End_F(S^{((1),(1))})$. 
We denote the other block algebra by $B$. 
Let $\lambda_1=((0),(1^2))$. 
Then, $B$ has the unique simple module $D^{\lambda_1}$ and 
$[P^{\lambda_1}:D^{\lambda_1}]=4$. If $P^{\lambda_1}$ is uniserial, 
then $B$ is finite by Lemma \ref{basic local}. 
However, $B$ is not finite by 
\cite[Theorem 1.4]{AM2}. Thus, the radical series of $P^{\lambda_1}$ 
has the following form. 
\begin{equation*}
\begin{split}
&D^{\lambda_1}\\
D^{\lambda_1} &\oplus D^{\lambda_1}\\
&D^{\lambda_1}
\end{split}
\end{equation*}
$\End_B(P^{\lambda_1})$ is the Kronecker algebra 
$F[X, Y]/(X^2, Y^2)$, which is a special biserial algebra. 
Thus, $B$ is tame. Applying Lemma \ref{complexity=2(2)}, 
we know that all simple $B$--modules have complexity $2$. 

Next assume that we are in (case 2b). Then $\H_2(q,Q)$ has one block. 
So we write $B$ for $\H_2(q,Q)$. 
We denote $((0),(1^2))$ and $((1),(1))$ by $\lambda_1$, 
$\lambda_2$ respectively. Then, the decomposition matrix of $B$ 
is as follows. 

\medskip
\begin{center}
\begin{tabular}{c|cc}
 & $\lambda_1$ & $\lambda_2$ \\\hline
$((0),(1^2))$      & 1 & 0 \\
$((0),(2))$      & 1 & 0 \\
$((1),(1))$      & 1 & 1 \\
$((1^2),(0))$  & 0 & 1 \\
$((2),(0))$  & 0 & 1 \\
\end{tabular}
\end{center}

\medskip
Thus, $D^{\lambda_1}$ and $D^{\lambda_2}$ are given by
\begin{alignat*}{2}
D^{\lambda_1}:\;& T_0\mapsto -1, \quad& T_1 \mapsto -1, \\
D^{\lambda_2}:\;& T_0\mapsto 1, \quad& T_1 \mapsto -1, \\
\end{alignat*}

Let us consider the following three representations. 
\begin{gather*}
T_0\mapsto 
\begin{pmatrix}
-1 & 0 \\
0 & -1 \end{pmatrix},\quad
T_1\mapsto 
\begin{pmatrix}
-1 & 1 \\
0 & -1 \end{pmatrix},\\[5pt]
T_0\mapsto 
\begin{pmatrix}
1 & 0 \\
0 & -1 \end{pmatrix},\quad
T_1\mapsto 
\begin{pmatrix}
-1 & 1 \\
0 & -1 \end{pmatrix},\\[5pt]
T_0\mapsto 
\begin{pmatrix}
1 & 0 \\
0 & 1 \end{pmatrix},\quad
T_1\mapsto 
\begin{pmatrix}
-1 & 1 \\
0 & -1 \end{pmatrix}.
\end{gather*}
These are indecomposable representations. So, 
$\Ext_B^1(D^{\lambda_i},D^{\lambda_j})\ne0$, 
for all $i$ and for all $j$. Hence, the radical structure of 
$P^{\lambda_1}$ and $P^{\lambda_2}$ are as follows. 

\begin{equation*}
\begin{split}
&D^{\lambda_1}\\
D^{\lambda_1}&\oplus D^{\lambda_2}\\
&D^{\lambda_1}
\end{split}
\qquad
\begin{split}
&D^{\lambda_2}\\
D^{\lambda_1}&\oplus D^{\lambda_2}\\
&D^{\lambda_2}
\end{split}
\end{equation*}

Observe that there is a uniserial submodule $U^i_j$ of $\Rad P^{\lambda_j}$ 
whose top is $D^{\lambda_i}$ and whose socle is 
$\Soc P^{\lambda_j}=D^{\lambda_j}$. 

Define $\alpha\in\End_B(P^{\lambda_1})$ and 
$\beta\in\End_B(P^{\lambda_2})$ by $\Im \alpha=U^1_1$, 
$\Im \beta=U^2_2$. Then, $\Ker \alpha=U^2_1$, $\Ker \beta=U^1_2$. 

Similarly, define $\mu\in\Hom_B(P^{\lambda_1},P^{\lambda_2})$ and 
$\nu\in\Hom_B(P^{\lambda_2},P^{\lambda_1})$ by 
$\Im \mu=U^1_2$, 
$\Im \nu=U^2_1$. Then, $\Ker \mu=U^1_1$, $\Ker \nu=U^2_2$. 

Therefore, $\End_B(P^{\lambda_1}\oplus P^{\lambda_2})\simeq FQ/I$ 
where $Q$ is the directed graph with adjacency matrix 
$\binom{1\;1}{1\;1}$ and the relations are given by
\begin{gather*}
\mu\alpha=0,\;\;\beta\mu=0,\;\;\alpha\nu=0,\;\;\nu\beta=0,\\
\alpha^2=\nu\mu,\;\;\beta^2=\mu\nu,
\end{gather*}
because both $\End_B(P^{\lambda_1}\oplus P^{\lambda_2})$ 
and $FQ/I$ are $8$--dimensional. 
Note that $FQ/I$ has the basis
\[
\{e_1,\;e_2,\;\alpha,\;\beta,\;
\mu,\;\nu,\;\alpha^2,\;\beta^2\}.
\]

Since $FQ/I$ is a special biserial algebra, $B$ is tame or finite. 
However, as $B$ is not finite by \cite[Theorem 1.4]{AM2}, $B$ is tame. 
Applying Lemma \ref{complexity=2(2)} again, we know that 
all simple $B$--modules have complexity $2$. 

Now, we assume that we are in (case 3b). Then, $\H_3(q,Q)$ has 
two blocks. One has $D^{((0),(2,1))}$ as the unique simple module 
and $[P^{((0),(2,1))}:D^{((0),(2,1))}]=2$. Thus, $P^{((0),(2,1))}$ 
is uniserial of length $2$. Then, Lemma \ref{basic local} implies that 
this block algebra is finite. 

We denote the other block algebra by $B$. 
We write $\lambda_1$ and $\lambda_2$ for 
$((0),(1^3))$ and $((1),(1^2))$ respectively. 
Then, by the computation of the canonical basis as before, 
$B$ has the following decomposition matrix. 

\medskip
\begin{center}
\begin{tabular}{c|cc}
 & $\lambda_1$ & $\lambda_2$ \\\hline
$((0),(1^3))$  & 1 & 0 \\
$((0),(3))$    & 1 & 0 \\
$((1),(1^2))$  & 1 & 1 \\
$((1),(2))$    & 1 & 1 \\
$((1^2),(1))$  & 1 & 1 \\
$((2),(1))$    & 1 & 1 \\
$((1^3),(0))$  & 1 & 0 \\
$((3),(0))$    & 1 & 0 \\
\end{tabular}
\end{center}

\medskip
As the dual of $S^{\lambda_2}$ is uniserial of length $2$ 
whose top is $D^{\lambda_1}$ 
and whose socle is $D^{\lambda_2}$, we have 
$\Ext_B^1(D^{\lambda_1},D^{\lambda_2})\ne0$. 
We shall show that 
$\Ext_B^1(D^{\lambda_1},D^{\lambda_1})=F^2$. 
Then, the Gabriel quiver of $B$ contains the directed graph $Q$
with adjacency matrix $\binom{2\;1}{0\;0}$. Thus, 
Lemma \ref{criterion for wildness-5} implies that $B$ is wild. 

The computation of $\Ext_B^1(D^{\lambda_1},D^{\lambda_1})$ 
is easy. In fact, $D^{\lambda_1}$ is given by
\[
T_0\mapsto 1,\;\;T_1\mapsto -1,\;\;T_2\mapsto -1,
\]
and the result is a family of representations 
\[
T_0\mapsto 
\begin{pmatrix}
1 & x \\
0 & 1 \end{pmatrix},\quad
T_1,\;T_2\mapsto
\begin{pmatrix}
-1 & y \\
0 & -1 \end{pmatrix},\quad\text{where $x, y\in F$.}
\]
Denote the corresponding module by $M(x,y)$. 
This is indecomposable if and only if $(x,y)\ne(0,0)$. 
Assume that $M(x,y)$ and $M(x',y')$ are indecomposable. 
Then, $M(x,y)\simeq M(x',y')$ if and only if 
$(x,y)$ and $(x',y')$ define the same element in 
the projective space $P^1(F)$. 

Assume that we are in (case 4b). Then, $\H_3(q,Q)$ has two blocks 
$B$ and $B'$, and if we denote 
$((0),(2,1))$ and $((1),(1^2))$ by $\lambda_1$, $\lambda_2$, and 
$((0),(1^3))$ and $((1^2),(1))$ by $\lambda_1'$, $\lambda_2'$, 
then the decomposition matrices are as follows. 

\medskip
\begin{center}
$B:\;$
\begin{tabular}{c|cc}
 & $\lambda_1$ & $\lambda_2$ \\\hline
$((0),(2,1))$  & 1 & 0 \\
$((1),(1^2))$  & 1 & 1 \\
$((1),(2))$    & 1 & 1 \\
$((1^3),(0))$  & 0 & 1 \\
$((3),(0))$    & 0 & 1 \\
\end{tabular}
\qquad
$B':\;$
\begin{tabular}{c|cc}
 & $\lambda_1'$ & $\lambda_2'$ \\\hline
$((0),(1^3))$  & 1 & 0 \\
$((0),(3))$    & 1 & 0 \\
$((1^2),(1))$  & 1 & 1 \\
$((2),(1))$    & 1 & 1 \\
$((2,1),(0))$  & 0 & 1 \\
\end{tabular}
\end{center}

\medskip
Thus, $D^{\lambda_2}=S^{((1^3),(0))}$ and $D^{\lambda_1'}=S^{((0),(1^3))}$ 
imply that 
\begin{alignat*}{3}
D^{\lambda_2}:\;& T_0\mapsto \;\,\;1, 
\quad& T_1 \mapsto -1, \quad& T_2 \mapsto -1,\\
D^{\lambda_1'}:\;& T_0\mapsto -1, 
\quad& T_1 \mapsto -1, \quad& T_2 \mapsto -1.\\\end{alignat*}

We have an algebra automorphism defined by 
$T_0\mapsto -T_0$ and $T_i\mapsto -2-T_i$, for $i=1,2$. 
We denote this automorphism by $\omega$ again. Then $\omega$ interchanges 
$D^{\lambda_2}$ and $D^{\lambda_1'}$, which implies that 
$\omega$ interchanges $B$ and $B'$. Hence, to show that 
$\H_3(q,Q)$ is tame, it is enough to consider $B$ only. 

As $S^{((1),(2))}$ is a submodule of $P^{\lambda_1}$, 
there is a uniserial module whose top is $D^{\lambda_2}$ and 
whose socle is $D^{\lambda_1}$. Hence, 
$\Ext_B^1(D^{\lambda_1},D^{\lambda_2})\ne0$. 

\begin{lem}
\label{brutal computation}
\begin{itemize}
\item[(1)]
Choosing a suitable basis of $D^{\lambda_1}$, 
$D^{\lambda_1}$ can be represented by the following matrices. 
\[
T_0\mapsto 
\begin{pmatrix}
-1 & 0 \\
0 & -1 \end{pmatrix},\quad
T_1\mapsto 
\begin{pmatrix}
-1 & -1 \\
0 & -1 \end{pmatrix},\quad
T_2\mapsto 
\begin{pmatrix}
0 & -1 \\
1 & -2 \end{pmatrix}.
\]
\item[(2)]
$\Ext_B^1(D^{\lambda_1},D^{\lambda_1})=0$ and 
$\Ext_B^1(D^{\lambda_2},D^{\lambda_2})\ne0$. 
\end{itemize}
\end{lem}
\begin{proof}
As $D^{\lambda_1}=S^{\lambda_1}$, we can obtain the matrix 
representation by reduction from the seminormal representation. 
To obtain the same matrix representation as above, we take 
a nonzero vector $v$ such that 
\[
L_1v=qv,\;L_2v=q^2v,\;L_3v=v,\;
T_1v=qv,
\]
where $L_1=T_0$, $L_2=q^{-1}T_1T_0T_1$ and 
$L_3=q^{-2}T_2T_1T_0T_1T_2$. We choose the basis 
$\{v, T_2v\}$. Then, after some computations, we obtain (1). 
Now we consider (2). As the second assertion is easy to verify, 
we focus on the first assetion. Verification of the details is left 
to the reader. Let $T_i$, for $i=0,1,2$, be the matrices 
given in (1) and define $\hat T_i$ as follows. 
\[
\hat T_0=\begin{pmatrix}
T_0 & X \\ 0 & T_0 \end{pmatrix},\quad
\hat T_1=\begin{pmatrix}
T_1 & Y \\ 0 & T_1 \end{pmatrix},\quad
\hat T_2=\begin{pmatrix}
T_2 & Z \\ 0 & T_2 \end{pmatrix},
\]
where $X, Y, Z\in M(2,2,F)$. 
We require the defining relations and write down 
the relations among $X, Y, Z$. The quadratic relations 
imply that $X$ is the zero matrix and $Y$, $Z$ are given by 
\begin{equation*}
Y=\begin{pmatrix} a & b \\ 0 & -a\end{pmatrix}\quad\text{and}\quad
Z=\begin{pmatrix} c & d \\ 2c+d & -c\end{pmatrix},\quad
\text{where $a,b,c,d\in F$.}
\end{equation*}
Further, the braid relation between $\hat T_1$ and $\hat T_2$ 
implies $2a+b=2c+d$. 
Then, using these equalities, 
we compute the socle of this $4$--dimensional matrix 
representation. Then, the $2$--dimensional subspaces 
\begin{equation*}
F\begin{pmatrix}1\\0\\0\\0\end{pmatrix}\;\bigoplus\;
F\begin{pmatrix}0\\1\\0\\0\end{pmatrix}\quad
\text{and}\quad
F\begin{pmatrix}0\\a\\1\\0\end{pmatrix}\;\bigoplus\;
F\begin{pmatrix}c-a\\b\\0\\1\end{pmatrix}
\end{equation*}
are submodules, both of which are isomorphic to $D^{\lambda_1}$. 
This shows that a short exact sequence 
\[
0\longrightarrow D^{\lambda_1}\longrightarrow M 
\longrightarrow D^{\lambda_1}\longrightarrow 0
\]
always splits. Hence the result. 
\end{proof}

We shall determine the heart $H(P^{\lambda_i})$ of 
$P^{\lambda_i}$, for $i=1,2$. More precisely, 
our goal is to show 
\begin{equation*}
H(P^{\lambda_1})=
\begin{split}
&D^{\lambda_2}\\
&D^{\lambda_1}\\
&D^{\lambda_2}
\end{split}
\quad\text{and}\quad
H(P^{\lambda_2})=
D^{\lambda_2}\bigoplus\;
\begin{split}
&D^{\lambda_1}\\
&D^{\lambda_2}\\
&D^{\lambda_1}
\end{split}\;.
\end{equation*}

We start with $P^{\lambda_1}$. 
Lemma \ref{brutal computation} implies that 
\[
[\Rad^2 P^{\lambda_1}]=2[D^{\lambda_1}]+[D^{\lambda_2}].
\]
As $S^{((1),(2))}$ is a submodule of $P^{\lambda_1}$, 
$\Ext_B^1(D^{\lambda_1},D^{\lambda_2})\ne0$ and 
the radical length of $\Rad^2 P^{\lambda_1}$ is greater than 
or equal to $2$. 
If the radical length of $\Rad^2 P^{\lambda_1}$ is $2$, then 
\[
\Soc^2 P^{\lambda_1}/\Soc P^{\lambda_1}=D^{\lambda_1}\oplus D^{\lambda_2}, 
\]
which implies that $D^{\lambda_1}$ appears in 
$\Rad P^{\lambda_1}/\Rad^2 P^{\lambda_1}$, a contradiction. 
Hence $P^{\lambda_1}$ is uniserial and $H(P^{\lambda_1})$ 
has the desired form. 

Next consider $P^{\lambda_2}$. 
$H(P^{\lambda_2})$ has a Specht filtration whose 
successive quotients are $\Soc S^{\lambda_2}=D^{\lambda_1}$, 
$S^{((1),(2))}$ and $S^{((1^3),(0))}=D^{\lambda_2}$. 
Recalling that $S^{((1),(2))}$ is a uniserial submodule of $P^{\lambda_1}$, 
we know that $H(P^{\lambda_2})$ is not semisimple. 

As $\Ext_B^1(D^{\lambda_1},D^{\lambda_2})\ne 0$, there exists a surjective $B$--module homomorphism which is the 
composition of maps given as follows. 
\[
\Rad P^{\lambda_2}\longrightarrow \Rad P^{\lambda_2}/\Rad^2 P^{\lambda_2}
\longrightarrow D^{\lambda_1}
\]
Hence, we have $\phi:P^{\lambda_1}\longrightarrow \Rad P^{\lambda_2}$ which is 
a lift of this homomorphism. Taking the radical structure of $P^{\lambda_1}$ 
and $[H(P^{\lambda_2}):D^{\lambda_1}]=2$ 
into consideration, its image in $H(P^{\lambda_2})$ is either 
\begin{equation*}
D^{\lambda_1}\quad\text{or}\qquad
\begin{split}
&D^{\lambda_1}\\
&D^{\lambda_2}\\
&D^{\lambda_1}
\end{split}
\end{equation*}
If the image is $D^{\lambda_1}$, then this appears both in the top and in 
the socle of $H(P^{\lambda_2})$. Thus we can write 
\[
H(P^{\lambda_2})=D^{\lambda_1}\oplus M.
\]
As $H(P^{\lambda_2})$ is self--dual and not semisimple, 
$[M]=[D^{\lambda_1}]+2[D^{\lambda_2}]$ implies that $M$ is uniserial 
with $\Top M=D^{\lambda_2}$, $\Soc M=D^{\lambda_2}$ and $H(M)=D^{\lambda_1}$. 
However, as 
$D^{\lambda_2}$ appears in $\Rad^3 P^{\lambda_1}/\Rad^4 P^{\lambda_1}$, 
Landrock's theorem \cite[Theorem 1.7.8]{B2} implies that 
$D^{\lambda_1}$ must appear in $\Rad^3 P^{\lambda_2}/\Rad^4 P^{\lambda_2}$,
a contradiction. We have proved that the image of $\phi$ in $H(P^{\lambda_2})$ 
is uniserial of length $3$. Consider the surjection 
$\psi:P^{\lambda_2}\longrightarrow \Rad P^{\lambda_1}$ again. 
This induces the surjection 
$H(P^{\lambda_2})\longrightarrow \Rad^2 P^{\lambda_1}$. 
If $D^{\lambda_2}$ does not appear in $\Soc H(P^{\lambda_2})$ then 
$\Soc H(P^{\lambda_2})=D^{\lambda_1}$, which must vanish under this 
surjection. This is a contradiction because 
$D^{\lambda_1}$ must appear in the image twice. 
Therefore, $D^{\lambda_2}$ appears in $\Soc H(P^{\lambda_2})$ and 
this $D^{\lambda_2}$ also appears in 
$\Top H(P^{\lambda_2})$. Thus, $H(P^{\lambda_2})$ is the direct sum 
of $D^{\lambda_2}$ and the image of $\phi$ in $H(P^{\lambda_2})$, 
proving that $H(P^{\lambda_2})$ has the desired form. 
As this $D^{\lambda_2}$ must vanish 
under the map $H(P^{\lambda_2})\longrightarrow \Rad^2 P^{\lambda_1}$, 
we have also proved that $\Ker \psi$ is uniserial of length $2$ 
whose top and socle are $D^{\lambda_2}$. 

Denote the submodule $\Im \phi$ of $P^{\lambda_2}$ by $N_1$ and 
the submodule $\Ker \psi$ of $P^{\lambda_2}$ by $N_2$. We define 
$\beta\in\End_B(P^{\lambda_2})$ by $\Im \beta=N_2$. Then 
$\Ker \beta=N_1$. Nextly, we define 
$\mu\in \Hom_B(P^{\lambda_1},P^{\lambda_2})$ and 
$\nu\in \Hom_B(P^{\lambda_2},P^{\lambda_1})$ by $\mu=\phi$, 
$\nu=\psi$. Thus, $\Im \mu=N_1$, $\Im \nu=\Rad P^{\lambda_1}$, 
$\Ker \mu=\Soc P^{\lambda_1}$ and $\Ker \nu=N_2$. 
Let $Q$ be the directed graph with adjacency matrix 
$\binom{0\;1}{1\;1}$ as in Lemma \ref{complexity=2(1)} and define 
$A=FQ/I$ by the relations
\[
\nu\beta=0,\;\;\beta\mu=0,\;\;\beta^2=(\mu\nu)^2.
\]
Note that $A$ has the following basis.
\[
\{e_1,\;e_2,\;\mu,\;\nu,\;\beta,\;\mu\nu,\;
\nu\mu,\;\beta^2,\;\mu\nu\mu,\;
\nu\mu\nu,\;(\nu\mu)^2\}.
\]
Then, by multiplying $\mu$ or $\nu$ with a suitable scalar, 
we obtain the surjective 
algebra homomorphism 
\[
A\longrightarrow \End_B(P^{\lambda_1}\oplus P^{\lambda_2}),
\]
which is an isomorphism because both are $11$ dimensional. 
As $A$ is a special biserial algebra, $B$ is tame. 
Note that $A$ is self--injective, and if we set 
$S=\oplus_{i\in Q_0} \Soc Ae_i$ then $A/S$ is the algebra 
$\Lambda$ in Lemma \ref{complexity=2(1)}. Therefore, 
using the self--injectivity of $A$, we know that 
the simple $A$--module corresponding to the node $2$ has 
complexity $2$. 

Finally, assume that we are in (case 5b). We consider 
the parabolic subalgebra $\H_2(q,Q)\otimes\H_2^A(q)$ of 
$\H_4(q,Q)$. Then, as in the proof of 
Corollary \ref{critical rank}, the Mackey decomposition theorem 
implies that Proposition \ref{reduction to critical rank} 
applies. Hence, it is enough to prove that 
$\H_2(q,Q)\otimes_F\H_2^A(q)$ is wild. In (case 2b) we proved that 
$\H_2(q,Q)$ is a symmetric tame special biserial algebra 
whose Gabriel quiver satisfies the assumptions of 
Lemma \ref{complexity=2(2)}. Thus, there is an 
$\H_2(q,Q)$--module $M$ with complexity $2$. On the other hand, 
as $\H_2^A(q)$ is not semisimple, Theorem \ref{complexity}(1) implies 
that there is an $\H_2^A(q)$--module $N$ whose complexity is greater than 
or equal to $1$. As the complexity of $M\otimes N$ is greater than 
or equal to $3$, Theorem \ref{complexity}(3) implies the result. 

\section{The case of one parameter Hecke algebras}
\subsection{}
Let $(W,S)$ be a finite Coxeter system. The 
{\sf Poincar\'e polynomial} $P_W(x)$ of $(W,S)$ is 
defined by 
\[
P_W(x)=\sum_{w\in W}x^{\ell(w)}
\]
where $\ell(w)$ is the length of $w$. According to 
$W=W(A_{n-1})$, $W(B_n)$ or $W(D_n)$, we denote the 
Poincar\'e polynomial $P_W(x)$ by $P_n^A(x)$, $P_n^B(x)$ and 
$P_n^D(x)$ respectively. 

It is well--known that the semisimplicity of a Hecke algebra 
is governed by its Poincar\'e polynomial \cite{GU}. 
Our aim in the subsequent sections is to show that 
the Poincar\'e polynomial governs other representation types also 
if the Hecke algebra is of classical type. 
In this section, we consider Hecke algebras associated with an 
irreducible Weyl group, and prove the following theorem. 
The general case where the Weyl group is not assumed to be irreducible 
will be considered in the next section. Note that 
the finiteness result, for $\H_n^X(q)$, $X=A,B,D$, was 
proved in \cite{A4}, and for $\H_n(q,Q)$, in \cite{AM3}. 

\begin{thm}
\label{one parameter case}
Let $e$ be the multiplicative order of $q\ne1$. 
\begin{itemize}
\item[(1)]
Assume that $e\ge 3$. Then, for $X=A$, $B$ or $D$, $\H_n^X(q)$ is 
\begin{itemize}
\item[--]
finite if $(x-q)^2$ does not divide $P_n^X(x)$.
\item[--]
wild otherwise.
\end{itemize}
\item[(2)]
Assume that $e=2$. Then, for $X=A$, $B$ or $D$, $\H_n^X(q)$ is 
\begin{itemize}
\item[--]
finite if $(x-q)^2$ does not divide $P_n^X(x)$.
\item[--]
tame if $(x-q)^2$ divides but $(x-q)^3$ does not divide $P_n^X(x)$.
\item[--]
wild otherwise.
\end{itemize}
\end{itemize}
\end{thm}

\subsection{}
First we consider $\H_n^A(q)$. Note that the result in this case is 
nothing but \cite[Proposition 3.7, Theorem 3.8]{U}. 
We give a proof for the sake of comleteness. 
$P_n^A(x)$ is given by 
\[
P_n^A(x)=\prod_{k=1}^n \frac{x^k-1}{x-1}.
\]
If $n\ge 2e$ then $x^e-1$ and $x^{2e}-1$ are divisible by $x-q$. 
If $n<2e$ then $q^k\ne 1$, for $1\le k\le e-1$ and $e+1\le k\le 2e-1$. 
Thus, $(x-q)^2$ does not divide $P_n^A(x)$ if and only if $n<2e$. 

Assume that $e\ge 3$. Then, Theorem \ref{Erdmann-Nakano theorem} implies 
that if $n<2e$ then $\H_n^A(q)$ is finite, and that if $n\ge 2e$ then 
$\H_n^A(q)$ is wild. Hence (1). 

Next assume that $e=2$. Then, 
the characteristic of $F$ is odd, and $P_n^A(x)$ is divisible by 
$(x+1)^3$ if and only if $\left[\frac{n}{2}\right]\ge 3$, that is, 
$n\ge 6$. By Theorem \ref{Erdmann-Nakano theorem} again, we have that 
$(x+1)^3$ divides $P_n^A(x)$ if and only if $\H_n^A(q)$ is wild. 
Hence (2). 

\subsection{}
Next consider $\H_n^B(q)$, which is $\H_n(q,Q)$ with $Q=q$. 
The result in this case is essentially \cite[Theorem 2.1]{A5}. Here, we prove 
this case as a corollary of 
Theorem \ref{separated parameter case} and 
Theorem \ref{two parameter case}. 
$P_n^B(x)$ is given by 
\[
P_n^B(x)=\prod_{k=1}^n \frac{x^{2k}-1}{x-1}.
\]

Assume that the multiplicative order $e$ of $q$ is odd. 
Then $-Q\not\in q^{\mathbb Z}$, and $(x-q)^2$ does not divide 
$P_n^B(x)$ if and only if $n<2e$. Thus, 
Theorem \ref{separated parameter case} implies the result. 

Next assume that $e$ is even. Then $-Q=q^{\frac{e}{2}+1}$ 
and $(x-q)^2$ does not divide 
$P_n^B(x)$ if and only if $n<e$. If $e\ge 4$ then $f\ne0$ also holds. 
If $e=2$ then $f=0$ and the characteristic of $F$ is odd. Further, 
$(x+1)^3$ divides $P_n^B(x)$ if and only if $n\ge 3$. 
Hence, Theorem \ref{two parameter case} implies the result. 

\subsection{}
Finally, we consider $\H_n^D(q)$. Note that $n\ge 4$ and $\H_n^D(q)$ is 
the $F$--algebra defined by generators $S_0$, $S_1,\dots,S_{n-1}$ 
and relations 
\begin{gather*}
(S_i-q)(S_i+1) =0\;\;(0\le i\le n-1),\\
S_0S_2S_0=S_2S_0S_2,\quad S_0S_i=S_iS_0\;\;(i\ne 2),\\
S_{i+1}S_iS_{i+1}=S_iS_{i+1}S_i\;\;(1\le i\le n-2),\\
S_iS_j=S_jS_i\;\;(1\le i<j-1\le n-2).
\end{gather*}

$P_n^D(x)$ is given by
\[
P_n^D(x)=\frac{x^n-1}{x-1}\;\prod_{k=1}^{n-1}\frac{x^{2k}-1}{x-1}.
\]

We embed $\H_n^D(q)$ into $\H_n(q,1)$ by the injective algebra 
homomorphism 
\begin{equation*}
\H_n^D(q) \longrightarrow \H_n(q,1)
\end{equation*}
defined by 
$S_0\mapsto T_0T_1T_0$ and $S_i\mapsto T_i$, 
for $1\le i\le n-1$, and we identify $\H_n^D(q)$ with its image. 
Then 
\begin{equation*}
\H_n(q,1)=\H_n^D(q)\oplus T_0\H_n^D(q)\text{\; and \;} 
T_0\H_n^D(q)=\H_n^D(q)T_0. 
\end{equation*}

First assume that $e$ is odd. If $n<2e$ then 
Theorem \ref{separated parameter case} implies that 
$\H_n(q,1)$ is finite. Thus, so is $\H_n^D(q)$ by 
\cite[Lemma 2.5]{AM2}. If $n=2e$ then $\H_n^A(q)$ is 
wild. Then, by the Mackey decomposition again, 
Proposition \ref{reduction to critical rank} is 
applicable and $\H_n^D(q)$ is wild. Thus, $\H_n^D(q)$, 
for $n\ge 2e$, is also wild by Corollary \ref{critical rank}(2). 

Next assume that $e$ is even. 
In particular, the characteristic of $F$ is odd again. 

\medskip
\noindent
(case $e=2$) 
As $P_n^D(x)$, for $n\ge 4$, is always divisible by $(x+1)^3$, 
our aim is to show that $\H_n^D(q)$, for $n\ge 4$, is wild. 
Recall that $\H_2^A(q)$ is not 
semisimple. Thus, $\H_2^A(q)\otimes \H_2^A(q)\otimes \H_2^A(q)$ is wild 
by Theorem \ref{complexity}. 
Applying Proposition \ref{reduction to critical rank} to 
$\H_2^A(q)\otimes \H_2^A(q)\otimes \H_2^A(q)\subset \H_4^D(q)$, 
we know that $\H_4^D(q)$ is wild. Thus, $\H_n^D(q)$, for $n\ge 4$, 
is wild by Corollary \ref{critical rank}. 

\medskip
\noindent
(case $e\ge 4$) Note that $\H_n(q,1)$ is of the form 
$\H_n(q,-q^f)$ with $f=\frac{e}{2}$. 

If $n<e$ then $x^{2k}-1$ with $k=\frac{e}{2}$ is the only term 
in $P_n^D(x)$ which may be divisible by $x-q$. On the other hand, 
$n<e$ implies that 
$\H_n(q,1)$ is finite by Theorem \ref{two parameter case}. 
Thus, $\H_n^D(q)$ is finite by \cite[Lemma 2.5]{AM2}. 

Next assume that $n=e$. Then both $x^n-1$ and 
$x^{2k}-1$ with $k=\frac{e}{2}$ are divisible by $x-q$, and 
our aim is to show that $\H_n^D(q)$ is wild in this case. 
To apply the Fock space theory, we set 
$\Lambda=\Lambda_0+\Lambda_{\frac{e}{2}}$ and consider 
the Fock space $\mathcal F_v(\Lambda)$. 
By Theorem \ref{Fock space theory}(1), simple $\H_n(q,1)$--modules are 
those $D^\lambda$ with $\lambda\vdash n$ and $\lambda\in B(\Lambda)$. 

Our assumption that $f=\frac{e}{2}$ and $e\ge 4$ imply that 
$f\ge 2$ and $e-f\ge 2$. So, we are in (case 1a) and (case 2a) of 
subsection (3.3). In both cases, 
\[
[S^{\lambda_{1,\frac{e}{2}-1}}]=
[D^{\lambda_{1,\frac{e}{2}-1}}]+[D^{\lambda_{2,\frac{e}{2}-1}}]
+[D^{\lambda_1}]
\]
and $S^{\lambda_{1,\frac{e}{2}-1}}$ is a submodule of $P^{\lambda_1}$. 
Thus, the radical series of $S^{\lambda_{1,\frac{e}{2}-1}}$ has 
the following form. 
\begin{equation*}
\begin{split}
&D^{\lambda_{1,\frac{e}{2}-1}}\\
&D^{\lambda_{2,\frac{e}{2}-1}}\\
&D^{\lambda_1}
\end{split}
\end{equation*}
In particular, we have the following uniserial 
$\H_n(q,1)$--modules. 
\begin{equation*}
\Rad S^{\lambda_{1,\frac{e}{2}-1}}=\;
\begin{split}
&D^{\lambda_{2,\frac{e}{2}-1}}\\
&D^{\lambda_1}
\end{split}\,,\quad
(S^{\lambda_{1,\frac{e}{2}-1}}/\Soc S^{\lambda_{1,\frac{e}{2}-1}})^*
=\;
\begin{split}
&D^{\lambda_{2,\frac{e}{2}-1}}\\
&D^{\lambda_{1,\frac{e}{2}-1}}\\
\end{split}\;\,.
\end{equation*}

Define an algebra automorphism $\tau$ of $\H_n(q,1)$ by 
\begin{equation*}
\tau(T_1)=T_0T_1T_0 \text{\; and \;} 
\tau(T_i)=T_i\text{\; for $i\ne 1$.}
\end{equation*}
Then $\tau^2=1$ and $\tau$ induces the Dynkin automorphism of 
$\H_n^D(q)$ given by $S_i\mapsto S_{1-i}$ for $i=0,1$ and 
$S_i\mapsto S_i$ for $i\ge 2$. 

Define another algebra automorphism $\sigma$ of $\H_n(q,1)$ by 
\begin{equation*}
\sigma(T_0)=-T_0 \text{\; and \;} 
\sigma(T_i)=T_i\text{\; for $i\ne 0$.}
\end{equation*}
Then $\H_n^D(q)$ is the fixed point subalgebra $\H_n(q,1)^\sigma$. 
Now we apply the Clifford theory to the pair 
$\H_n^D(q)$ and $\H_n(q,1)$. Note that $\sigma\tau=\tau\sigma$. 

Let $\lambda\vdash n$ be a Kleshchev bipartition. 
We define another Kleshchev bipartition 
$h(\lambda)$ by $(D^\lambda)^\sigma=D^{h(\lambda)}$. 
Then the Clifford theory tells us the following. Recall that 
we are in the case where the characteristic of $F$ is odd. 

\begin{thm}[{\cite[4.3, 4.4]{Hu1}}]
\label{Clifford theory}
Assume that the characteristic of $F$ is odd. 

\begin{itemize}
\item[(1)]
If $h(\lambda)\ne\lambda$ then $D^\lambda\!\!\downarrow_{\H_n^D(q)}$ 
remains irreducible and $D^\lambda\!\!\downarrow_{\H_n^D(q)}
\simeq D^{h(\lambda)}\!\!\downarrow_{\H_n^D(q)}$. 
\item[(2)]
If $h(\lambda)=\lambda$ then $D^\lambda\!\!\downarrow_{\H_n^D(q)}$ 
is a direct sum of two pairwise non--isomorphic simple 
$\H_n^D(q)$--modules. If we denote them by $D^{\lambda_\pm}$ then 
$(D^{\lambda_+})^\tau=D^{\lambda_-}$. 
\item[(3)]
The simple $\H_n^D(q)$--modules in (1) and (2) constitute 
a complete set of simple $\H_n^D(q)$--modules. 
\end{itemize}
\end{thm}

The following result of Hu enables us to 
compute $h(\lambda)$ explicitly. 

\begin{thm}[{\cite[Theorem 1.5]{Hu2}}]
\label{Hu's result}
Let $\lambda$ be a Kleshchev bipartition. Take 
a path from the empty bipartition $\emptyset$ 
to $\lambda$ in the crystal graph 
$B(\Lambda)$ and write 
\[
\lambda=\tilde f_{i_1}\cdots\tilde f_{i_n}\emptyset. 
\]
Then $h(\lambda)$ is given by 
\[
h(\lambda)=\tilde f_{i_1+\frac{e}{2}}\cdots\tilde f_{i_n+\frac{e}{2}}\emptyset.
\]
in $B(\Lambda)$. 
\end{thm}

If $h(\lambda)=\lambda$ then we denote the projective cover 
of $D^\lambda_+$ and $D^\lambda_-$ by $P^\lambda_+$ and 
$P^\lambda_-$ respectively. 
If $h(\lambda)\ne\lambda$ then we denote 
the projective cover of 
$D^\lambda\!\!\downarrow_{\H_n^D(q)}$ by $\overline P^\lambda$. 

\begin{lem}
\label{PIM in type D}
\begin{itemize}
\item[(1)]
Suppose that $h(\lambda)=\lambda$. Then, 
$P^\lambda\!\!\downarrow_{\H_n^D(q)}=
P^\lambda_+\oplus 
P^\lambda_-$. 
\item[(2)]
Suppose that $h(\lambda)\ne\lambda$. Then, 
$\overline P^\lambda=P^\lambda\!\!\downarrow_{\H_n^D(q)}$.
\end{itemize}
\end{lem}
\begin{proof}
(1) Theorem \ref{Clifford theory}(2) implies that 
$(P^\lambda_+)^\tau=
P^\lambda_-$ and, by 
the Mackey decomposition theorem, 
$P^\lambda\!\!\downarrow_{\H_n^D(q)}$ is a 
projective $\H_n^D(q)$--module. 
As we have surjective $\H_n^D(q)$--module homomorphisms 
\[
P^\lambda\!\!\downarrow_{\H_n^D(q)}
\longrightarrow 
D^\lambda_+,\quad
P^\lambda\!\!\downarrow_{\H_n^D(q)}
\longrightarrow 
D^\lambda_-,
\]
and $D^\lambda_+
\not\simeq D^\lambda_-$, 
$P^\lambda_+\oplus P^\lambda_-$ is 
a direct summand of $P^\lambda\!\!\downarrow_{\H_n^D(q)}$. 
Thus, we have 
\[
\operatorname{dim}_F P^\lambda_+
+\operatorname{dim}_F P^\lambda_-
\le 
\operatorname{dim}_F P^\lambda.
\]
On the other hand, we have a surjective $\H_n(q,1)$--module 
homomorphism
\[
P^\lambda_+\!\!\uparrow^{\H_n(q,1)}
\longrightarrow 
D^\lambda_+\!\!\uparrow^{\H_n(q,1)}=
D^\lambda,
\]
and this implies that $P^\lambda$ is a direct summand 
of $P^\lambda_+\!\!\uparrow^{\H_n(q,1)}$. Thus, 
\[
\operatorname{dim}_F P^\lambda
\le 
2\operatorname{dim}_F P^\lambda_+.
\]
Therefore, 
\[
2\operatorname{dim}_F P^\lambda_+
=\operatorname{dim}_F P^\lambda_+
+\operatorname{dim}_F P^\lambda_-
\le 
\operatorname{dim}_F P^\lambda
\le 
2\operatorname{dim}_F P^\lambda_+,
\]
which proves the result. 

\noindent
(2) As we have a surjective $\H_n(q,1)$--module homomorphism
\[
\overline P^\lambda\!\!\uparrow^{\H_n(q,1)}
\longrightarrow 
(D^\lambda\!\!\downarrow_{\H_n^D(q)})\!\!\uparrow^{\H_n(q,1)}
=D^\lambda\oplus D^{h(\lambda)},
\]
and $\lambda\ne h(\lambda)$, 
$P^\lambda\oplus P^{h(\lambda)}$ is a direct summand 
of $\overline P^\lambda\!\!\uparrow^{\H_n(q,1)}$. As 
$P^\lambda$ and $P^{h(\lambda)}$ have the same dimension, 
this implies that 
\[
\operatorname{dim}_F P^\lambda\le 
\operatorname{dim}_F \overline P^\lambda. 
\]
On the other hand, we have a surjective $\H_n^D(q)$--module 
homomorphism
\[
P^\lambda\!\!\downarrow_{\H_n^D(q)}
\longrightarrow 
D^\lambda\!\!\downarrow_{\H_n^D(q)},
\]
which implies that $\overline P^\lambda$ is a direct summand 
of $P^\lambda\!\!\downarrow_{\H_n^D(q)}$. The result follows. 
\end{proof}

We apply Theorem \ref{Hu's result} in our setting. For our purposes, 
it suffices to 
compute $h(\lambda_1)$, $h(\lambda_2)$, $h(\lambda_{2,\frac{e}{2}-1})$ 
and $h(\lambda_{1,\frac{e}{2}-1})$. Then, since 
\begin{equation*}
\begin{split}
\lambda_1&=\tilde f_{\frac{e}{2}+1}\cdots \tilde f_{e-1}
\tilde f_0\cdots \tilde f_{\frac{e}{2}}\emptyset,\\
\lambda_2&=\tilde f_{\frac{e}{2}+2}\cdots \tilde f_{e-1}
\tilde f_0\tilde f_{\frac{e}{2}+1}
\tilde f_1\cdots \tilde f_{\frac{e}{2}}\emptyset,\\
\lambda_{2,\frac{e}{2}-1}&=\tilde f_1\cdots \tilde f_{\frac{e}{2}-1}
\tilde f_{\frac{e}{2}+1}\tilde f_{\frac{e}{2}}
\tilde f_{\frac{e}{2}+2}\cdots\tilde f_{e-1}\tilde f_0\emptyset,\\
\lambda_{1,\frac{e}{2}-1}&=\tilde f_1\cdots 
\tilde f_{e-1}\tilde f_0\emptyset,
\end{split}
\end{equation*}
we have 
\begin{equation*}
\begin{split}
h(\lambda_1)&=\tilde f_1\cdots \tilde f_{e-1}\tilde f_0\emptyset
=\lambda_{1,\frac{e}{2}-1},\\
h(\lambda_2)&=\tilde f_2\cdots \tilde f_{\frac{e}{2}}\tilde f_1
\tilde f_{\frac{e}{2}+1}\cdots \tilde f_{e-1}\tilde f_0\emptyset
=\lambda_{1,\frac{e}{2}-2},\\
h(\lambda_{2,\frac{e}{2}-1})&=
\tilde f_{\frac{e}{2}+1}\cdots \tilde f_{e-1}
\tilde f_1\tilde f_0
\tilde f_2\cdots \tilde f_{\frac{e}{2}}\emptyset=\lambda_{2,\frac{e}{2}-1},\\
h(\lambda_{1,\frac{e}{2}-1})&=\tilde f_{\frac{e}{2}+1}\cdots 
\tilde f_{e-1}\tilde f_0 \cdots 
\tilde f_{\frac{e}{2}}\emptyset=\lambda_1. 
\end{split}
\end{equation*}

Therefore, Theorem \ref{Hu's result} implies that 
$D^{\lambda_1}\!\!\downarrow_{\H_n^D(q)}$, 
$D^{\lambda_2}\!\!\downarrow_{\H_n^D(q)}$, 
$D^{\lambda_{1,\frac{e}{2}-1}}\!\!\downarrow_{\H_n^D(q)}$ are 
simple $\H_n^D(q)$--modules, and that 
\[
D^{\lambda_{2,\frac{e}{2}-1}}\!\!\downarrow_{\H_n^D(q)}=
D^{\lambda_{2,\frac{e}{2}-1}}_+\oplus D^{\lambda_{2,\frac{e}{2}-1}}_-.
\]

Since $\Rad S^{\lambda_{1,\frac{e}{2}-1}}$ is a uniserial 
$\H_n(q,1)$--module whose top is $D^{\lambda_{2,\frac{e}{2}-1}}$ and 
whose socle is $D^{\lambda_1}$, there is 
a surjective $\H_n^D(q)$--module homomorphism 
\[
P^{\lambda_{2,\frac{e}{2}-1}}\!\!\downarrow_{\H_n^D(q)}
=
P^{\lambda_{2,\frac{e}{2}-1}}_+\oplus 
P^{\lambda_{2,\frac{e}{2}-1}}_-
\longrightarrow 
\Rad S^{\lambda_{1,\frac{e}{2}-1}}\!\!\downarrow_{\H_n^D(q)}.
\]
$\Rad S^{\lambda_{1,\frac{e}{2}-1}}\!\!\downarrow_{\H_n^D(q)}$ 
is not semisimple: if otherwise, then the simple 
$\H_n^D(q)$--module $D^{\lambda_1}\!\!\downarrow_{\H_n^D(q)}$ 
would appear in 
$\Top(\Rad S^{\lambda_{1,\frac{e}{2}-1}}\!\!\downarrow_{\H_n^D(q)})$. 
However, this would imply that the morphism is not surjective. 
Hence, the radical series of 
$\Rad S^{\lambda_{1,\frac{e}{2}-1}}\!\!\downarrow_{\H_n^D(q)}$ 
is one of the following. 
\begin{equation*}
D^{\lambda_{2,\frac{e}{2}-1}}_+\bigoplus\;\;
\begin{split}
&D^{\lambda_{2,\frac{e}{2}-1}}_-\\
&D^{\lambda_1}\!\!\downarrow_{\H_n^D(q)}
\end{split},
\qquad
\begin{gathered}
D^{\lambda_{2,\frac{e}{2}-1}}_+\oplus 
D^{\lambda_{2,\frac{e}{2}-1}}_-\\
D^{\lambda_1}\!\!\downarrow_{\H_n^D(q)}
\end{gathered},
\qquad
\begin{split}
&D^{\lambda_{2,\frac{e}{2}-1}}_+\\
&D^{\lambda_1}\!\!\downarrow_{\H_n^D(q)}
\end{split}
\bigoplus\;\;
D^{\lambda_{2,\frac{e}{2}-1}}_-.
\end{equation*}

In any case, by applying $\tau$ if necessary, we conclude that 
there are uniserial $\H_n^D(q)$--modules of the following form. 
\begin{equation*}
\begin{split}
&D^{\lambda_{2,\frac{e}{2}-1}}_+\\
&D^{\lambda_1}\!\!\downarrow_{\H_n^D(q)}
\end{split}\qquad\text{and}\qquad
\begin{split}
&D^{\lambda_{2,\frac{e}{2}-1}}_-\\
&D^{\lambda_1}\!\!\downarrow_{\H_n^D(q)}
\end{split}
\end{equation*}

If we consider the dual of 
$S^{\lambda_{1,\frac{e}{2}-1}}/\Soc S^{\lambda_{1,\frac{e}{2}-1}}$, 
and arguing in the same way, we prove the existence of 
$\H_n^D(q)$--modules with 
the following radical structure. 
\begin{equation*}
\begin{split}
&D^{\lambda_{2,\frac{e}{2}-1}}_+\\
&D^{\lambda_{1,\frac{e}{2}-1}}\!\!\downarrow_{\H_n^D(q)}
\end{split}\qquad\text{and}\qquad
\begin{split}
&D^{\lambda_{2,\frac{e}{2}-1}}_-\\
&D^{\lambda_{1,\frac{e}{2}-1}}\!\!\downarrow_{\H_n^D(q)}
\end{split}
\end{equation*}

Consider the $\H_n(q,1)$--module $S^{\lambda_2}$. 
$S^{\lambda_2}$ is uniserial of length $2$ 
whose top is $D^{\lambda_2}$ and whose socle is $D^{\lambda_1}$. 
Then, we have a surjective $\H_n^D(q)$--module homomorphism 
\[
\overline P^{\lambda_2}=
P^{\lambda_2}\!\!\downarrow_{\H_n^D(q)}
\longrightarrow S^{\lambda_2}\!\!\downarrow_{\H_n^D(q)}.
\]
By the same reasoning as above, $S^{\lambda_2}\!\!\downarrow_{\H_n^D(q)}$ 
is not semisimple. Thus, the radical series of 
$S^{\lambda_2}\!\!\downarrow_{\H_n^D(q)}$ is as follows. 
\begin{equation*}
\begin{split}
&D^{\lambda_2}\!\!\downarrow_{\H_n^D(q)}\\
&D^{\lambda_1}\!\!\downarrow_{\H_n^D(q)}
\end{split}
\end{equation*}

As a consequence, by considering 
the dual of these simple $\H_n^D(q)$--modules, we know that 
the Gabriel quiver of $\H_n^D(q)$ contains 
the following quiver as a subquiver. 

\setlength{\unitlength}{16pt}
\begin{picture}(15,3.5)(-3,-1)
\put(-1.1,0.8){$\lambda_{1,\frac{e}{2}-1}$}
\put(1,1){\vector(3,1){2}}
\put(3.2,1.6){$(\lambda_{2,\frac{e}{2}-1})_+$}
\put(1,1){\vector(3,-1){2}}
\put(3.2,0){$(\lambda_{2,\frac{e}{2}-1})_-$}
\put(8,1){\vector(-3,1){2}}
\put(8,1){\vector(-3,-1){2}}
\put(8.2,0.8){$\lambda_1$}
\put(9,1){\vector(1,0){2}}
\put(11.3,0.8){$\lambda_2$}
\end{picture}

Therefore, Lemma \ref{criterion for wildness-1} implies that 
$\H_n^D(q)$ is wild when $n=e$. So, 
$\H_n^D(q)$, for $n\ge e$, is wild by Corollary \ref{critical rank}. 

\section{The general case}

We begin by analysing the case where $\H_n^X(q)$, for $X=A, B, D$, 
is finite but not semisimple. 

\begin{lem}
\label{finite case e=2}
Suppose that $e=2$. If a non--semisimple block algebra of 
$\H_n^A(q)$, $\H_n^B(q)$ or $\H_n^D(q)$ is finite, then it is 
Morita--equivalent to $F[X]/(X^2)$. 
\end{lem}
\begin{proof}
If $\H_n^A(q)$ is finite but not semisimple, then 
$n=2$ or $n=3$ and non--semisimple blocks are 
Morita--equivalent to the algebra $F[X]/(X^2)$. 
If $\H_n(q,Q)$ is finite but not semisimple, 
then Theorem \ref{separated parameter case}(2) and 
Theorem \ref{two parameter case}(2) imply that 
the same holds. In particular, we have the result 
for $\H_n^B(q)$. $\H_n^D(q)$ cannot be finite. 
Thus, the result follows. 
\end{proof}

\begin{lem}
\label{finite case results}
Suppose that $e\ge 3$. 
\begin{itemize}
\item[(1)]
If $\H_n^A(q)$ is finite but not semisimple 
and $B$ is a non--semisimple block algebra of $\H_n^A(q)$, 
then 
\begin{itemize}
\item[--]
the number of pairwise non--isomorphic simple $B$--modules 
is $e-1$,
\item[--]
$B$ is a Brauer tree algebra with the tree being a 
straight line such that all indecomposable projective 
$B$--modules have radical length $3$. 
\end{itemize}
\item[(2)]
If $\H_n(q,Q)$ is finite but not semisimple 
and $B$ is a non--semisimple block algebra of $\H_n(q,Q)$, 
then 
\begin{itemize}
\item[--]
if $-Q\not\in q^{\mathbb Z}$ then 
the number of pairwise non--isomorphic simple $B$--modules 
is $e-1$,
\item[--]
if $-Q=q^f$, for $0\le f\le e-1$, then 
the number of pairwise non--isomorphic simple $B$--modules 
is either $e-f+1$ or $f+1$,
\item[--]
$B$ is a Brauer tree algebra with the tree being a 
straight line such that all indecomposable projective 
$B$--modules have radical length $3$ unless $f=0$ 
and either $e=3$ and $n=1,2$ or $e\ge 4$ and $n=1,2,3$. 
\end{itemize}
\end{itemize}
\end{lem}
\begin{proof}
(1) The statements follow from the computation of the decomposition numbers of 
$B$ 
by using the Jantzen--Schaper sum formula. See \cite[Exercise 5.10]{Ma} 
for example. 
We may prove these in another way: 
it is shown in \cite[Theorem 8.2]{Jo} that $B$ is Morita--equivalent 
to the non--semisimple block algebra of $\H_e^A(q)$ when the characteristic 
of $F$ is zero, and it is observed in \cite{EN} that this result is 
valid without the assumption on the characteristic. 
Hence, the results follow from the computation of the radical 
structure of indecomposable projective $\H_e^A(q)$--modules given in \cite{U}. 
See also \cite{Ge1} for a more general result 
in the characteristic zero case. 

\noindent
(2) When $e$ is odd, then they follow from 
Theorem \ref{Morita theorem} and (1). 
When $e$ is even, they were proved in the proof 
of \cite[Theorem 6.30]{AM2}. 
\end{proof}

The similar results hold for $\H_n^D(q)$ as follows. 

\begin{prop}
\label{PIM structure in type D}
Let $e\ge 3$ as before, and 
assume that $\H_n^D(q)$ is finite but not semisimple. Then, 
the radical series of an indecomposable 
projective $\H_n^D(q)$--module is one of the following form. 

\begin{equation*}
\begin{split}
&D^\lambda\\
&D^\lambda\\
&D^\lambda
\end{split}\quad
\text{where $h(\lambda)\ne\lambda$,}
\qquad
\begin{split}
&D^\lambda\\
&D^\mu\\
&D^\lambda
\end{split}\quad
\text{where $\begin{cases}h(\lambda)\ne\lambda,\;\mu,\\ h(\mu)\ne\mu.\\ 
\end{cases}$}
\end{equation*}

\begin{equation*}
\begin{split}
&D^\lambda\\
D^\mu_+&\oplus D^\mu_-\\
&D^\lambda
\end{split}\quad
\text{where $\begin{cases}h(\lambda)\ne\lambda,\\ h(\mu)=\mu.\end{cases}$}
\end{equation*}

\begin{equation*}
\begin{split}
&D^\lambda_+\\
&D^\mu_+\\
&D^\lambda_+
\end{split}\quad\text{and}\quad
\begin{split}
&D^\lambda_-\\
&D^\mu_-\\
&D^\lambda_-
\end{split}\quad\text{or}\quad
\begin{split}
&D^\lambda_+\\
&D^\mu_-\\
&D^\lambda_+
\end{split}\quad\text{and}\quad
\begin{split}
&D^\lambda_-\\
&D^\mu_+\\
&D^\lambda_-
\end{split}\quad
\text{where $\begin{cases} h(\lambda)=\lambda,\\ h(\mu)=\mu.
\end{cases}$}
\end{equation*}

\begin{equation*}
\begin{split}
&D^\lambda\\
D^\mu &\oplus D^\mu\\
&D^\lambda
\end{split}\quad
\text{where $\begin{cases} h(\lambda)\ne\lambda,\;\mu,\\ 
h(\mu)\ne\mu.
\end{cases}$}\quad
\begin{split}
&D^\lambda\\
D^\lambda &\oplus D^\nu\\
&D^\lambda
\end{split}\quad
\text{where $\begin{cases} h(\lambda)\ne\lambda,\;\nu,\\ h(\nu)\ne\nu.
\end{cases}$}
\end{equation*}

\begin{equation*}
\begin{split}
&D^\lambda\\
D^\mu &\oplus D^\nu\\
&D^\lambda
\end{split}\quad
\text{where $\begin{cases} h(\lambda)\ne\lambda,\;\mu,\;\nu,\\
h(\mu)\ne\mu,\;\nu\\
h(\nu)\ne\nu.
\end{cases}$}
\end{equation*}

\begin{equation*}
\begin{split}
&D^\lambda\\
D^\lambda \oplus &D^\nu_+ \oplus D^\nu_-\\
&D^\lambda
\end{split}\quad
\text{where $h(\lambda)\ne\lambda,\;
h(\nu)=\nu$.}
\end{equation*}

\begin{equation*}
\begin{split}
&D^\lambda\\
D^\mu \oplus &D^\nu_+ \oplus D^\nu_-\\
&D^\lambda
\end{split}\quad
\text{where $\begin{cases} h(\lambda)\ne\lambda,\;\mu\\
h(\mu)\ne\mu,\\
h(\nu)=\nu.\end{cases}$}
\end{equation*}

\begin{equation*}
\begin{split}
&D^\lambda\\
D^\mu_+ \oplus D^\mu_- &\oplus D^\nu_+ \oplus D^\nu_-\\
&D^\lambda
\end{split}\quad
\text{where $\begin{cases} h(\lambda)\ne\lambda,\\
h(\mu)=\mu,\;\;h(\nu)=\nu,\\
\mu\ne\nu.\end{cases}$}
\end{equation*}

\begin{equation*}
\begin{split}
&D^\lambda_+\\
&D^\mu\\
&D^\lambda_+
\end{split}\quad\text{and}\quad
\begin{split}
&D^\lambda_-\\
&D^\mu\\
&D^\lambda_-
\end{split}\quad
\text{where $\begin{cases}h(\lambda)=\lambda,\\ h(\mu)\ne\mu.\end{cases}$}
\end{equation*}

\begin{equation*}
\begin{split}
&D^\lambda_\pm\\
D^\mu_\pm&\oplus D^\nu_\pm\\
&D^\lambda_\pm
\end{split}\quad\text{or}\quad
\begin{split}
&D^\lambda_\pm\\
D^\mu_\pm&\oplus D^\nu_\mp\\
&D^\lambda_\pm
\end{split}\quad
\text{where $\begin{cases}h(\lambda)=\lambda,\; h(\mu)=\mu,\;
h(\nu)=\nu,\\
\text{$\lambda$, $\mu$, $\nu$ are pairwise distinct.}\end{cases}$}
\end{equation*}
Here, we write $D^\lambda$ instead of $D^\lambda\!\!\downarrow_{\H_n^D(q)}$ 
for short. 
\end{prop}
\begin{proof}
As $\H_n^D(q)$ is finite, $\H_n(q,1)$ is also finite. 
So, Lemma \ref{finite case results} implies that 
if $B$ is a block algebra of $\H_n(q,1)$ then the number of pairwise 
non--isomorphic $B$--modules is greater than or equal to $2$ and 
that if $\lambda$ is a Kleshchev bipartition then 
the radical structure of $P^\lambda$ is 
either 
\begin{equation*}
\begin{split}
&D^\lambda\\
&D^\mu\\
&D^\lambda
\end{split}\quad\text{or}\quad
\begin{split}
&D^\lambda\\
D^\mu &\oplus D^\nu\\
&D^\lambda
\end{split}
\end{equation*}
where $\lambda$, $\mu$ and $\nu$ are pairwise distinct. 

Firstly, we consider the case where the $\H_n(q,1)$--module 
$P^\lambda$ is uniserial. 
Assume that $h(\lambda)\ne\lambda$. 
Then, 
Lemma \ref{PIM in type D}(2) implies that 
$P^\lambda\!\!\downarrow_{\H_n^D(q)}=\overline P^\lambda$. Thus, 
$\overline P^\lambda$ has radical length $3$ and 
$H(\overline P^\lambda)=D^\mu\!\!\downarrow_{\H_n^D(q)}$. 
If $\mu=h(\lambda)$ then we are in the first case of the list above. 
If $h(\mu)\ne\mu$ and $\mu\ne h(\lambda)$ then we are in the second 
case of the list. If $h(\mu)=\mu$ then we are in the third case. 
Note that this case does not appear if the number of the isomorphism 
classes of simple $\H_n(q,1)$--modules in the block is $2$. 
Assume that $h(\lambda)=\lambda$. Then, $(P^\lambda_+)^\tau=P^\lambda_-$ 
implies that $h(\mu)\ne\mu$ cannot happen. 
Thus, $h(\mu)=\mu$ and we are in the fourth case or in the fifth case. 

Nextly, we consider the case where the $\H_n(q,1)$--module 
$P^\lambda$ is not uniserial. Then, the similar argument 
gives us the remaining cases. 
\end{proof}

A consequence of this proposition is the 
following. 

\begin{cor}
Assume that $e\ge3$ and $\H_n^D(q)$ is finite but not semisimple. 
Then, for simple $\H_n^D(q)$--modules $S$ and $T$, we have 
$\Ext_{\H_n^D(q)}^1(S,T)\ne 0$ if and only if 
$\Ext_{\H_n^D(q)}^1(T,S)\ne 0$. 
\end{cor}

\begin{prop}
\label{possible subquiver}
Let $B$ be a block algebra of $\H_n^A(q)$, $\H_n^B(q)$ or $\H_n^D(q)$, 
and suppose that $e\ge3$ and that 
$B$ is finite but not semisimple. Then one of the following 
holds. 
\begin{itemize}
\item[(a)]
The Gabriel quiver of $B$ contains a straight line of length $2$ 
(with $3$ nodes), or 
a directed graph with adjacency matrix 
$\binom{1\;1}{1\;0}$ or $\binom{0\;2}{2\;0}$. 
\item[(b)]
There are indecomposable projective $B$--modules of the following form. 
\begin{equation*}
\begin{split}
&S\\
&T\\
&S\end{split}\quad\text{and}\quad
\begin{split}
&T\\
&S\\
&T\end{split}
\end{equation*}
where $S$ and $T$ are non--isomorphic simple $B$--modules. 
\end{itemize}
\end{prop}
\begin{proof}
Assume that $B$ is a block algebra of $\H_n^A(q)$. Then 
Lemma \ref{finite case results}(1) implies that 
the number of pairwise non--isomorphic simple $B$--modules 
is equal to $e-1\ge 2$. Hence, we are in 
the case (a) if the number is greater than $2$ and 
in the case (b) if the number is equal to $2$, by 
Lemma \ref{finite case results}(1). 

Assume that $B$ is a block algebra of $\H_n^B(q)$. If $e$ is odd 
then Theorem \ref{Morita theorem} implies that 
$B$ is Morita--equivalent to a block algebra of $\H_k^A(q)$, 
for some $k$. If $e$ is even then 
Lemma \ref{finite case results}(2) implies that 
the number of pairwise non--isomorphic simple $B$--modules 
is equal to $\frac{e}{2}+2$ or $\frac{e}{2}\ge 2$. Hence, 
by the same argument as above, we are either in 
the case (a) or in the case (b) again. 

Assume that $B$ is a block algebra of $\H_n(q,1)$. 
Then, by Lemma \ref{finite case results}(2) again, 
the number of pairwise non--isomorphic simle $B$--modules 
is equal to $\frac{e}{2}+1\ge 3$. 
Thus, the Gabriel quiver of $B$ contains a straight line of 
length $2$. We denote the nodes of the line as follows. 

\setlength{\unitlength}{16pt}
\begin{picture}(15,2.3)(-3,0)
\put(0.3,0.8){$\lambda$}
\put(1,1){\line(1,0){2}}
\put(3.2,0.8){$\mu$}
\put(4,1){\line(1,0){2}}
\put(6.2,0.8){$\nu$}
\end{picture}

Assume that $h(\lambda)\ne\lambda$, $h(\mu)\ne\mu$ and 
$h(\nu)\ne\nu$. If $h(\lambda)\ne\mu, \nu$ and 
$h(\mu)\ne \nu$ then Lemma \ref{PIM in type D} 
and Proposition \ref{PIM structure in type D} imply that 
the Gabriel quiver of $\H_n^D(q)$ also contains a straight 
line of length $2$. If $h(\lambda)=\mu$ or $h(\lambda)=\nu$ 
then we obtain one of the other two quivers listed in the case (a). 

Assume that $h(\lambda)=\lambda$ and $h(\mu)\ne\mu$ or 
$h(\lambda)\ne\lambda$ and $h(\mu)=\mu$. 
Then the Gabriel quiver of $\H_n^D(q)$ 
contains the straight line whose 
nodes are $\lambda_+, \mu, \lambda_-$ or 
$\mu_+, \lambda, \mu_-$. 

Assume that $h(\lambda)=\lambda$, $h(\mu)=\mu$ and 
$h(\nu)\ne\nu$. Then the Gabriel quiver of $\H_n^D(q)$ 
contains the straight line whose 
nodes are $\mu_+$, $\nu$ and either $\lambda_+$ or $\lambda_-$. 

Assume that $h(\lambda)=\lambda$, $h(\mu)=\mu$ and 
$h(\nu)=\nu$. Then, the similar argument shows that 
the Gabriel quiver of $\H_n^D(q)$ 
contains a straight line of length $2$. 
\end{proof}

\begin{lem}
\label{finite tensor finite}
Assume that $A$ and $B$ are block algebras of 
$\H_n^A(q)$, $\H_n^B(q)$, $\H_n^D(q)$ or $\H_n(q,Q)$ and 
that they are finite but not semisimple. 
\begin{itemize}
\item[(1)]
Suppose that $e\ge 3$. Then $A\otimes B$ is wild. 
\item[(2)]
Suppose that $e=2$. Then $A\otimes B$ is tame. 
\end{itemize}
\end{lem}
\begin{proof}
(1) Assume that the Gabriel quiver of $A$ contains one of 
the directed graphs listed in the case (a) of 
Proposition \ref{possible subquiver}. As 
the Gabriel quiver of $B$ always contains 
a straight line of length $1$ by 
Proposition \ref{possible subquiver}, 
the Gabriel quiver of $A\otimes B$ contains 
the product of these quivers as a subquiver. 
In any of these three cases, 
Lemma \ref{criterion for wildness-6} 
implies that $A\otimes B$ is wild. 

Now assume that both $A$ and $B$ are as in the case (b) of 
Proposition \ref{possible subquiver}. Let $S$ and $T$ be 
the simple $A$--modules. We denote the projective 
covers of $S$ and $T$ by $P(S)$, $P(T)$ respectively. 
Similarly, we denote the simple $B$--modules by $S'$ and $T'$, 
the projective covers $P(S')$ and $P(T')$. We 
define indecomposable projective $A\otimes B$--modules $P_1$ and $P_2$ by 
\[
P_1=P(S)\otimes P(S'),\;\;
P_2=P(T)\otimes P(T').
\]
We shall show that $\End_{A\otimes B}(P_1\oplus P_2)$ 
is wild. This implies the result because its opposite 
algebra may be identified with an algebra of 
the form $p(A\otimes B)p$ where $p$ is an idempotent of 
$A\otimes B$. 
Define a factor algebra $R$ of 
$\End_{A\otimes B}(P_1\oplus P_2)$ by 
\[
R=\End_{A\otimes B}(P_1/\Rad^3 P_1\oplus P_2/\Rad P_2).
\]
Note that $P_1/\Rad^3 P_1$ has the 
radical series of the following form. 
\begin{equation*}
\begin{split}
&S\otimes S'\\
S\otimes T'&\;\;\,\oplus\;\; S\otimes T'\\
S\otimes S'\;\;\oplus\;\; &T\otimes T'\;\;\oplus\;\; S\otimes S'
\end{split}
\end{equation*}
Hence, $R$ has the form $FQ/I$ where $Q_0=\{1,2\}$ and $Q_1$ consists of 
two loops on the node $1$ and the arrow $1\leftarrow 2$. This implies 
that $R$ is wild by Lemma \ref{criterion for wildness-5}. 

\noindent
(2) By Lemma \ref{finite case e=2}, we may assume that 
both $A$ and $B$ are isomorphic to $F[X]/(X^2)$. 
Thus, $A\otimes B\simeq F[X, Y]/(X^2, Y^2)$, which is tame. 
\end{proof}

Let $(W,S)$ be a finite Weyl group of classical type, and let 
$P_W(x)$ be its Poincar\'e polynomial. Fix $q\in F^\times$ and 
denote the Hecke algebra associated with $(W,S)$ and $q$ by $\H_W(q)$. 
We assume that $q$ is a primitive $e^{th}$ root of unity with $e\ge2$. 
Note that if $q=1$ then $\H_W(q)=FW$ and we can also tell the representation 
type. 

\begin{thm}
\label{final result}
Let $P_W(x)$, $\H_W(q)$ and $e\ge 2$ be as above. Then, 
\begin{itemize}
\item[(1)]
Assume that $e\ge 3$. Then $\H_W(q)$ is 
\begin{itemize}
\item[--]
finite if $(x-q)^2$ does not divide $P_W(x)$.
\item[--]
wild otherwise.
\end{itemize}
\item[(2)]
Assume that $e=2$. Then, $\H_W(q)$ is 
\begin{itemize}
\item[--]
finite if $(x-q)^2$ does not divide $P_W(x)$.
\item[--]
tame if $(x-q)^2$ divides but $(x-q)^3$ does not divide $P_W(x)$.
\item[--]
wild otherwise.
\end{itemize}
\end{itemize}
\end{thm}
\begin{proof}
We write $W=W_1\times\cdots\times W_s$ where $W_i$ are irreducible 
Weyl groups. 

\noindent
(1) If $(x-q)^2$ does not divide $P_W(x)$, then 
$x-q$ divides at most one $P_{W_i}(x)$ and $(x-q)^2$ does 
not divide $P_{W_i}(x)$. Thus, $\H_{W_i}(q)$ is finite 
by Theorem \ref{one parameter case}, and 
all the other $\H_{W_j}(q)$, for $j\ne i$, are semisimple. 
Hence, the result follows. If 
$x-q$ divides $P_{W_i}(x)$ and $P_{W_j}(x)$, for $i\ne j$, 
then $\H_{W_i}(q)$ and $\H_{W_j}(q)$ are finite but not 
semisimple. Thus, Lemma \ref{finite tensor finite}(1) implies 
the result. The result 
for the case where $(x-q)^2$ divides some $P_{W_i}(x)$ was 
proven in Theorem \ref{one parameter case}. 

\noindent
(2) We have already proved that if $\H_n(q,Q)$ is tame then there is 
an $\H_n(q,Q)$--module with complexity $2$. 
Assume that $\H_n^A(q)$ is tame. 
Thus, $e=2$ and $n=4$ or $n=5$. 
\cite[Proposition A]{EN} asserts that $\H_4^A(q)$ is 
Morita--equivalent to the path algebra which is considered 
in (case 4b) of the proof of Theorem \ref{two parameter case}(2). 
Thus, there is an $\H_4^A(q)$--module with comlexity $2$. 
Similarly, \cite[Proposition B]{EN} asserts that the unique 
non--semisimple block algebra of $\H_5^A(q)$ has the directed graph $Q$ 
with adjacency matrix $\binom{1\;1}{1\;1}$ as its Gabriel 
quiver, and if we denote the loops on the nodes $1$ and $2$ by 
$\alpha$, $\beta$, and $1\rightarrow 2$, $1\leftarrow 2$ by 
$\mu$ and $\nu$, then the relations are given by
\[
\alpha^2=0,\;\;\beta^2=0,\;\;
\mu\alpha=0,\;\;\alpha\nu=0,\;\;\beta\mu=0,\;\;\nu\beta=0.
\]
This is a special biserial algebra which satisfies the 
assumptions of Lemma \ref{complexity=2(2)}. Thus, 
there is an $\H_5^A(q)$--module with complexity $2$. 
$\H_n^D(q)$ cannot be tame. 

Hence, if $\H_n^X(q)$, for some $n$ and $X$, is tame 
then there is an $\H_n^X(q)$--module with complexity $2$. 

Assume that $(x-q)^2$ does not divide $P_W(x)$. Then, 
$\H_{W_i}(q)$ is semisimple for all but at most one $i$, say $i_0$. 
Then, $\H_{W_{i_0}}(q)$ is finite by Theorem \ref{one parameter case} and 
so is $\H_W(q)$. Assume that $(x-q)^2$ divides $P_W(x)$ but $(x-q)^3$ does not 
divide $P_W(x)$. If $(x-q)^2$ divides $P_{W_i}(x)$, for some $i$, 
then all the other $\H_{W_j}(q)$, $j\ne i$, are semisimple and 
the result follows from Theorem \ref{one parameter case}. If 
$x-q$ divide $P_{W_i}(x)$ and $P_{W_j}(x)$ for distinct $i$ and $j$, 
then $\H_{W_i}(q)$ and $\H_{W_j}(q)$ are finite and not semisimple. 
All the other $\H_{W_k}(q)$, $k\not\in\{i,j\}$, are semisimple. 
Thus, Lemma \ref{finite tensor finite}(2) implies that 
$\H_W(q)$ is tame. 
Assume that $(x-q)^3$ divides $P_W(x)$. If there are distinct 
$i,j,k$ such that $x-q$ divides $P_{W_i}(x)$, $P_{W_j}(x)$ and 
$P_{W_k}(x)$ then $\H_{W_i}(q)$, $\H_{W_j}(q)$ and $\H_{W_k}(q)$ 
are finite but not semisimple. Thus, there are modules $M_i$, $M_j$ and 
$M_k$ of these algebras such that they have complexity greater than or 
equal to $1$. 
As $M_i\otimes M_j\otimes M_k$ has complexity greater than or 
equal to $3$, 
the result follows from Theorem \ref{complexity}(3). 
If $(x-q)^2$ divides $P_{W_i}(x)$ but $(x-q)^3$ does not divide 
$P_{W_i}(x)$ and $x-q$ divides $P_{W_j}(x)$, for 
$j\ne i$, then, the same complexity argument works: $\H_{W_i}(q)$ is 
tame by Theorem \ref{one parameter case}(2) and there is an 
$\H_{W_i}(q)$--module with complexity $2$, as is remarked above. 
The case where $(x-q)^3$ divides 
$P_{W_i}(x)$, for some $i$, reduces to 
Theorem \ref{one parameter case} as before. 
\end{proof}

\section{Appendix}
The following is a list of corrections and remarks for 
\cite{A1}. The author is 
grateful to Professor Kashiwara for pointing out many of these. 

\begin{itemize}
\item[(p.11)]
In Assertion 4. If we consider $U\rightarrow U(\mathfrak g)$ and use 
the universal property of $U(\mathfrak g)$, the proof would be much 
shorter. 
\item[(p.20)]
As $\tilde M(\lambda)$ is a $T(V)$--module, 
it is also a ${\mathfrak g}(V)$--module. Similarly, 
the algebra homomorphism $T(V_-)\rightarrow U(\tilde{\mathfrak n}_-)$ 
induced by $V_-\subset \tilde{\mathfrak n}_-$ 
gives the surjection ${\mathfrak g}(V_-)\rightarrow \tilde{\mathfrak n}_-$. 
\item[(p.25)]
In the proof of Proposition 4.5. After \lq\lq We can choose $m\ne0$\rq\rq, add 
\lq\lq in $M_0\otimes_K\overline K$\rq\rq. After \lq\lq $c=\pm v^l$.\rq\rq, 
add \lq\lq Hence, we can choose $m$ in $M_0$.\rq\rq. 
\item[(p.25)]
In Definition 4.6. $vt-v^{-1}t^{-1}\rightarrow vt+v^{-1}t^{-1}$ and 
$\pm\frac{[l+1]}{v-v^{-1}}\rightarrow 
\pm\frac{v^{l+1}+v^{-l-1}}{(v-v^{-1})^2}$.
\item[(p.89)]
in l.-10; $R_i(\lambda)
\setminus\{y_1,\dots,y_s,x_1,\dots,x_t\} \rightarrow 
R_i(\lambda)\setminus\{x_1,\dots,x_t\}$.
\item[(p.92)]
in l.5; $\lambda(\sigma,k,S)\setminus\{x, \tilde x\} \rightarrow 
\{\lambda(\sigma,k,S)\setminus\{x,\tilde x\}\,|\,S \subset N_iR_\lambda\}$.
\item[(p.98)]
Just above (12.1); $g_{\mu,\mu}(v)=1$ $\rightarrow$ 
$g_{\mu,\mu}(v)=1,\;g_{\lambda,\mu}(0)=\delta_{\lambda,\mu}$.
\item[(p.98)]
in l.-5; $\mathbb Q[v,v^{-1}] \rightarrow 
\mathbb Q[v]\oplus v^{-1}\mathbb Q[v^{-1}]$ must be 
$\mathbb Q[v,v^{-1}] \rightarrow v\mathbb Q[v]\oplus \mathbb Q[v^{-1}]$.
\item[(p.108)]
the denominator of (1.3); 
$v_{c_i}-q^{k_i}v_{c_{i-1}}\rightarrow q^{k_i}v_{c_i}-v_{c_{i-1}}$ 
and 
the $(1,1)$--entry; $(q-1)v_{c_i}\rightarrow q^{k_i}(q-1)v_{c_i}$, 
the $(2,2)$--entry; $(1-q)v_{c_i}\rightarrow (1-q)v_{c_{i-1}}$. 
\item[(p.108)]
In the definition of the $\lambda$-separatedness; 
$v_{c_i}-q^{k_i}v_{c_{i-1}} \rightarrow q^{k_i}v_{c_i}-v_{c_{i-1}}$. 
\item[(p.109)]
in l.14; Delete the sentence \lq\lq To prove this, ...\rq\rq, and 
change \lq\lq If we consider ... requirement.\rq\rq to 
\lq\lq Consider $w\in S_n$ such that $w{\bf t}={\bf s}$. 
We argue by induction on $l(w)$. Choose $s_i=(i-1,i)$ with 
$s_iw<w$ so that ${\bf s}'=s_i{\bf s}$ is standard. 
By the induction hypothesis, ${\bf s}'\in W$. Then $a_i{\bf s}'\in W$ 
because ${\bf s}$ appears in $a_i{\bf s}'$ with a nonzero coefficient.\rq\rq 

Note that if there is no such $s_i$ then the entry, say $j$, in ${\bf s}$ 
of the node $n$ of ${\bf t}$ must be $n$: otherwise we could choose 
$s_{j+1}$. We delete $n$ and continue the same argument to conclude that 
${\bf s}={\bf t}$. 
\item[(p.110)]
in l.5; $v_{c_i}-q^{k_i}v_{c_{i-1}} \rightarrow q^{k_i}v_{c_i}-v_{c_{i-1}}$. 
\item[(p.110)]
In Proposition 13.10(2); 
Add \lq\lq and $V^\lambda_{\mathbb K}$ is irreducible.\rq\rq
\item[(p.112)]
Let $L''=\pi^{N-1}L'+L$. Then, we have 
$\pi L''\subset L\subset L''$ and $\pi^{N-1}L'\subset L''\subset L'$. 
Hence, we can also prove the result by induction on $N$. 
\item[(p.113)]
The statement of Lemma 13.19 is not accurate. We must consider 
the subcategory of $\hat H_{n,{\mathbb K}}-mod$ 
whose objects are those modules which admits an 
$\mathbb S$--lattice here. In practice, we work with the category 
$\mathcal C$ (Corollary 13.26). As any module in 
$\mathcal C$ admits an $\mathbb S$--lattice, 
this inaccuracy does not affect the rest of the chapters. 
\item[(p.115)]
Murphy's lemma, which says that ${\bf s}\trianglerighteq{\bf t}$ if and 
only if $d({\bf s})\le d({\bf t})$, is used to show that 
$m_{{\bf s}{\bf t}^\lambda}{\bf t}^\lambda\in V^\lambda_{\mathbb K}$ 
belongs to $\oplus_{{\bf u}\trianglerighteq{\bf s}}\mathbb K{\bf u}$. 
The lemma easily follows from Ehresmann's characterization of 
the Bruhat order, which says that $w_1\le w_2$ if and only if 
$|\{i\,|\,i\le p, w_1(i)\le q\}|\ge |\{i\,|\,i\le p, w_2(i)\le q\}|$ 
for all $p$ and $q$. 
\item[(p.115)]
In Corollary 13.24; Replace $\hat H_{n,{\mathbb K}}-mod$ with 
the subcategory of modules which admits an 
$\mathbb S$--lattice.
\item[(p.121)]
in l.5; the Grothendieck group $V(\mathbb F)$ $\rightarrow$ 
the Grothendieck groups $K_0(\mathcal H_{n,\mathbb F}-mod)$
\item[(p.129)]
At the bottom of Definition 14.15; Add 
\lq\lq where we fix a square root $\sqrt{p}$ of $p$ and 
we choose the square root $p^{e/2}$ of $p^e$ 
to be $-\sqrt{p}^e$.\rq\rq. 
\item[(p.132)]
Proof of Lemma 14.25 is absurd. To save this, 
we let $C\subset \overline{\mathbb Q}_l^\times$ be the 
multiplicative group of pure numbers, 
$\mathbb ZC$ its group algebra, and change the 
statement of Lemma 14.25 as follows. 
Then the same proof (but we do not assert 
the \lq\lq(a root of unity)\rq\rq part) 
works. 

\noindent
\lq\lq Then there exists an element 
$F^{\underline m''}_{\underline m,\underline m'}\in 
\mathbb ZC$ such that if we write 
$F^{\underline m''}_{\underline m,\underline m'}=
\sum n_i[c_i]$, where $n_i\in\mathbb Z$ and $c_i\in C$, 
then the following holds for all $e$: 
$F^{\underline m''}_{\underline m,\underline m'}(\mathbb F_{p^e})
=\sum n_ic_i^e$.\rq\rq
\item[(p.133)]
In Definition 14.26; Change 
\lq\lq Let $H_A$ be\rq\rq to 
\lq\lq Let $v=-[\sqrt{p}^{-1}]$, $A=\mathbb Z[v,v^{-1}]$ 
and define $H_A$ to be\rq\rq, and change 
\lq\lq Define an $A$-bilinear map\rq\rq to 
\lq\lq Then let $H_C=H_A\otimes_A \mathbb ZC$ 
and define a $\mathbb ZC$-bilinear map\rq\rq, 
and replace 
\lq\lq $F^{\underline m''}_{\underline m,\underline m'}(v)$\rq\rq 
with 
\lq\lq $F^{\underline m''}_{\underline m,\underline m'}$\rq\rq. 

\noindent
\lq\lq $H_A$ is an $A$-algebra.\rq\rq $\rightarrow$ 
\lq\lq $H_C$ is a $\mathbb ZC$-algebra.\rq\rq 

\noindent
At the end of the definition, add 
\lq\lq We call $H_A$ {\bf the Hall module}.\rq\rq
\item[(p.133)]
In Lemma 14.27; Change \lq\lq $H_A$\rq\rq 
to \lq\lq $H_C$\rq\rq in (1) and 
add \lq\lq In particular, 
$H_A$ is a right $U_A^-$-module.\rq\rq to the end 
of (2). 

As a $\mathbb ZC$-module, $H_C$ is nothing but 
the Grothendieck group of equivariant mixed Weil sheaves. 
By Lemma 14.27, $H_A$ is a $U_A^-$-module. 

Note that we only need $H_A$ at $v=1$ to 
prove Theorem 12.5, and that the statements 
Lemma 14.27(2) and Proposition 14.28 for 
$H_A$ remain unchanged. 
\item[(p.137)]
Lemma 14.34; \lq\lq ... makes $\mathcal H$ into 
a unital associative algebra, 
which is isomorphic to $H_A$.\rq\rq 
$\rightarrow$ 
\lq\lq ... restricts to an isomorphism between 
the $A$-subalgebras of $\mathcal H$ and $H_A$  
generated by $f_i^{(n)}$'s. More precisely, 
each object $F^.$ which appears in $U_A^-\subset\mathcal H$ is 
given a unique mixed structure such that 
the identification of $[F^.]\in\mathcal H$ 
with $\sum_{j\in\mathbb Z}(-1)^j[\mathcal H^j(F^.)]
\in H_C$ gives the identification of $U_A^-\subset\mathcal H$ 
with $U_A^-\subset H_A$.\rq\rq
\item[(p.137)]
in l.-8 and l.-7; $\mathcal H^i(R{p_3}_!C^.) \rightarrow 
\hphantom{}^p\mathcal H^i(R{p_3}_!C^.)$. 
\item[(p.137)]
in l.-4; \lq\lq take ${\operatorname{Fr}^e}^*(F^.)\simeq F^.$\rq\rq 
$\rightarrow$ 
\lq\lq take ${\operatorname{Fr}^e}^*(F^.)\simeq F^.$ as in 
[cb-E, Theorem 5.2]\rq\rq
\item[(p.138)]
At the end of the proof of Lemma 14.34; 
\lq\lq To summarize, ...\rq\rq 
$\rightarrow$ 
\lq\lq Thus, if we expand the product of canonical 
basis elements $b_1b_2$ into a linear combination 
of canonical basis elements $\sum c_{b_1b_2}^{b_3}(v)b_3$ in $\mathcal H$, 
[cb-E, Theorem 5.4] shows that the polynomial 
$c_{b_1b_2}^{b_3}(v)$ in the shift $v$ corresponds to 
the polynomial $c_{b_1b_2}^{b_3}(v)$ in $v=-[\sqrt{p}^{-1}]$ under 
the identification. Thus we get the isomorphism of 
the subalgebras $U_A^-\subset\mathcal H$ and 
$U_A^-\subset H_A$.\rq\rq
\item[(p.139)]
$R{p_3}_!\mathbb C\simeq \oplus_{i\in\mathbb Z}\mathcal H^i[-i] \rightarrow 
[R{p_3}_!\mathbb C]=\sum_{i\in\mathbb Z}
[\mathcal H^i(R{p_3}_!\mathbb C)[-i]]$. 
\item[(p.141)]
After Theorem 14.41; Add 
\lq\lq The second part follows from Lemma 14.27(2) and 
[CG, Theorem 8.6.23]. Recall that 
canonical basis elements are given 
mixed structure as in Lemma 14.34.\rq\rq
\item[(p.142)]
\lq\lq{\bf specialized Hall algebra}\rq\rq $\rightarrow$ 
\lq\lq{\bf specialized Hall module}\rq\rq 
\end{itemize}

\end{document}